\documentclass[oneside, 11pt]{article}
\usepackage{setspace}
\usepackage{titlesec}
\usepackage[utf8]{inputenc}
\usepackage{tabularx}	
\usepackage{etoolbox}
\usepackage{mathtools}
\usepackage{amsmath}
\usepackage{amssymb} 
\usepackage{amsfonts}
\usepackage{mathrsfs} 
\usepackage{amsthm,pifont}
\usepackage{graphicx}
\usepackage{bbm}
\usepackage[colorlinks=true, pdfstartview=FitV, linkcolor=blue,
citecolor=blue, urlcolor=blue, pagebackref=false]{hyperref}
\usepackage{cleveref}
\usepackage[nottoc,numbib]{tocbibind}
\usepackage{xcolor}
\usepackage{tikz}
\usepackage{cite}
\usepackage{braket}
\usepackage{leftidx,tensor}
\usepackage{dsfont}
\usepackage{calc}
\usepackage{caption}
\usepackage{subcaption}
\usepackage{color}
\usepackage[printonlyused,nohyperlinks]{acronym}
\usetikzlibrary{decorations.markings}
\usetikzlibrary{shapes.geometric}
\usetikzlibrary{patterns}
\tikzstyle{vertex}=[circle, draw, inner sep=0pt, minimum size=6pt]
\newcommand{\vertex}{\node[vertex]}
\definecolor{superlightgray}{RGB}{235,235,235}
\definecolor{wellgray}{RGB}{170,170,170}
\definecolor{lightergray}{RGB}{225,225,225}

\usepackage[utf8]{inputenc}
\usepackage[english]{babel}
\usepackage{enumitem}

\allowdisplaybreaks

\newtheorem{theorem}{Theorem}[section]
\makeatletter
\patchcmd{\ttlh@hang}{\parindent\z@}{\parindent\z@\leavevmode}{}{}
\patchcmd{\ttlh@hang}{\noindent}{}{}{}
\makeatother
\titleformat*{\section}{\large\bfseries}
\titleformat*{\subsection}{\small\bfseries}
\titleformat*{\subsubsection}{\small\bfseries}
\titleformat*{\paragraph}{\small\bfseries}
\titleformat*{\subparagraph}{\small\bfseries}

\newcommand{\N}{\mathbb{N}}
\newcommand{\R}{\mathbb{R}}
\newcommand{\Z}{\mathbb{Z}}
\newcommand{\E}{\mathbb{E}}
\newcommand{\p}{\mathbb{P}}
\newcommand{\md}{\ensuremath{\mathrm{d}}}
\newcommand{\dxp}{\theta}
\newcommand{\e}{\mathcal{E}}
\newcommand{\eps}{\varepsilon}
\newcommand{\dia}{\text{Diam}}

\newcommand{\mo}{\mathbf{1}}
\newcommand{\mz}{\mathbf{0}}
\newcommand{\cA}{\mathcal{A}}
\newcommand{\cB}{\mathcal{B}}
\newcommand{\cC}{\mathcal{C}}

\newcommand{\cG}{\mathcal{G}}
\newcommand{\cH}{\mathcal{H}}
\newcommand{\cN}{\mathcal{N}}
\newcommand{\cP}{\mathcal{P}}
\newcommand{\cS}{\mathcal{S}}
\newcommand{\cX}{\mathcal{X}}

\newtheorem{lemma}[theorem]{Lemma}
\newtheorem{remark}[theorem]{Remark}
\newtheorem{proposition}[theorem]{Proposition}
\newtheorem{definition}[theorem]{Definition}
\newtheorem{corollary}[theorem]{Corollary}

\newtheorem{example}[theorem]{Example}

\usepackage{geometry}
\begin{document}
	
\title{Distances in $\frac{1}{|x-y|^{2d}}$ percolation models for all dimensions}

\author{Johannes B\"aumler\footnote{ \textsc{Department of Mathematics, TU Munich, Germany}. E-Mail: \href{mailto:johannes.baeumler@tum.de}{johannes.baeumler@tum.de}} }

\maketitle
	
	\begin{center}
		\parbox{13cm}{ \textbf{Abstract.} We study independent long-range percolation on $\mathbb{Z}^d$ for all dimensions $d$, where the vertices $u$ and $v$ are connected with probability 1 for $\|u-v\|_\infty=1$ and with probability $p(\beta,\{u,v\})=1-e^{-\beta \int_{u+\left[0,1\right)^d} \int_{v+\left[0,1\right)^d} \frac{1}{\|x-y\|_2^{2d}}d x d y } \approx \frac{\beta}{\|u-v\|_2^{2d}}$ for  $\|u-v\|_\infty \geq 2$. Let $u \in \mathbb{Z}^d$ be a point with $\|u\|_\infty=n$. We show that both the graph distance $D(\mathbf{0},u)$ between the origin $\mathbf{0}$ and $u$ and the diameter of the box $\{0 ,\ldots, n\}^d$ grow like $n^{\theta(\beta)}$, where $0<\theta(\beta ) < 1$. We also show that the graph distance and the diameter of boxes have the same asymptotic growth when two vertices $u,v$ with $\|u-v\|_2 > 1$ are connected with a probability that is close enough to $p(\beta,\{u,v\})$.
		Furthermore, we determine the asymptotic behavior of $\theta(\beta)$ for large $\beta$, and we discuss the tail behavior of $\frac{D(\mathbf{0},u)}{\|u\|_2^{\theta(\beta)}}$.}
	\end{center}

\let\thefootnote\relax\footnotetext{{\sl MSC Class}: 05C12, 60K35, 82B27, 82B43}
\let\thefootnote\relax

\hypersetup{linkcolor=black}
\tableofcontents
\hypersetup{linkcolor=blue}

\section{Introduction}

Consider independent long-range bond percolation on $\Z^d$ where all edges $\{u,v\}$ with $\|u-v\|_\infty=1$ are open and an edge $\{u,v\}$ with $\|u-v\|_\infty \geq 2$ is open with probability
\begin{equation*}
	p(\beta,\{u,v\}) \coloneqq 1-e^{-\beta \int_{u+\cC} \int_{v+\cC} \frac{1}{\|x-y\|^{2d}} \md x \md y } ,
\end{equation*}
where $\cC\coloneqq \left[0,1\right)^d$, $\beta \geq 0$, and where we use $\|\cdot\|$ for the Euclidean norm of a vector. We call the corresponding probability measure $\p_\beta$ and denote its expectation by $\E_\beta$. The resulting graph is clearly connected and the graph distance $D(u,v)$ between two points $u,v\in \Z^d$ satisfies $D(u,v) \leq \|u-v \|_\infty$. We are interested in the scaling of the typical  distance of two points $u,v \in \Z^d$ and the scaling of the diameter of boxes $\left\{0,\ldots,N\right\}^d$. In \cite{coppersmith2002diameter} it is proven that the typical diameter of some box grows at most polynomially with some power strictly smaller than $1$. More precisely, Coppersmith, Gamarnik, and Sviridenko proved that for all $\beta >0$ there exists an exponent $\dxp^\prime = \dxp^\prime(\beta) < 1 $ such that $\lim_{N \to \infty} \p_\beta \left( \dia \left(\left\{0,\ldots,N\right\}^d\right) \leq N^{\dxp^\prime}\right) = 1$. However, the authors do not give any polynomial lower bound for dimensions $d\geq 2$. An analogous lower bound was already conjectured in \cite{benjamini2001diameter, biskup2004scaling}, and an exact lower bound was later proven to hold in one dimension: In \cite{ding2013distances} Ding and Sly showed that for the connection probability $p(\beta, \{u,v\})$ given by $p(\beta, \{u,v\})= \frac{\beta}{|u-v|^2}\wedge 1$ for $|u-v|\geq 2$ and $p(\beta, \{u,v\})=  1$ for $|u-v|= 1$ the typical distance between the two points $0,n \in \N$ and the diameter of $\{0, \ldots ,n\}$ both grow like $n^\dxp$ for some $\dxp \in \left(0,1\right]$, where $\dxp=1$ if and only if $\beta=0$. More precisely, they proved that
\begin{align*}
n^\dxp \approx_P D(0,n) \approx_P \dia \left(\{0, \ldots ,n\}\right) \approx_P \E\left[D(0,n)\right]
\end{align*}
where the notation $A(n)\approx_P B(n)$ means that for all $\eps >0 $ there exist $0<c<C<\infty$ such that $\p\left(cB(n) \leq A(n) \leq CB(n)\right)> 1-\eps$ for all $n \in \N$. In this paper, we prove an analogous result for all dimensions. 

\subsection{Main results}

\begin{theorem}\label{theo:exponent}
	For all dimensions $d$ and all $\beta > 0$, there exists an exponent $\dxp = \dxp(d,\beta) \in (0,1)$ such that
	\begin{align}\label{eq:theo:exponent1}
		\|u\|^\dxp \approx_P D\left(\mz,u\right)  \approx_P \E_\beta\left[D(\mz,u)\right]
	\end{align}
	and
	\begin{align}\label{eq:theo:exponent2}
	k^\dxp \approx_P \dia\left( \left\{0,\ldots, k\right\}^d \right)  \approx_P \E_\beta\left[ \dia\left( \left\{0,\ldots, k\right\}^d \right) \right]\text.
	\end{align}
\end{theorem}
As we consider the dimension $d$ as fixed, we also write $\theta(\beta)$ for $\theta(d,\beta)$, although $\theta(d,\beta)$ depends on the dimension.
We write $\mz$ for the vector with all entries equal to $0$ and the notation $A(u)\approx_P B(u)$ means that for all $\eps >0 $ there exist $0<c<C<\infty$ such that $\p_\beta \left(cB(u) \leq A(u) \leq CB(u)\right)> 1-\eps$ for all $u \in \Z^d$. 
The inclusion probability $p(\beta,\{u,v\}) \coloneqq 1-e^{-\beta \int_{u+\cC} \int_{v+\cC} \frac{1}{\|x-y\|^{2d}} \md x \md y }$ is only one possible choice for a function which asymptotically grows like $\frac{\beta}{\|u-v\|^{2d}}$. In Section \ref{section:comparison}, we will show the same results for other possible choices of such functions. Examples of inclusion probabilities we  consider are $\frac{\beta}{\|u-v\|^{2d}}\wedge 1$ and $1-e^{-\tfrac{\beta}{\|u-v\|^{2d}}}$.\\

The exponent $\dxp=\dxp(\beta)$ defined in \Cref{theo:exponent} arises through a subadditivity argument (see \cref{sec:submultiplicativity} below) and its precise value is not known to us. However, we determine the asymptotic behavior of the function $\dxp(\beta)$ for large $\beta$.
\begin{theorem}\label{theo:largebeta}
	For all dimensions $d$, there exist constants $0<c<C<\infty$ such that
	\begin{align*}
	\frac{c}{\log(\beta)} \leq \dxp(\beta) \leq \frac{C}{\log(\beta)}
	\end{align*}
	for all $\beta \geq 2$.
\end{theorem}

\subsection{The continuous model}\label{subsec:cts}

For $\beta>0$, the described discrete percolation model has a self-similarity that comes from a coupling with the underlying continuous model, that we will now describe for any dimension. This will also explain our, at first sight complicated, choice of the connection probability. Consider a Poisson point process $\tilde{\e}$ on $\R^d \times \R^d$ with intensity $\frac{\beta}{2 \|t-s\|_2^{2d}}$. Define the symmetrized version $\e$ by $\e \coloneqq \{(t,s) \in \R^d \times \R^d : (s,t) \in \tilde{\e}\} \cup \tilde{\e}$. For $u,v \in \Z^d$ with $\|u-v\|_\infty \geq 1$ we put an edge between $u$ and $v$ if and only if $\left(\left(u+\cC\right) \times \left(v+\cC\right)\right) \cap \e \neq \emptyset$, where we use the notation $\cC = \left[0,1\right)^d$. The cardinality of $\left(\left(u+\cC\right) \times \left(v+\cC\right)\right) \cap \tilde{\e}$ is a random variable with Poisson distribution and parameter $\int_{u+\cC} \int_{v+\cC} \frac{\beta}{2 \|t-s\|^{2d}} \md s \md t$. Thus, by the properties of Poisson processes, the probability that $u \nsim v$ equals
\begin{align*}
	&\p\left(\left(\left(u+\cC\right) \times \left(v+\cC\right)\right) \cap \e = \emptyset\right) = \p\left(\left(\left(u+\cC\right) \times \left(v+\cC\right)\right) \cap \tilde{\e} = \emptyset\right)^2\\
	& = \left(e^{-\int_{u+\cC} \int_{v+\cC} \frac{\beta}{2 \|t-s\|^{2d}} \md s \md t}\right)^2
	= e^{-\int_{u+\cC} \int_{v+\cC} \frac{\beta}{ \|t-s\|^{2d}} \md s \md t}
	= 1 - p(\beta, \{u,v\})
\end{align*}
which is exactly the probability that $u \nsim v$ under the measure $\p_\beta$.
Note that for $\{u,v\}$ with $\|u-v\|_\infty=1$ we have $\int_{u+\cC} \int_{v+\cC} \frac{\beta}{ \|t-s\|^{2d}} \md s \md t = \infty$. So we really get that all edges of the form $\{u,v\}$ with $\|u-v\|_\infty=1$ are open.
The construction with the Poisson process also implies that the presence of different bonds is independent and thus the resulting measure of the random graph constructed above equals $\p_\beta$.
The chosen inclusion probabilities have many advantages. First of all, the resulting model is invariant under translation and invariant under the reflection of coordinates, i.e.,  when we change the $i$-th component $p_i(x)$ of all $x\in \Z^d$ to $-p_i(x)$. Furthermore, the model has the following self-similarity: For some vector $u= (p_1(u),\ldots,p_d(u))\in \Z^d$ and $n\in \N_{>0}$ we define the translated boxes  $V_u^n \coloneqq \prod_{i=1}^{d} \{p_i(u) n , \ldots, (p_i(u)+1)n-1\} = nu  + \prod_{i=1}^{d} \{0,\ldots, n-1\} $. Then for all points $u,v \in \Z^d$, and all $n\in \N_{>0}$ one has
\begin{align*}
	&\p_\beta\left(V_u^n \nsim V_v^n\right) = \prod_{x\in V_u^n} \prod_{ y \in V_v^n} \p_\beta \left( x \nsim y \right)
	= \prod_{x\in V_u^n} \prod_{ y \in V_v^n} 
	e^{-\int_{x+\cC} \int_{y+\cC} \frac{\beta}{ \|t-s\|^{2d}} \md s \md t}\\
	&
	= e^{- \sum_{x\in V_u^n} \sum_{ y \in V_v^n} \int_{x+\cC} \int_{y+\cC} \frac{\beta}{ \|t-s\|^{2d}} \md s \md t}
	= e^{-\int_{nu+\left[0,n\right)^d} \int_{nv+\left[0,n\right)^d} \frac{\beta}{ \|t-s\|^{2d}} \md s \md t}\\
	&
	= e^{-\int_{u+\cC} \int_{v+\cC} \frac{\beta}{ \|t-s\|^{2d}} \md s \md t}
	= \p_\beta \left( u \nsim v \right)
\end{align*}
which shows the self-similarity of the model. Also observe that for any $\alpha \in \R_{> 0}$ the process $\alpha \tilde{\e} \coloneqq \left\{(x,y) \in \R^d \times \R^d : \left(\frac{1}{\alpha}x , \frac{1}{\alpha}y \right)\in \tilde{\e} \right\}$ is again a Poisson point process with intensity $\frac{\beta}{2 \|x-y\|^{2d}}$.

\subsection{Notation}
We use the notation $e_i$ for the $i$-th standard unit vector in $\R^d$. For a vector $y \in \R^d$, we write $p_i(y)$ for the $i$-th coordinate of $y$, i.e., $p_i(y) = \langle e_i, y \rangle$. We also use the notation $\mz$ for the vector with all entries equal to 0 and the notation $\mo$ for the vector with all entries equal to 1. When we write $\|u\|$ we always mean the $2$-norm of the vector $u$. We write $\cS_k$ for the set of points $\left\{x \in \Z^d : \|x\|_\infty = k \right\}$ and $\cS_{\geq k}$ for the set $\left\{x \in \Z^d : \|x\|_\infty \geq k \right\}$. For the closed ball of radius $r$ around $x \in \Z^d$ in the $\infty$-norm we use the notation $B_r(x)$, i.e.,  $B_r(x)=\left\{y \in \Z^d : \|x-y\|_\infty \leq r \right\}$. For a vector $u\in \Z^d$ and $n\in \N$, we write 
\begin{align*}
	V_u^n = n\cdot u + \{0,\ldots, n-1\}^d = \prod_{i=1}^d \{n p_i(u),\ldots, np_i(u)+n-1\}
\end{align*}
for the box of side length $n$ shifted by $nu$. When we want to emphasize that we work on certain subgraphs  $A\subset\Z^d$ we will write $D_A \left(x,y\right)$ for the graph distance inside the set $A$, i.e., when we consider edges with both endpoints inside $A$ only. Whenever we write $\dia(A)$ for some set $A\subset \Z^d$ we always mean the inner diameter of this set, i.e., $\dia(A) = \max_{x,y  \in A} D_A(x,y)$. For a graph $(V,E)$ we think of the percolation configuration as a random element $\omega \in \{0,1\}^E$, where we say that the edge $e$ exists or is open or present if $\omega(e)=1$. For $\omega \in \left\{0,1\right\}^E$ and $e \in E$, we define the configuration $\omega^{e-} \in \{0,1\}^E$ by
\begin{equation*} 
\omega^{e-}(\tilde{e}) = \begin{cases}
0 & \tilde{e} = e \\
\omega(e) &  \tilde{e} \neq e
\end{cases} \ ,
\end{equation*}
so this is the configuration where we deleted the edge $e$. For $\omega \in \{0,1\}^E$, we also write $D(u,v;\omega)$ for the graph distance between $u$ and $v$ in the environment represented by $\omega$.

\subsection{Related work} 
The scaling of the graph distance, also called chemical distance or hop-count distance, is a central characteristic of a random graph and has also been examined for many different models of percolation, see for example \cite{hao2021graph, berger2004lower, biskup2004scaling, biskup2011graph, benjamini2001diameter, coppersmith2002diameter, nachmias2008critical, addario2012continuum, ding2013distances, antal1996chemical, ding2018chemical, drewitz2014chemical, garet2007large, hernandez2021chemical, dembin2022variance, biskup2019sharp}. For all dimensions $d$, one can also consider the long-range percolation model with connection probability asymptotic to $\frac{\beta}{\|u-v\|^s}$. When varying the parameter $s$, there are a total of $5$ different regimes, with the transitions happening at $s=d$ and $s=2d$. The value of the first transition $s=d$ is very natural, as the resulting random graph is locally finite if and only if $s>d$. For $s<d$ the graph distance between two points is at most $\lceil \frac{d}{d-s} \rceil$ \cite{benjamini2011geometry}, whereas for $s=d$, the diameter of the box $\{0,\ldots,n\}^d$ is of order $\frac{\log(n)}{\log(\log(n))}$ \cite{coppersmith2002diameter, wu2022sharp}.
In \cite{benjamini2001diameter, biskup2004scaling, biskup2011graph, biskup2019sharp} the authors proved that for $d<s<2d$ the typical distance between two points of Euclidean distance $n$ grows like $\log(n)^\Delta$, where $\Delta^{-1} = \log_2 \left(\frac{2d}{s}\right)$. 
The behavior of the typical distance for long-range percolation on $\Z^d$ also changes at $s=2d$. The reason why $s=2d$ is a critical value is that for $s=2d$ the graph is self-similar, as described in \cref{subsec:cts}.
 For $s>2d$ the graph distance grows at least linearly in the Euclidean distance of two points, as proven in \cite{berger2004lower}.
In \cite{ding2013distances} it is shown that the typical distance for $d=1, s=2$ grows like $n^\dxp$ for some $\dxp\in (0,1)$. For $d\geq 2$ and $s=2d$ the authors in \cite{coppersmith2002diameter} proved a polynomial upper bound on the graph distance, but no lower bound. In this paper, we show a matching polynomial lower bound for all dimensions $d$, similar to the results of \cite{ding2013distances} in one dimension. For fixed dimension $d$, different characteristics of the function $\beta \mapsto \dxp(\beta)$ are considered in the companion paper \cite{baeumler2022behaviour}, where it is shown that the exponent is continuous and strictly decreasing in $\beta$. Together with the fact $\dxp(0)=1$ and Theorem \ref{theo:largebeta} this shows that the long-range percolation model for $s=2d$ interpolates between linear growth and subpolynomial growth as $\beta$ goes from $0$ to $+\infty$.\\
Another line of research is to investigate what happens when one drops the assumption that $p(\beta,\{u,v\})=1$ for all nearest neighbor edges $\{u,v\}$, but assigns i.i.d. random variables to the nearest neighbor edges instead. For $d=1$, there is a change of behavior at $s=2$. As proven in \cite{aizenman1986discontinuity, schulman1983long} or \cite{duminil2020long}, an infinite cluster can not emerge for $s > 2$ and for $s=2, \beta \leq 1$, no matter how small $\p \left(k \nsim k+1\right)$ is. On the other hand, an infinite cluster can emerge for $s<2$ and $s=2,\beta > 1$ (see \cite{newman1986one}). In \cite{aizenman1986discontinuity} the authors proved that there is a discontinuity in the percolation density for $s=2$, contrary to the situation for $s<2$, as proven in \cite{berger2002transience,hutchcroft2021power}. For models, for which an infinite cluster exists the behavior of the percolation model at and near criticality is also a well-studied question (cf. \cite{barsky1991percolation, borgs2005random, berger2002transience, hutchcroft2021power, hutchcroft2021critical, damron2017chemical, hutchcroft2022sharp,baumler2022isoperimetric}).
It is not known up to now how the typical distance in long-range percolation grows for $s=2d$ and $p(\beta,\{u,v\})<1$ for nearest-neighbor edges $\{u,v\}$, but we conjecture also a polynomial growth in the Euclidean distance, whenever an infinite cluster exists.

\noindent
\paragraph*{Acknowledgements.} I thank Noam Berger for  many useful discussions. I thank Yuki Tokushige for helpful comments on an earlier version of this paper. I thank Christophe Garban and Tom Hutchcroft for a discussion leading to \Cref{theo:tail behavior}. I thank an anonymous referee for very many helpful remarks and comments. This work is supported by TopMath, the graduate program of the Elite Network of Bavaria and the graduate center of TUM Graduate School.

\section{Asymptotic behavior of the distance exponent for large $\beta$}

In this chapter, we prove Theorem \ref{theo:largebeta}. On the way, in Section \ref{subsec:bounds on connection p}, we prove several elementary bounds on connection probabilities between certain points and boxes in the long-range percolation graph that will also be used in the following sections. In Section \ref{sec:submultiplicativity}, we prove a submultiplicative structure of the expected distance between two points, leading to the existence of a distance exponent, and also to the inverse logarithmic upper bound in Theorem \ref{theo:largebeta}. In Section \ref{subsec:spacing}, we show that vertices inside a box are not connected to more than one box that is far away, typically. This is necessary in order to prove strict positivity of the distance exponent $\dxp(\beta)$ in Section \ref{subsec:lower bound proof}, and the lower bound on $\dxp(\beta)$ in Theorem \ref{theo:largebeta}.
 
\subsection{Bounds on connection probabilities}\label{subsec:bounds on connection p}

\begin{lemma}
	For all $\beta \geq 0$, all $n\in \N$, and all $u,v\in \Z^d$ with $\|u-v\|_\infty \geq 2$, one has the upper bound
	\begin{align}\label{eq:connectionupper bound}
		\p_\beta \left(u\sim v\right) = \p_\beta \left(V_u^n \sim V_v^n\right)  \leq \frac{2^{2d}\beta}{\|u-v\|_\infty^{2d}},
	\end{align}
	and one has the lower bound
	\begin{align}\label{eq:lowerbound connection prob}
	\p_\beta \left(u\sim v\right)  = \p_\beta \left(V_u^n \sim V_v^n\right)  \geq \frac{(4d)^{-2d}\beta}{\|u-v\|_\infty^{2d}} \wedge \frac{1}{2} \text .
	\end{align}
	For all $k\geq 2$ one has
	\begin{align}\label{eq:point to distance greater k}
		\p_\beta \left( \mz \sim \cS_{\geq k} \right) \leq \beta 50^d k^{-d},
	\end{align}
	and for $m\in \N$, any vertex $x \in V_\mz^m$, and a box $V_w^m$ with $\|w\|_\infty \geq 2$ one has
	\begin{align}\label{eq:point to box bound}
		\p_\beta \left( x \sim V_w^m \right) \leq \frac{\beta 4^{2d} }{\|w\|_\infty^{2d}m^d} \text .
	\end{align}
\end{lemma}

\begin{proof}
The equality $\p_\beta \left(u\sim v\right) = \p_\beta \left(V_u^n \sim V_v^n\right)$ is clear from the discussion about the underlying continuous model in Section \ref{subsec:cts}. We start with the proof of \eqref{eq:connectionupper bound}. Applying the inequalities $1-e^{-x} \leq x$ and $\|\cdot\|_2 \geq \|\cdot\|_\infty$, we get that for two vertices $u,v$ with $\|u-v\|_\infty \geq 2$
\begin{align}\label{eq:bound gleich}
	&\notag \p_\beta \left( u \sim v \right) 
	 = 1-e^{-\beta \int_{u+\cC} \int_{v+\cC} \frac{1}{\|x-y\|^{2d}} \md x \md y }
	\leq \beta \int_{u+\cC} \int_{v+\cC} \frac{1}{\|x-y\|^{2d}} \md x \md y \\
	& \leq \beta \int_{u+\cC} \int_{v+\cC} \frac{1}{\|x-y\|_\infty^{2d}} \md x \md y \leq \frac{\beta}{\left(\|u-v\|_\infty -1 \right)^{2d}}
	 \leq \frac{2^{2d}\beta}{ \|u-v\|_\infty^{2d}} \text .
\end{align}
In order to bound the connection probability between $u$ and $v$ from below, first observe that $\|x\|_2 \leq \|x\|_1 \leq d \|x\|_\infty$ for all $x\in \R^d$. Thus we have 
\begin{align*}
& \int_{u+\cC} \int_{v+\cC} \frac{1}{\|t-s\|^{2d}} \md s \md t
\geq 
d^{-2d}\int_{u+\cC} \int_{v+\cC} \frac{1}{\|t-s\|_\infty^{2d}} \md s \md t\\
& \geq
d^{-2d} \int_{u+\cC} \int_{v+\cC} \frac{1}{\left(\|u-v\|_\infty+1 \right)^{2d}} \md s \md t
\geq
(2d)^{-2d} \frac{1}{\|u-v\|_\infty^{2d}} 
\end{align*}
and this already gives
\begin{align*}
\p_\beta (u\sim v) \geq 1-e^{-(2d)^{-2d} \frac{\beta }{\|u-v\|_\infty^{2d}}}
\geq \frac{(4d)^{-2d} \beta }{\|u-v\|_\infty^{2d}} \wedge \frac{1}{2}
\end{align*}
as $1-e^{-x} \geq \frac{x}{2}\wedge \frac{1}{2}$ for all $x \in \R_{\geq 0}$. So we showed \eqref{eq:lowerbound connection prob}.\\

For each point $x\in \cS_k = \{z \in \Z^d : \|z\|_\infty = k\}$, at least one of its coordinates $p_i(x)$ equals $-k$ or $+k$. All other coordinates can be any integer between $-k$ and $+k$. Thus we can bound the cardinality of the set $\cS_k$ by $\left|\cS_k\right|\leq 2d (2k+1)^{d-1}$. In \eqref{eq:bound gleich} we showed that $\p_\beta \left( \mz \sim x  \right) \leq \frac{\beta}{(\|x\|_\infty-1)^{2d}}$. This already implies that for $k \geq 2$
\begin{align*}
	\p_\beta \left(\mz \sim \cS_k \right) 
	\leq \sum_{x \in \cS_k} \p_\beta ( \mz \sim x)
	\overset{\eqref{eq:bound gleich}}{\leq} \sum_{x \in \cS_k} \frac{\beta}{\left(\|x\|_\infty-1 \right)^{2d}}
	\leq 2d (2k+1)^{d-1} \frac{\beta}{\left(k-1 \right)^{2d}}
\end{align*}
and thus also 
\begin{align}\label{eq:point to distance greater k helpful bound}
	& \notag \p_\beta \left(\mz \sim \cS_{\geq k} \right) \leq \sum_{k^\prime =k}^{\infty} \p_\beta \left(\mz \sim \cS_{k^\prime} \right)
	\leq \sum_{k^\prime =k}^{\infty} 2d (2k^\prime +1)^{d-1} \frac{\beta}{\left(k^\prime -1 \right)^{2d}}
	\\
	& \leq \sum_{k^\prime =k}^{\infty} 2^d 3^d (k^\prime)^{d-1} \frac{\beta 2^{2d}}{(k^\prime)^{2d}}
	 = \beta 24^d \sum_{k^\prime =k}^{\infty} (k^\prime)^{-d-1} \leq \beta 50^d k^{-d},
\end{align}
which already proves \eqref{eq:point to distance greater k}. For $m \in \N$, a vertex $x \in V_\mz^m$, and a box $V_w^m$ with $\|w\|_\infty \geq 2$, we have for all $z \in V_w^m$ that $\|x-z\|_\infty \geq (\|w\|_\infty-1)m$. This implies
\begin{align*}
	\p_\beta\left( x \sim V_w^m \right) 
	\overset{\eqref{eq:connectionupper bound}}{\leq} \sum_{z \in V_w^m} \frac{2^{2d}\beta}{ \|x-z\|_\infty^{2d}}
	\leq \sum_{z \in V_w^m} \frac{2^{2d} \beta}{\left( (\|w\|_\infty-1) m\right)^{2d}} \leq \frac{\beta 4^{2d}}{\|w\|_\infty^{2d} m^d} \text , 
\end{align*}
which shows \eqref{eq:point to box bound}. 
\end{proof}

We will often condition on the event that two blocks $V_u^n, V_v^n$ are connected. So if we write $X$ for the number of edges between them, we condition on the event $X\geq 1$. This conditioning clearly increases the expected number of edges between $V_u^n$ and $V_v^n$, but by at most $+1$, as shown in the next lemma.

\begin{lemma}\label{lem:number of edges conditioned}
	Let $u,v\in \Z^d$ with $u\neq v$ and let $X$ be the number of edges between the blocks $V_u^n$ and $V_v^n$. Then for all $\beta >0$
	\begin{align*}
		\E_\beta\left[X | X\geq 1\right] \leq 1+ \E_\beta \left[ X \right].
	\end{align*}
\end{lemma}

\begin{proof}
	The random variable $X$ is a sum of independent Bernoulli random variables and we prove the statement for all random variables of this type. We use the notation $X= \sum_{i=1}^{m} X_i$, where $m\in \N$, and $(X_i)_{i\in \{1,\ldots, m\}}$ are independent Bernoulli random variables. For $i \in\{1,\ldots,m\}$, let $A_i$ be the event that $X_i = 1$ and $X_j = 0$ for all $j \in \{1,\ldots, i-1\}$. As $\{X\geq 1\}$ implies that there is a first index $i$ such that $X_i=1$, we get that
	\begin{align*}
		\left\{X \geq 1 \right\} = \bigsqcup_{i=1}^m A_i \text , 
	\end{align*}
	where the symbol $\bigsqcup$ means a disjoint union.
	On the event $A_i$, we know that all the random variables $X_j$ with $j<i$ equal $0$, but we have no information about random variables $X_j$ with $j>i$. Thus we get that
	\begin{align*}
		\frac{\E_\beta\left[X \mathbbm{1}_{A_i}\right]}{\p_\beta \left(A_i\right)}
		=
		\E_\beta \left[X | A_i\right] = \E_\beta \left[ 1+  \sum_{j=i+1}^{m} X_j \Big| A_i\right] = 1 + \E_\beta \left[  \sum_{j=i+1}^{m} X_j\right] \leq 1 + \E_\beta \left[X\right] \text .
	\end{align*}
	Multiplying by $\p_\beta \left(A_i\right)$ on both sides of this inequality we get that $\E_\beta\left[X \mathbbm{1}_{A_i}\right] \leq  \p_\beta\left( A_i \right)  \left(1+\E_\beta\left[X\right]\right)$. As the events $(A_i)_{i\in \{1,\ldots,m\}}$ are disjoint, we finally get that
	\begin{align*}
		\E_\beta\left[X | X\geq 1\right] & = 
		\frac{ \E_\beta\left[X \mathbbm{1}_{\left\{X\geq 1\right\}}\right] }{\p_\beta \left(X \geq 1\right)}
		=
		\frac{ \sum_{i=1}^{m} \E_\beta\left[X \mathbbm{1}_{A_i}\right] }{\p_\beta \left(X \geq 1\right)} \\
		&
		\leq
		\frac{ \sum_{i=1}^{m} \p_\beta\left( A_i \right) \left(1+\E_\beta\left[X\right]\right) }{\p_\beta \left(X \geq 1\right)} = 1 +\E_\beta\left[X\right] \text .
	\end{align*}
\end{proof}

\subsection{Submultiplicativity and the upper bound in Theorem \ref{theo:largebeta}}\label{sec:submultiplicativity}

In this section, we prove the submultiplicative structure in the model in \Cref{lem:submultiplicativity}. This allows us to define the distance growth exponent $\dxp(\beta)$ and also helps to prove the upper bound on $\dxp(\beta)$ in \Cref{theo:largebeta}.

\begin{lemma}\label{lem:submultiplicativity}
	For all dimensions $d$ and all $\beta \geq 0$ the sequence
	\begin{equation}\label{eq:Lambda}
	\Lambda(n) = \Lambda(n,\beta)\coloneqq\max_{ u,v \in \left\{0, \ldots ,n-1\right\}^d} \E_\beta \left[D_{V_\mz^n}(u,v)\right] +1
	\end{equation}
	is submultiplicative and for all $\beta \geq 0$
	\begin{align*}
	\dxp(\beta) = \inf_{n\geq 2} \frac{\log\left( \Lambda(n, \beta) \right) }{\log(n)} = \lim_{n\to \infty} \frac{\log\left( \Lambda(n, \beta) \right) }{\log(n)} \text.
	\end{align*}
\end{lemma}

\begin{proof}
	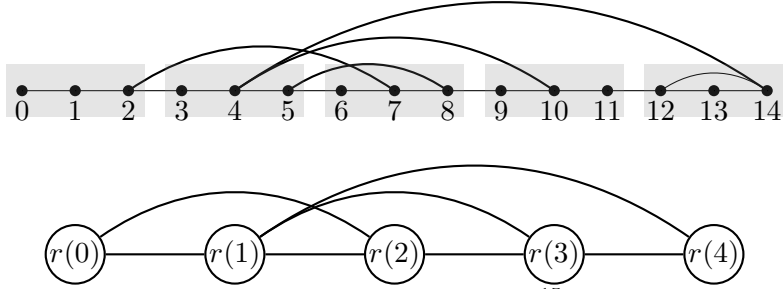
\begin{figure}
		
		\begin{tikzpicture}[scale=0.7]
			\foreach \x in {0,1,...,14} {
				\fill (\x,0) circle (3pt);
			}
			\draw[-] (0,0) -- (14,0) node[right] {};
			\draw[thick] (2,0) to[out=35,in=145] (7,0);
			\draw[thick] (4,0) to[out=35,in=145] (10,0);
			\draw[thick] (5,0) to[out=35,in=145] (8,0);
			\draw[thick] (4,0) to[out=35,in=145] (14,0);
			\draw[black] (12,0) to[out=35,in=145] (14,0);
			\filldraw[gray,fill opacity=0.2, draw=none] (-0.3,-0.5) rectangle (2.3,0.5);
			\filldraw[gray,fill opacity=0.2, draw=none] (2.7,-0.5) rectangle (5.3,0.5);
			\filldraw[gray,fill opacity=0.2, draw=none] (5.7,-0.5) rectangle (8.3,0.5);
			\filldraw[gray,fill opacity=0.2, draw=none] (8.7,-0.5) rectangle (11.3,0.5);
			\filldraw[gray,fill opacity=0.2, draw=none] (11.7,-0.5) rectangle (14.3,0.5);
			\foreach \x in {0,1,...,14}
			\draw (\x,0pt) -- (\x,-0pt) node[anchor=north] {\x};
			
			\vertex[thick](0) at (1,-3.1) {$r(0)$};
			\vertex[thick](1) at (4,-3.1) {$r(1)$};
			\vertex[thick](2) at (7,-3.1) {$r(2)$};
			\vertex[thick](3) at (10,-3.1) {$r(3)$};
			\vertex[thick](4) at (13,-3.1) {$r(4)$};
			
			\path[thick] (0) edge (1) (1) edge (2) (2) edge (3) (3) edge (4);
			\path[thick,out=35,in=145] (0) edge (2);
			\path[thick,out=35,in=145] (1) edge (3);
			\path[thick,out=35,in=145] (1) edge (4);
		\end{tikzpicture}
		\centering
		\parbox{11cm}{\caption{Let $d=1$. The graph $V_\mz^{mn}=V_\mz^{15}$ with $m=5, n=3$ (top) and the graph $G^\prime$ (below). To find a path between two points, say $0$ and $13$, first find the shortest path between $r(0)$ and $r(4)$ in $G^\prime$, then fill the gaps.} \label{fig:renorm}}
	\end{figure}

	We show \eqref{eq:Lambda} using a renormalization argument. As before, we define $V_u^n= \prod_{i=1}^d \left\{p_i(u) n, \ldots ,(p_i(u)+1)n-1\right\}$. The graph $G^\prime$ obtained by identifying all the vertices in $V_u^n$ to one vertex $r(u)$ has the same connection probabilities as the original model. An example of this construction is given in Figure \ref{fig:renorm}. For $x,y \in \left\{0, \ldots ,mn-1\right\}^d$, say with $x \in V_u^n$ and $y\in V_w^n$, we create a path from $x$ to $y$ as follows. First we consider the shortest path $\mathcal{P}= \left( r(u_0)=r(u), r(u_1), \ldots , r(u_{l-1}), r(u_l) = r(w) \right)$ from $r(u)$ to $r(w)$ in $G^\prime$, where $l=D_{G^\prime}(r(u),r(w))$ is the distance between $r(u)$ and $r(w)$ in the renormalized model. Inside $V_{u_i}^n$, we first fix two vertices $z_i$ and $v_i$ such that $z_i \sim V_{u_{i-1}}^n$ and $v_i \sim V_{u_{i+1}}^n$; for $i=0$ set $z_0=x$ and for $i=l$ set $v_l=y$. In case there are several such vertices $z_i$ and $v_i$, we choose the one with smallest coordinates, where we weigh the coordinates in decreasing order (any deterministic rule that does not depend on the environment would work here). For each $i$, there clearly is a path between $z_i$ and $v_i$ that is completely inside $V_{u_i}^n$. As no information has been revealed up to now about the edges with both endpoints inside $V_{u_i}^n$, the expected distance between $v_i$ and $z_i$ inside $V_{u_i}^n$ is at most
	\begin{align*}
		\max_{ a,b \in  V_{u_i}^n } \E_\beta \left[D_{V_{u_i}^n}(a,b)\right] = \Lambda(n,\beta)-1 \text .
	\end{align*}
	Now we glue all these paths together to get a path from $x$ to $y$. To bound the total distance between $x$ and $y$ note that we have $l+1$ sets $V_{u_i}^n$ in which we need to find a path between two vertices. Additionally, we need to add $+l$ for the steps that we make from $V_{u_i}^n$ to $V_{u_{i+1}}^n$ for $i=0,\ldots,l-1$. Thus we get that
	\begin{align*}
		\E_\beta\left[D_{V_\mz^{mn}}(x,y) \ \big| \ D_{G^\prime}(r(u),r(w)) = l\right] \leq (l+1) \max_{ a,b \in \left\{0, \ldots ,n-1\right\}^d}\E_\beta \left[D_{V_\mz^n}(a,b)\right] + l \text .
	\end{align*}
	Taking expectations on both sides of this inequality yields
	\begin{align*}
	&\E_\beta\left[D_{V_\mz^{mn}}(x,y)\right] \\
	& \leq \left(\E_\beta \left[D_{G^\prime}(r(u),r(w))\right] +1\right) 
	\max_{ a,b \in \left\{0, \ldots ,n-1\right\}^d}\E_\beta \left[D_{V_\mz^n}(a,b)\right] + \E_\beta \left[D_{G^\prime}(r(u),r(w))\right] 
	\\
	&
	= \left(\E_\beta \left[D_{V_\mz^m}(u,w)\right] +1\right) \left( \max_{ a,b \in \left\{0, \ldots ,n-1\right\}^d} \E_\beta \left[ D_{V_\mz^n}(a,b)\right] +1 \right) - 1
	\\
	& 
	\leq\Lambda(m) \Lambda(n) - 1 \text.
	\end{align*}
	As $x,y \in \left\{0, \ldots ,nm-1\right\}^d$ were arbitrary we obtain
	\begin{equation}\label{eq:submultiplicativedist}
	\Lambda(mn) \leq \Lambda(m)\Lambda(n) \text,
	\end{equation}
	and as the sequence is submultiplicative we can define
	\begin{align*}
		\dxp=\dxp(\beta) = \lim_{k \to \infty} \frac{\log\left(\Lambda(2^k,\beta)\right)}{\log(2^k)} \text .
	\end{align*}
	Actually, this limit exists not just along dyadic points of the form $2^k$, for $k\in \N$, but even when taking a limit along the integers, i.e.,
	\begin{align*}
	\dxp=\dxp(\beta) = \lim_{n \to \infty} \frac{\log\left(\Lambda(n,\beta)\right)}{\log(n)} \text,
	\end{align*}
	which follows from \Cref{propo:scaling} below.
	As a next step, we want to show that $\Lambda(n) \geq n^\dxp$ for all $n$. We do this using a proof by contradiction. So assume the contrary, i.e., that for some $\beta \geq 0$ there exists a natural number $N$ and a $c<1$ with $\Lambda(N,\beta)= c N^{\dxp(\beta)}$. Using \eqref{eq:submultiplicativedist} we get that for every integer $k$
	\begin{equation*}
		\Lambda(N^k,\beta) \leq \Lambda(N)^k = c^k N^{\dxp(\beta) k}
	\end{equation*}
	and thus
	\begin{align*}
		\dxp(\beta) = \lim_{k \to \infty} \frac{\log\left(\Lambda(N^k,\beta)\right)}{\log(N^k)} \leq \limsup_{k\to \infty} \frac{\log(c^k N^{\dxp(\beta) k})}{\log (N^k)} = \frac{\log(c) + \dxp \log(N)}{\log(N)} < \dxp(\beta)
	\end{align*}
	which is a contradiction. Knowing this already gives us that for all positive numbers $K$ we have
	\begin{equation}\label{limitisinf}
	\dxp= \lim_{n\to \infty} \frac{\log(\Lambda(n))}{\log(n)} = \inf_{n \geq 2} \frac{\log(\Lambda(n))}{\log(n)}= \inf_{n \geq K} \frac{\log(\Lambda(n))}{\log(n)} \text .
	\end{equation}
\end{proof}
This lemma and its proof already have several interesting applications. First, we emphasize that $\Lambda(mn,\beta) \geq \Lambda(n,\beta)$ for all $m,n \in \N_{>0}$. This holds, as for arbitrary $x,y \in V_\mz^n$, the distance $D_{V_\mz^{mn}}(u,v)$ between $u \in V_x^m$ and $v \in V_y^m$ is at least the distance between $r(x)$ and $r(y)$ in $G^\prime$. Using the self-similarity and taking expectations we thus get that
\begin{align*}
	\E_\beta \left[ D_{V_\mz^{mn}}(u,v) \right] \geq 	\E_\beta \left[ D_{G^\prime} (r(x),r(y)) \right]
	=
	\E_\beta \left[ D_{V_\mz^n} (x,y) \right]
\end{align*}
which shows our claim. For $n=3$, we have for all $u,v \in \{0,1,2\}^d$ with $u\neq v$, and for all $\beta >0$ that
\begin{align*}
	\E_\beta \left[D_{[0,2]^d} (u,v)\right] = 1\cdot \p_\beta(u\sim v) + 2 \cdot \p_\beta(u\nsim v) < 2
\end{align*} 
and this already implies that $\Lambda(3) \eqqcolon 3^{\dxp^\prime } < 3$ for some $\dxp^\prime = \dxp^\prime(\beta) < 1$. Inductively, with a renormalization at scale $3$, we get that
\begin{align}\label{eq:iterative bound 3}
	\Lambda(3^k N) \leq \Lambda(3)^k \Lambda(N) = 3^{k\dxp^\prime} \Lambda(N)
\end{align}
for all $k,N \in \N_{>0}$. This inequality already gives the upper bound in expectation for $s=2d$, that was already observed in \cite{coppersmith2002diameter} with a very similar technique. Next, we do a renormalization at scale $\sqrt[2d]{\beta}$ instead of scale $3$ in order to get the inverse logarithmic upper bound stated in Theorem \ref{theo:largebeta}.

\begin{proof}[Proof of the upper bound in Theorem \ref{theo:largebeta}]
	We want to show that for each dimension $d$ there exists a constant $C<\infty$ such that for all $\beta \geq 2$
	\begin{align*}
		\dxp(\beta)\leq \frac{C}{\log(\beta)} \text .
	\end{align*}
	As the connection probability  $\p_\beta \left(u \sim v\right)$  between any two vertices $u,v\in \Z^d$ is increasing in $\beta$, the distance exponent $\dxp : \R_{\geq 0} \to \left[0,1\right]$ is clearly decreasing by the {\sl Harris coupling}, see for example \cite{heydenreich2017progress}. Thus it suffices to show the upper bound for $\beta$ large enough with $\sqrt[2d]{\beta} \in \N$. For such a $\beta $ and all $u,v \in \{0, \ldots, \sqrt[2d]{\beta} - 1 \}^d$, we have for all $y \in u + \cC$ and $x \in v + \cC$ that
	\begin{align}\label{eq:inequality for points}
		\|x-y\|^{2d} \leq d^{2d} \|x-y\|_\infty^{2d}  
		\leq d^{2d} \sqrt[2d]{\beta}^{2d}
		= d^{2d} \beta
	\end{align}
	and this already implies
	\begin{align*}
		\int_{u+\cC} \int_{v+ \cC} \frac{1}{\|x-y\|^{2d}} \md x \md y \geq \frac{1}{d^{2d} \beta} \text .
	\end{align*}
	Inserting this into the definition $p(\beta, \{u,v\})$ and using that $1-e^{-x}\geq \frac{x}{2}$ for all $x\leq 1$ we get that for large enough $\beta$ that satisfies $\frac{1}{d^{2d} \beta } \leq 1$ 
	\begin{align}\label{eq:zconnection upper bound}
		\p_\beta (u\sim v) = 1-e^{-\beta 	\int_{u+\cC} \int_{v+ \cC} \frac{1}{\|x-y\|^{2d}} \md x \md y } \overset{\eqref{eq:inequality for points}}{\geq}
		1-e^{-d^{-2d} } \geq \frac{1}{2} d^{-2d} \geq (2d)^{-2d}
	\end{align}
	for all $u,v \in \{0,\ldots,\sqrt[2d]{\beta}-1\}^d$. Next, we bound the expected graph distance between $u$ and $v$. We do this by comparing the distance to a geometric random variable. Let $(u=u_0,u_1,\ldots, u_k=v)$ be a deterministic self-avoiding path from $u$ to $v$ inside $V_{\mz}^{\sqrt[2d]{\beta}}$, with $k \leq  \sqrt[2d]{\beta}$ and $\|u_{i}- u_{i-1}\|_\infty =1$ for all $i \in \{1,\ldots,k\}$. Starting from this, we build a shorter path between $u$ and $v$ as follows. We start at $u_0=u$. Then for $i=0,\ldots,k-1$, if $u_i \sim v$, directly go to $v$. If $u_i \nsim v$, then go to $u_{i+1}$. This gives a path $P$ between $u$ and $v$, and for $l\in \{1,\ldots,k\}$ this path has length of at least $l$ if and only if $u_{i} \nsim v$ for all $i\in \{0,\ldots,l-2\}$. As the connections between $v$ and different $u_i$-s are independent we get that
	\begin{align*}
		\E_\beta \left[D_{V_{\mz}^{\sqrt[2d]{\beta}} } (u,v) \right] &
		= \sum_{l=1}^{k} \p_\beta \left(D_{V_{\mz}^{\sqrt[2d]{\beta}} } (u,v) \geq l \right) \leq \sum_{l=1}^{k} \p_\beta \left(u_i \nsim v \text{ for all } i\leq l-2\right)\\
		& \overset{\eqref{eq:zconnection upper bound}}{\leq} \sum_{l=1}^{k} \left(1-(2d)^{-2d}\right)^{l-1} \leq \frac{1}{1- \left(1-(2d)^{-2d}\right) } = (2d)^{2d} .
	\end{align*}
	This already implies that $\Lambda\left(\sqrt[2d]{\beta} , \beta \right) \leq (2d)^{2d} + 1 \leq (3d)^{2d}$.	Applying the submultiplicativity of $\Lambda$ iteratively we get that
	\begin{align*}
		\dxp(\beta) & = \lim_{k\to \infty} \frac{\log\left(\Lambda\left( \sqrt[2d]{\beta}^k , \beta \right)\right)}{\log\left(\sqrt[2d]{\beta}^k\right)}
		\leq \limsup_{k\to \infty} \frac{\log\left(\Lambda\left( \sqrt[2d]{\beta} , \beta \right)^k\right)}{\log\left(\sqrt[2d]{\beta}^k\right)}\\
		& = \frac{\log\left( \Lambda\left( \sqrt[2d]{\beta} , \beta \right) \right)}{\log\left(\sqrt[2d]{\beta} \right)} \leq \frac{ \log \left( (3d)^{2d} \right) }{ \frac{1}{2d} \log(\beta)} = \frac{4d^2 \log \left( 3d \right) }{\log(\beta)}
	\end{align*}
	which finishes the proof.
\end{proof}

\subsection{Spacing between points with long bonds}\label{subsec:spacing}

In this section, we investigate certain geometric properties of the cluster inside certain boxes. Mostly, we want to get upper bounds on the probability that a vertex is connected to two different long edges. As we will need it at a later point, namely in Section \ref{sec:distance of long edge points}, we will prove the statements for $\|x-y\|_\infty \leq 1$ instead of $x=y$. This does not cause major difficulties, as for each point $x \in \Z^d$, there are only $3^d$ many points $y \in \Z^d$ with $\|x-y\|_\infty \leq 1$. We start with showing that the probability that two vertices $x,y$ with $\|x-y\|_\infty \leq 1$ are both connected to far away boxes is very low.

\begin{lemma}\label{lem:separation both far}
	For blocks $V_u^m, V_v^m, V_w^m$ with $\|u-v\|_\infty , \|v-w\|_\infty \geq 2$, there exists a constant $C_d < \infty$ such that for all $\beta \geq 0$
	\begin{align*}
		\p_\beta \left( \exists x,y \in V_v^m : \|x-y\|_\infty \leq 1, x \sim V_u^m, y \sim V_w^m \right) \leq \frac{C_d \beta^2}{\|u-v\|_\infty^{2d} \|w-v\|_\infty^{2d} m^d} \text .
	\end{align*}
\end{lemma}
\begin{proof}
	By translational invariance we can assume that $v=\mz$, and thus $\|u\|_\infty , \|w\|_\infty \geq 2$. For each $x \in V_{\mz}^m$ there are at most $3^d$ vertices $y\in V_{\mz}^m$ with $\|x-y\|_\infty \leq 1$. For $x,y \in V_{\mz}^m$ the probability that $y \sim V_w^m$  is bounded by $\frac{\beta4^{2d}}{\|w\|_\infty^{2d}m^d}$, and the probability that $x\sim V_u^m$ is bounded by $\frac{\beta4^{2d}}{\|u\|_\infty^{2d}m^d}$, using \eqref{eq:point to box bound}. Thus 
	\begin{align*}
		& \notag \p_\beta \left( \exists x,y \in V_{\mz}^m : \|x-y\|_\infty \leq 1, x \sim V_u^m, y \sim V_w^m \right) 
		\leq \sum_{x \in V_{\mz}^m} \p_\beta \left( x \sim V_u^m \right) 3^d \frac{\beta 4^{2d}}{\|w\|_\infty^{2d}m^d}\\
		&
		= \frac{48^d \beta}{\|w\|_\infty^{2d} m^d} \sum_{x \in V_{\mz}^m} \p_\beta \left( x \sim V_u^m \right)
		\overset{\eqref{eq:point to box bound}}{\leq}  \frac{48^d \beta}{\|w\|_\infty^{2d} m^d} m^d \frac{\beta 4^{2d}}{\|u\|_\infty^{2d} m^d}
		\leq
		\frac{\beta^2 1000^d}{\|w\|_\infty^{2d} \|u\|_\infty^{2d} m^d}
	\end{align*}
	which finishes the proof.
\end{proof}

\begin{lemma}\label{lem:separation one far}
	For blocks $V_u^m, V_v^m, V_w^m$ with $ \|v-w\|_\infty \geq 2$ and $\|u-v\|_\infty = 1$, there exists a constant $C_d < \infty$ such that for all $\beta \geq 0$
	\begin{align*}
	\p_\beta \left( \exists x,y \in V_v^m : \|x-y\|_\infty \leq 1, x \sim V_u^m, y \sim V_w^m \right) \leq 
	 \begin{cases}
	\frac{C_d \beta \lceil \beta \rceil \log(m)}{\|v-w\|_\infty^{2d} m} & \text{ for } d=1\\
	\frac{C_d \beta \lceil \beta \rceil}{\|v-w\|_\infty^{2d} m} & \text{ for } d \geq 2
	\end{cases}  \ .
	\end{align*}
\end{lemma}
\begin{proof}
	By translational invariance we can assume that $v=\mz$, and thus $\|u\|_\infty = 1, \|w\|_\infty \geq 2$. For each $x \in V_{\mz}^m$ there are at most $3^d$ vertices $y\in V_{\mz}^m$ with $\|x-y\|_\infty \leq 1$. For each vertex $y \in V_{\mz}^m$ the probability that $y \sim V_w^m$ is bounded by $\frac{\beta4^{2d}}{\|w\|_\infty^{2d} m^d}$ by \eqref{eq:point to box bound}.
	Thus 
	\begin{align}\label{eq:near and far bound 1}
	& \notag \p_\beta \left( \exists x,y \in V_{\mz}^m : \|x-y\|_\infty \leq 1, x \sim V_u^m, y \sim V_w^m \right) 
	\leq \sum_{x \in V_{\mz}^m} \p_\beta \left( x \sim V_u^m \right) 3^d \frac{\beta 4^{2d}}{\|w\|_\infty^{2d} m^d}\\
	&
	= \frac{48^d \beta}{\|w\|_\infty^{2d} m^d} \sum_{x \in V_{\mz}^m} \p_\beta \left( x \sim V_u^m \right) \text .
	\end{align}
	As $\|u\|_\infty=1$ we have $D_\infty(x,V_u^m) \leq m$ for all $x \in V_{\mz}^m$, where $D_\infty$ is the distance with respect to the $\infty$-norm. We furthermore have the inequality
	\begin{align*}
	\left| \left\{ x \in V_{\mz}^m : D_\infty \left(x, V_u^m\right) = k \right\} \right| \leq 6^d m^{d-1}
	\end{align*}
	for all $k \in \N$. This is clear for $k>m$, as the relevant set is empty in this case. For $k \leq m$ the set $ \left\{ x \in \Z^d : D_\infty \left(x, V_u^m\right) = k \right\}$ is just the boundary of the box 
	\begin{equation*}
	\prod_{i=1}^{d} \{p_i(u) m - k, \ldots, (p_i(u)+1) m - 1 + k\} \ ,
	\end{equation*}
	which is a box of side length $m+2k \leq 3m$. Thus the boundary has a cardinality of at most $2d (3m)^{d-1} \leq 6^d m^{d-1}$. Using this observation we get that
	\begin{align}\label{eq:near box expected adjacents}
	&\notag \sum_{x \in V_{\mz}^m} \p_\beta \left( x \sim V_u^m \right) = \sum_{k=1}^{m} \sum_{\substack{x \in V_{\mz}^m : \\ D_\infty(x,V_u^m)=k}} \p_\beta \left(x\sim V_u^m \right) \leq \sum_{k=1}^{m} 6^d m^{d-1} \p_\beta \left(x \sim x+\cS_{\geq k} \right) \\
	&
	\overset{\eqref{eq:point to distance greater k}}{\leq}
	6^d m^{d-1} 
	+ 6^d m^{d-1} \sum_{k=2}^{m} \beta 50^d k^{-d} 
	\leq \begin{cases}
	m^{d-1} \log(m) \lceil \beta \rceil 400^d & \text{ for } d = 1\\
	m^{d-1} \lceil \beta \rceil 400^d & \text{ for } d\geq 2
	\end{cases} \ .
	\end{align}
	Inserting this into \eqref{eq:near and far bound 1}, we get that
	\begin{align*}
	& \notag \p_\beta \left( \exists x,y \in V_{\mz}^m : \|x-y\|_\infty \leq 1, x \sim V_u^m, y \sim V_w^m \right) 
	\leq \sum_{x \in V_{\mz}^m} \p_\beta \left( x \sim V_u^m \right) 3^d \frac{\beta 4^{2d}}{\|w\|_\infty^{2d} m^d}\\
	& 
	\leq \begin{cases}
	\frac{20000^d \beta \lceil \beta \rceil\log(m)}{\|w\|_\infty^{2d}m}
	& \text{ for } d = 1\\
	\frac{20000^d \beta \lceil \beta \rceil}{\|w\|_\infty^{2d}m}
	& \text{ for } d \geq 1
	\end{cases} \ 
	\end{align*}
	which finishes the proof.	
\end{proof}

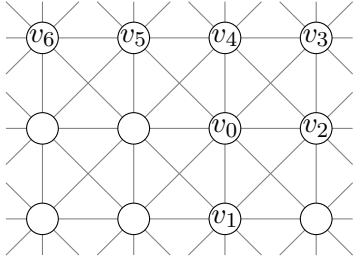
\begin{figure}
	\[\begin{tikzpicture}[scale = 1.2]
		
		\draw[gray, thin] (-0.4,-0.4) grid (3.4,2.4);
		
		\draw[gray, thin] (0-0.4,0+0.4) -- (0+0.4,0-0.4);
		\draw[gray, thin] (0-0.4,1+0.4) -- (1+0.4,0-0.4);
		\draw[gray, thin] (0-0.4,2+0.4) -- (2+0.4,0-0.4);
		\draw[gray, thin] (1-0.4,2+0.4) -- (3+0.4,0-0.4);
		\draw[gray, thin] (2-0.4,2+0.4) -- (3+0.4,1-0.4);
		\draw[gray, thin] (3-0.4,2+0.4) -- (3+0.4,2-0.4);

		\draw[gray, thin] (0-0.4,2-0.4) -- (0+0.4,2+0.4);
		\draw[gray, thin] (0-0.4,1-0.4) -- (1+0.4,2+0.4);
		\draw[gray, thin] (0-0.4,0-0.4) -- (2+0.4,2+0.4);
		\draw[gray, thin] (1-0.4,0-0.4) -- (3+0.4,2+0.4);
		\draw[gray, thin] (2-0.4,0-0.4) -- (3+0.4,1+0.4);
		\draw[gray, thin] (3-0.4,0-0.4) -- (3+0.4,0+0.4);
		
		\vertex[fill=white] (00) at (0,0) {$\textcolor{white}{v_0}$};
		\vertex[fill=white] (01) at (0,1) {$\textcolor{white}{v_0}$};
		\vertex[fill=white] (00) at (0,2) {$v_6$};
		\vertex[fill=white] (10) at (1,0) {$\textcolor{white}{v_0}$};
		\vertex[fill=white] (11) at (1,1) {$\textcolor{white}{v_0}$};
		\vertex[fill=white] (12) at (1,2) {$v_5$};
		\vertex[fill=white] (20) at (2,0) {$v_1$};
		\vertex[fill=white] (21) at (2,1) {$v_0$};
		\vertex[fill=white] (22) at (2,2) {$v_4$};
		\vertex[fill=white] (30) at (3,0) {$\textcolor{white}{v_0}$};
		\vertex[fill=white] (31) at (3,1) {$v_2$};
		\vertex[fill=white] (32) at (3,2) {$v_3$};
		
	\end{tikzpicture}\]
	\centering
	\parbox{11cm}{\caption{An example of points $v_0,\ldots,v_{l+1}$ with $l=5$ as described in Lemma \ref{lem:separation wanders}. $\|v_i-v_{i+1}\|_\infty = 1$ for all $i\in \{0,\ldots,5\}$ and $v_{l+1}$ is the first point with $\|v_{l+1}-v_0\|_\infty > 1.$} \label{fig:points}}
	
\end{figure}

\begin{lemma}\label{lem:separation wanders}
	Let $m\in \N, l\in \{1,\ldots, 3^d-1\} $, and let $v_0,v_1,\ldots, v_{l+1} \in \Z^d$ be distinct with $\|v_{i+1}-v_i\|_\infty = 1$ for all $i\in \{0,\ldots, l\}$, $\|v_{i}-v_0\|_\infty = 1$ for all $i\in \{1,\ldots, l\}$ and $\|v_{l+1}-v_0\|_\infty=2$. (An example of such a sequence of points is given in Figure \ref{fig:points}). Then there exists a constant $C_d<\infty$ such that the two probabilities
	\begin{align*}
		&\p_\beta \left( \exists i \in \{1,\ldots,l\} \exists x,y \in V_{v_i}^m \text{ with } \| x-y\|_\infty \leq 1, x\sim V_{v_{i-1}}^m, y \sim V_{v_{i+1}}^m \cap \left( y + \mathcal{S}_{\geq \frac{m}{6^d}} \right) \right), \\
		&\p_\beta \left( \exists i \in \{1,\ldots,l\} \exists x,y \in V_{v_i}^m \text{ with } \| x-y\|_\infty \leq 1, x\sim V_{v_{i+1}}^m, y \sim V_{v_{i-1}}^m \cap \left( y + \mathcal{S}_{\geq \frac{m}{6^d}} \right) \right)
	\end{align*}
	are both bounded by $\frac{C_d \beta \lceil \beta \rceil \log(m)}{m} $ for $d=1$, respectively by $\frac{C_d \beta \lceil \beta \rceil }{m} $ for $d\geq 2$.
\end{lemma}

\begin{proof}
	By a union bound we have that
	\begin{align*}
		& \p_\beta \left( \exists i \in \{1,\ldots,l\} \exists x,y \in V_{v_i}^m \text{ with } \| x-y\|_\infty \leq 1, x\sim V_{v_{i-1}}^m, y \sim V_{v_{i+1}}^m \cap \left( y + \mathcal{S}_{\geq \frac{m}{6^d}} \right) \right)\\
		& \leq \sum_{i \in \{1,\ldots,l\}} \sum_{\substack{x,y \in V_{v_i}^m : \\ \|x-y\|_\infty\leq 1}} \p_\beta \left(  x\sim V_{v_{i-1}}^m, y \sim V_{v_{i+1}}^m \cap \left( y + \mathcal{S}_{\geq \frac{m}{6^d}} \right) \right)\\
		& \leq \sum_{i \in \{1,\ldots,l\}} \sum_{\substack{x,y \in V_{v_i}^m : \\ \|x-y\|_\infty\leq 1}} \p_\beta \left(  x\sim V_{v_{i-1}}^m \right) \p_\beta \left( y \sim  \left( y + \mathcal{S}_{\geq \frac{m}{6^d}} \right) \right)\\
		& \overset{\eqref{eq:point to distance greater k}}{\leq} \beta 50^d \left(\frac{m}{6^d}\right)^{-d} \sum_{i \in \{1,\ldots,l\}} \sum_{\substack{x,y \in V_{v_i}^m : \\ \|x-y\|_\infty\leq 1}} \p_\beta \left(  x\sim V_{v_{i-1}}^m \right)\\
		&
		\leq
		\beta 150^d \left(\frac{m}{6^d}\right)^{-d} \sum_{i \in \{1,\ldots,l\}} \sum_{x \in V_{v_i}^m } \p_\beta \left(  x\sim V_{v_{i-1}}^m \right) \text .
	\end{align*}
		The sum $ \sum_{x \in V_{v_i}^m } \p_\beta \left(  x\sim V_{v_{i-1}}^m \right)$ was already upper bounded in \eqref{eq:near box expected adjacents}. Using this upper bound, $l\leq 3^d$, and inserting this into the line above we get that
	\begin{align*}
		\beta 150^d \left(\frac{m}{6^d}\right)^{-d} \sum_{i \in \{1,\ldots,l\}} \sum_{x \in V_{v_i}^m } \p_\beta \left(  x\sim V_{v_{i-1}}^m \right)
		\leq 
		\begin{cases}
		\frac{\beta \lceil \beta \rceil \left(6^d 10^6\right)^d\log(m)}{m} & \text{ for } d= 1 \\
		\frac{\beta \lceil \beta \rceil \left(6^d 10^6\right)^d}{m} &
		\text{ for } d \geq 2
		\end{cases}
		\ 
	\end{align*}
	which finishes the proof for the first item in the statement of the lemma. The estimate for the second term works analogously.
\end{proof}

\subsection{The lower bound in Theorem \ref{theo:largebeta}}\label{subsec:lower bound proof}

Finally, we develop all the necessary techniques in order to show the lower bound in Theorem \ref{theo:largebeta}, i.e., that for all dimensions $d$, there exists a constant $c>0$ such that $\dxp(\beta) \geq \frac{c}{\log(\beta)}$ for all $\beta \geq 2$.

\begin{proof}[Proof of the lower bound in Theorem \ref{theo:largebeta}]
	Inequality \eqref{eq:lowerbound connection prob} and Lemma \ref{lem:separation both far} show that for all dimensions $d$ there exists a constant $C_d < \infty$ such that for all $\beta \geq 2$
	and all $u,v,w$ with $\|u-v\|_\infty, \|v-w\|_\infty \geq 2$
	\begin{align}\label{eq:separation bound 1}
		\p_\beta \left( \exists x \in V_v^m : x \sim V_u^m , x \sim V_w^m \ | \ V_u^m \sim V_v^m \sim V_w^m \right) \leq \frac{C_d\beta^2}{m^{1/2}} \ .
	\end{align}
	Analogously, Lemma \ref{lem:separation one far} shows that there exists a constant $C_d < \infty$ such that for all $\beta \geq 2$ and all $u,v,w$ with $\|u-v\|_\infty \geq 2$ and $\|v-w\|_\infty=1$
	\begin{align}\label{eq:separation bound 2}
		\p_\beta \left( \exists x \in V_v^m : x \sim V_u^m , x \sim V_w^m \ | \ V_u^m \sim V_v^m \sim V_w^m \right) \leq \frac{C_d\beta^2}{m^{1/2}} 
	\end{align}
	where we also used that $\frac{\log(m)}{m} = \mathcal{O} \left( m^{-1/2}\right)$.
	Lemma \ref{lem:separation wanders} implies that for every $l\in \{1,\ldots,3^d-1\}$ and $v_0,v_1,\ldots, v_{l+1} \in \Z^d$ distinct with $\|v_{i+1}-v_i\|_\infty = 1$ for all $i\in \{0,\ldots, l\}$, $\|v_{i}-v_0\|_\infty = 1$ for all $i\in \{1,\ldots, l\}$, and $\|v_{l+1}-v_0\|_\infty=2$, one has the bound
	\begin{align}\label{eq:separation bound 3}
		\p_\beta \left( \exists x_1, \ldots, x_l : x_i \in V_{v_i}^m , V_{v_0}^m \sim x_1 \sim x_2 \sim \ldots \sim x_l \sim V_{v_{l+1}}^m \right) \leq \frac{C_d \beta^2}{m^{1/2}}
	\end{align}
	as a path from $V_{v_0}$ to $V_{v_{l+1}}$ in $l+1 \leq 3^d$ steps needs to contain at least one edge $\{x_i,x_{i+1}\}$ with $\|x_i-x_{i+1}\|_\infty \geq \frac{m}{3^d}$ and thus $x_i \sim x_i+\cS_{\geq \frac{m}{6^d}}$ in this case. We will now show that
	\begin{align}\label{eq:lowerboundtoshow1}
		\E_\beta \left[ D_{V_{\mz}^{mM}} \left(\mz, (mM-1)e_1\right) \right] \geq 
		\left( 1 + \frac{1}{3^{d+4}} \right)
		\E_\beta \left[ D_{V_{\mz}^{M}} \left(\mz, (M-1)e_1\right) \right]
	\end{align}
	for $m \geq \left(2000 \cdot \lceil \beta \rceil^3  3^{5d} C_d\right)^{\left(3^{4d}\right)}$ and all large enough $M$. We will see later where this condition on $m$ comes from.  To see \eqref{eq:lowerboundtoshow1}, we use a renormalization. For $u \in V_{\mz}^{M}$, we identify the blocks $V_u^m$ to vertices $r(u)$ and call the resulting graph $G^\prime$. Then we will prove that
	\begin{align*}
		\E_\beta \left[ D_{V_{\mz}^{mM}} \left(\mz, (mM-1)e_1\right) \right] \geq 
		\left( 1 + \frac{1}{3^{d+4}} \right)
		\E_\beta \left[ D_{G^\prime} \left(r(\mz), r((M-1)e_1)\right) \right]
	\end{align*}
	for large enough $M$.
	This implies \eqref{eq:lowerboundtoshow1}, as the random graphs $G^\prime$ and $V_{\mz}^{M}$ have the same distribution, as shown in \cref{subsec:cts}. Now we condition on the graph $G^\prime$, i.e.,  we already have the knowledge which blocks of the form $V_u^m$ are connected in the original graph. Let $P^\prime = (r(v_0),\ldots, r(v_k))$ be a self-avoiding path in $G^\prime$ starting at the origin vertex, i.e.,  $v_0=\mz$. Let $k \geq 3^{d+3}$ and let $l=\big\lfloor \frac{k}{3^{d+1}}\big\rfloor$. For $j \in \{0,\ldots, l\}$, we call the subsequence $R_j \coloneqq \left(r(v_{2j3^{d}}), r(v_{2j3^{d}+1}) , \ldots , r(v_{(2j+2)3^{d}})\right) $ {\sl separated} if there does {\sl not} exist a sequence $\left(x_i\right)_{i=2j3^d + 1}^{(2j+2) \cdot 3^{d} - 1}$ such that $x_i \in V_{v_i}^m \text{ for all } i \in \{2j3^d + 1,\ldots, (2j+2) 3^{d} - 1 \} $ and
	\begin{align*}
		V_{v_{2j3^d}}^m \sim x_{2j3^d+1} \sim x_{2j3^d+2} \sim \ldots  \sim x_{2(j+2)3^d-1} \sim	V_{v_{2(j+2)3^d}}^m \text .
	\end{align*}
	For a given self-avoiding path $P^\prime \subset G^\prime$ and different values of $j \in \{0,\ldots,l\}$, it is independent whether the subsequences $\left(r(v_{2j3^{d}}), r(v_{2j3^{d}+1}) , \ldots , r(v_{(2j+2)3^{d}})\right)$ are separated, and the probability that a specific subsequence $\left(r(v_{2j3^{d}}), r(v_{2j3^{d}+1}) , \ldots , r(v_{(2j+2)3^{d}})\right)$ is not separated is bounded by $\frac{C_d \beta^2}{m^{1/2}}$, as for every sequence $\left(v_{2j3^{d}}, v_{2j3^{d}+1} , \ldots , v_{(2j+2)3^{d}} \right)$ at least one of the situations of \eqref{eq:separation bound 1}, \eqref{eq:separation bound 2} or \eqref{eq:separation bound 3} holds, as we will argue below. Here we say that the situation of $\eqref{eq:separation bound 1}$ holds if there exists an index $i \in \{2j3^d+1,\ldots,(2j+2)3^d-1\}$ such that $\|v_i-v_{i+1}\|_\infty, \|v_i-v_{i-1}\|_\infty \geq 2$, the situation of $\eqref{eq:separation bound 2}$ holds if there exists an index $i \in \{2j3^d+1,\ldots,(2j+2)3^d-1\}$ such that $\|v_i-v_{i+1}\|_\infty = 1, \|v_i-v_{i-1}\|_\infty \geq 2$ or $\|v_i-v_{i-1}\|_\infty = 1, \|v_i-v_{i+1}\|_\infty \geq 2$, and the situation of \eqref{eq:separation bound 3} holds if there exists $l\in \{1,\ldots,3^d-1\}$ such that $\|v_{i+1}-v_i\|_\infty = 1$ for all $i\in \{2j3^d,\ldots, 2j3^d+l\}$, $\|v_{i}-v_{2j3^d}\|_\infty = 1$ for all $i\in \{2j3^d+1,\ldots, 2j3^d+l\}$ and $\|v_{2j3^d+l+1}-v_{2j3^d}\|_\infty=2$.
	If none of the situations in \eqref{eq:separation bound 1},\eqref{eq:separation bound 2} holds, then the path makes only nearest neighbor-jumps within the subsequence $R_j$. However, as that there are only $3^d-1$ many points $v \in \Z^d$ with $\|v-v_{2j3^d}\|_\infty = 1$, the situation of \eqref{eq:separation bound 3} must occur within the subsequence $R_j$ for some $l$. So in total we see that
	\begin{align*}
		\p_\beta \left(R_j \text{ not separated } \big| \ G^\prime  \right) \leq \frac{C_d \beta ^2}{m^{1/2}}\text.
	\end{align*}
	The reason why we consider separated subsequences is that in a separated subsequence, the walk on the original graph $V_{\mz}^{mM}$ needs to take at least one additional step. For a fixed path $P^\prime \subset G^\prime$ of length $k $ and $l=\big \lfloor \frac{k}{3^{d+1}} \big\rfloor$ we have that
	\begin{align*}
		& \p_\beta \left( \left| \left\{ j \in \{0,\ldots,l\}: R_j \text{ not separated} \right\} \right| > \frac{l}{2}  \ \big| \ G^\prime \right)\\
		& = \p_\beta \left( \bigcup_{\substack{U \subset \{0,\ldots, l\}\\ |U|> l/2}} \left\{ R_j \text{ not separated for all } j\in U \right\}   \ \Big| \ G^\prime \right)\\
		&
		\leq \sum_{\substack{U \subset \{0,\ldots, l\}\\ |U|> l/2}} \p_\beta \left(  \left\{ R_j \text{ not separated for all } j\in U \right\}   \ \Big| \ G^\prime \right)
		\leq 2^{l} \left( \frac{C_d\beta^2}{m^{1/2}} \right)^{l/2} \text.
	\end{align*}
	Next, we want to bound the expected degree of vertices in the long-range percolation graph from above. With the bound on the connection probability $\p_\beta(\mz\sim u)$ \eqref{eq:connectionupper bound}, we get that
	\begin{align}\label{eq:expected degree upper bound}
		\notag \E_\beta \left[ \deg(\mz) \right] & 
		= \sum_{u\in \Z^d\setminus\{\mz\}} \p_\beta(\mz\sim u)
		\leq 3^d + \sum_{k=2}^{\infty} \sum_{u \in \cS_k} \frac{2^{2d} \beta }{k^{2d}}
		\leq 3^d + \sum_{k=2}^{\infty} 
		2d (2k+1)^{d-1} \frac{2^{2d} \beta }{k^{2d}}\\
		& \leq 3^d + \beta 2^{3d} 3^d \sum_{k=2}^{\infty} k^{-d-1} \leq 3^d + \beta 3^{4d} \leq \lceil \beta \rceil 3^{5d}.
	\end{align}
	Let $\mathcal{P}_k^\prime$ be the set of self-avoiding paths of length $k$ in $G^\prime$ starting at $r(\mz)$. By a comparison to the case of a Galton-Watson tree, inequality \eqref{eq:expected degree upper bound} already gives that $\E_\beta \left[ | \mathcal{P}_k^\prime | \right]\leq \left( \lceil \beta \rceil 3^{5d} \right)^k$. As $\big\lfloor \frac{k}{3^{d+2}} \big\rfloor \leq \frac{\lfloor \frac{k}{3^{d+1}} \rfloor}{2}$, we see that
	\begin{align*}
		& \p_\beta \left( \exists P^\prime \in  \mathcal{P}_k^\prime \text{ with less than } \big\lfloor \frac{k}{3^{d+2}} \big\rfloor \text{ separated subpaths } R_j \right)\\
		& = \E_{\beta } \left[\p_\beta  \left( \exists P^\prime \in  \mathcal{P}_k^\prime \text{ with less than } \big\lfloor \frac{k}{3^{d+2}} \big\rfloor \text{ separated subpaths } R_j \Big| G^\prime \right) \right]\\
		& \leq  \E_{\beta } \left[ \left|\cP_k^\prime\right| 2^k \left(\frac{C_d\beta^2}{m^{1/2}}\right)^{\big\lfloor \frac{k}{3^{d+2}} \big\rfloor} \right]
		\leq 
		\left(\lceil \beta \rceil 3^{5d}\right)^k 2^k \left(\frac{C_d\beta^2}{m^{1/2}}\right)^{\big\lfloor \frac{k}{3^{d+2}} \big\rfloor}\\
		&
		\leq \left(\lceil \beta \rceil^3 3^{5d}\right)^k 2^k C_d^k \frac{1}{m^{\frac{k}{3^{4d}}}} \leq 0.01^k
	\end{align*}
	by the choice of $m \geq \left(2000 \cdot 3^{5d} C_d \lceil\beta\rceil^3 \right)^{\left(3^{4d}\right)}$. Next, we want to translate this bound on the probability of certain events to bounds on the expectation of the distances. For this, let $\mathcal{G}_{k}$ be the event that all self-avoiding paths $P^\prime \subset G^\prime$ starting at the origin and of length $\tilde{k} \geq k$ contain at least $\big{\lfloor} \frac{\tilde{k}}{3^{d+2}} \big{\rfloor}$ separated subpaths $R_j$. With the preceding inequality we directly get $\p_\beta \left(\mathcal{G}_k\right) \geq 1 - 0.1^k$. On the event $\mathcal{G}_k$, each path $P \subset V_{\mz}^{mM}$ starting at the origin, for which the loop-erased projection on $G^\prime$ goes through $\tilde{k}+1$ different blocks of the form $V_u^m$, needs to have a length of at least $\tilde{k} + \big{\lfloor} \frac{\tilde{k}}{3^{d+2}} \big{\rfloor} \geq \left(1+\frac{1}{3^{d+3}} \right) \tilde{k}$. Furthermore, if we have $D_{G^\prime}\left(r(\mz), r\left((M-1)e_1\right) \right) = \tilde{k}$, then  every path connecting $\mz$ to $(mM-1)e_1$ in the original model $V_\mz^{mM}$ goes through at least $\tilde{k}+1$ different blocks of the form $V_u^m$, with $u\in V_\mz^M$. So we get that
	\begin{align*}
		\E_\beta & \left[D_{V_{\mz}^{mM}} \left(\mz , (mM-1)e_1 \right)\right] 
		\geq 
		\sum_{k=3^{d+3}}^{\infty} \E_\beta\left[ D_{V_{\mz}^{mM}} \left(\mz , (mM-1)e_1 \right) \mathbbm{1}_{\left\{D_{G^\prime}\left(r(\mz), r\left((M-1)e_1\right) \right) = k \right\}} \right]\\
		& \geq \sum_{k=3^{d+3}}^{\infty} \E_\beta\left[ D_{V_{\mz}^{mM}} \left(0 , (mM-1)e_1 \right) \mathbbm{1}_{\left\{D_{G^\prime}\left(r(\mz), r\left((M-1)e_1\right) \right) = k \right\}} \mathbbm{1}_{\mathcal{G}_k} \right]\\
		& \geq \sum_{k=3^{d+3}}^{\infty} \E_\beta\left[ k \left(1+\frac{1}{3^{d+3}}\right) \mathbbm{1}_{\left\{D_{G^\prime}\left(r(\mz), r\left((M-1)e_1\right) \right) = k \right\}} \mathbbm{1}_{\mathcal{G}_k} \right]\\
		& \geq \sum_{k=3^{d+3}}^{\infty} \left( \E_\beta\left[ k \left(1+\frac{1}{3^{d+3}}\right) \mathbbm{1}_{\left\{D_{G^\prime}\left(r(\mz), r\left((M-1)e_1\right) \right) = k \right\}}  \right]
		-
		\E_\beta\left[ 2k \mathbbm{1}_{\mathcal{G}_k^C} \right]
		\right)\\
		& \geq \left(1+\frac{1}{3^{d+3}}\right)
		\sum_{k=3^{d+3}}^{\infty}
		k
		 \E_\beta\left[  \mathbbm{1}_{\left\{D_{V_{\mz}^M}\left(\mz,  (M-1)e_1 \right) = k \right\}}  \right]
		-
		2 \sum_{k=1}^{\infty} 0.1^k k\\
		& \geq \left(1+\frac{1}{3^{d+3}}\right) \E_\beta \left[ D_{V_{\mz}^M}\left(\mz,  (M-1)e_1 \right) \right] - 3^{d+4} \p_\beta \left( D_{V_{\mz}^M}\left(\mz,  (M-1)e_1 \right) < 3^{d+3} \right) - 1\\
		& \geq \left(1+\frac{1}{3^{d+4}}\right) \E_\beta \left[ D_{V_{\mz}^M}\left(\mz,  (M-1)e_1 \right) \right]
	\end{align*}
	where the last inequality holds for all large enough $M$, as the probability of the event $\left\{  D_{V_{\mz}^M}\left(\mz,  (M-1)e_1 \right) < 3^{d+3} \right\}$ tends to $0$ as $M \to \infty$. Say that it holds for all $M \geq m^N$, where $m= \Big\lceil \left(2000 \cdot 3^{5d} C_d \lceil\beta\rceil^3 \right)^{\left(3^{4d}\right)} \Big\rceil$. The important property about the choice of $m$ is, that its size is polynomial in $\beta$. This already implies that 
	\begin{align*}
		\dxp(\beta) &\geq  \lim_{n \to \infty} \frac{\log \left( \E_\beta \left[ D_{V_{\mz}^{m^n}}\left(\mz,  (m^n-1)e_1 \right) \right]\right)}{\log\left(m^n\right)}
		\geq
		\lim_{n \to \infty} \frac{\log \left(  \left(1+\frac{1}{3^{d+4}}\right)^{n-N} \right) }{\log\left(m^n\right)}\\
		& = \frac{\log\left(1+\frac{1}{3^{d+4}} \right)}{\log(m)} \geq \frac{c}{\log(\beta)}
	\end{align*}
	for some small $c>0$ and all $\beta \geq 2$.
\end{proof}

\section{Connected sets in graphs}

The expected number of open paths in the long-range percolation model, of length $k$, and starting at $\mz$, grows at most like $\E\left[\deg(\mz)\right]^k$, which can be easily proven by a comparison with a Galton-Watson tree. However, it is a priori not clear how the number of connected subsets of $\Z^d$ containing the origin grows. In particular, because the maximal degree of vertices is unbounded. In this chapter, we prove several results about the structure of connected sets in the long-range percolation graph. Mostly, we want to prove that with exponentially high probability in $k$, all connected sets of size $k$ in the graph have not too many edges. First, we need to define what we mean by a connected set. Formally, we define the a connected set as follows. For a graph $G=(V,E)$ we say that a subset $Z \subset V$ is {\sl connected} if the graph $(Z , E^\prime)$ with edge set $E^\prime = \left\{  \{x,y\} \in E : x,y \in Z \right\}$ is connected. As a first step, we bound the expected number of connected sets of certain size in Galton-Watson trees. This counting of connected sets plays an important role in Section \ref{sec:prooftheo1} below and in the companion paper \cite{baeumler2022behaviour}. 

\begin{lemma}\label{lem:galtonwatson_consets}
	Let $\cX$ be a countable set with a total ordering and a minimal element, let $X$ be a countable sum of independent Bernoulli-distributed random variables over this set, i.e.,  $X=\sum_{i\in \cX} X_i$, and let $\mu$ be the expectation value of $X$. Say that $q(k)=\p(X_k=1)$. Let $T$ be a Galton-Watson tree with offspring distribution $\mathcal{L}(X)$. We denote the set of all subtrees of $T$ of size $k$ containing the origin by $\mathcal{T}_k$. Then
	\begin{align*}
	\E \left[|\mathcal{T}_k|\right] \leq 4^k \mu^k  .
	\end{align*}
\end{lemma}

\begin{proof}
	The choice of the set $\cX$ and the total ordering on it do not influence the outcome, so we will always work with $\cX=\N$ from here on.
	We can think of the Galton-Watson tree as a model of independent bond percolation on the graph with vertex set $L=\bigcup_{n=0}^\infty L_n$, where $L_n = \N^n$, and with edge set $S=\left\{\{v,(v \ m)\}: v \in L, m\in \N \right\}$ where some edge of the form $\{v,(v \ m)\}$ is open with probability $q(m)$. Note that the graph $G=(L,S)$ is a tree, so in particular there exists a unique path from the origin $\emptyset$ to every vertex; this tree is also known as the Ulam-Harris tree. For a vertex $v \in L$, the number of open edges of the form $\{v,(v \ m)\}$ has the same law as $X$ and thus we can identify the open cluster connected to the root $\emptyset$ with a Galton-Watson tree with offspring distribution $\mathcal{L}(X)$. So in particular, the expected number of subtrees of the Galton-Watson tree $T$ of size $k$ is the same as the expected number of connected sets of size $k$ in $(L,S)$.  For a vertex $v \in L$, the number of open edges of the form $\{v,(v \ m)\}$ has the same law as $X$ and thus we can identify the open cluster connected to the root $\emptyset$ with a Galton-Watson tree with offspring distribution $\mathcal{L}(X)$. For a vertex $v\in L$, we call the vertices of the form $(v \ m)$ that are connected to $v$ by an open bond the {\sl children} of $v$. Vice versa, we say that $v$ is the {\sl parent} of the vertex $(v \ m )$, if $(v \ m)$ is connected to $v$. For a connected set $L^\prime \subset L$ of size $k$, we now describe an exploration process $\left(Y_i\right)_{i\in\{1,\ldots,2k-1\}}$ of it:
	\begin{enumerate}\addtocounter{enumi}{-1}
		\item Start with $Y_1 = \emptyset$.
		\item For $i=1 , \ldots , 2k-1$ 
		\begin{enumerate}
			\item If there exists $m\in \N$ for which $\left(Y_i \ m\right) \in L^\prime$ and $Y_j \neq (Y_i \ m)$ for all $j< i$, let $m^\prime$ be the minimal among those $m\in \N$ and set $Y_{i+1} = (Y_i \ m^\prime)$.
			\item If such an $m$ does not exist, let $Y_{i+1}$ be the parent of $Y_i$.
		\end{enumerate}
	\end{enumerate}

	\begin{figure}
		\[\begin{tikzpicture}[scale = 1.5]
		\vertex[label = $Y_1 \text{,}  Y_7 \text{,} Y_{15}$ ] (01) at (0,3) {$Y_1^\prime$};
		
		\vertex[label = $Y_2 \text{,} Y_4  \text{,} Y_6   \ \ \ \  \text{ }$] (11) at (-1.5,2) {$Y_2^\prime$};
		\vertex[label = $ \ \ \ \ Y_{8} \text{,} Y_{12}  \text{,} Y_{14}$] (12) at (1.5,2) {$Y_5^\prime$};
		
		\vertex[label = $Y_3$] (21) at (-2.5,1) {$Y_3^\prime$};
		\vertex[label = $ \ Y_5$] (22) at (-1.2,1) {$Y_4^\prime$};
		\vertex[label = $Y_{9} \text{,} Y_{11} \ \  \ \ \ \ $] (24) at (1,1) {$Y_6^\prime$};
		\vertex[label = $\ \ Y_{13}$] (25) at (2,1) {$Y_8^\prime$};
		
		\vertex[label = $\ \ \ Y_{10}$] (31) at (1.3,0) {$Y_7^\prime$};
		
		\path[thick]
		(01) edge (11) (01) edge (12)
		(11) edge (21) (11) edge (22)
		(12) edge (24) (12) edge (25)
		(24) edge (31);
		
		\end{tikzpicture}\]
		\centering
		\parbox{11cm}{\caption{In the above tree, the process $\left(Y_i^\prime\right)_{i \in \{1,\ldots,9\}}$ is written inside the vertices and the process $\left(Y_i\right)_{i \in \{1,\ldots,15\}}$ is written above the vertices. For this tree we have $(a_1,\ldots,a_{14})=(d,d,u,d,u,u,d,d,d,u,u,d,u,u)$.} \label{fig:tree}}

	\end{figure}
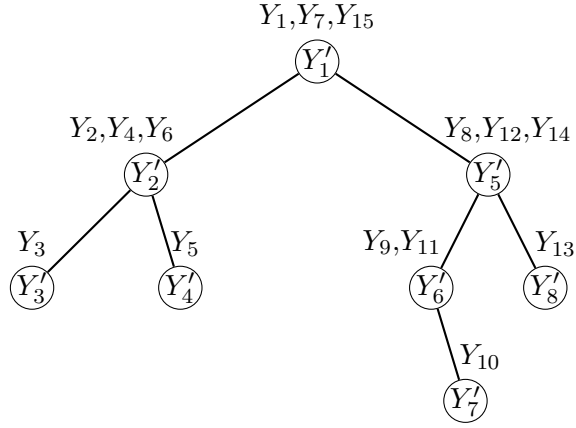

	An example of this procedure is given in Figure \ref{fig:tree}. This exploration process traverses every edge exactly twice in opposite directions and starts and ends at the origin of the tree. We also say that the exploration process $Y_i$ goes (one level) down if $(a)$ occurs in the algorithm above and otherwise we say that the process goes (one level) up. We also define a different process $\left(Y^\prime_i\right)_{i\in\{1,\ldots,k\}}$, where $Y^\prime_i$ is the unique point $Y_l$ such that $\left|\{Y_1,\ldots, Y_{l-1}\}\right|< i$ and $\left|\{Y_1,\ldots, Y_{l}\}\right| = i$. So the process $\left(Y^\prime_i\right)_{i\in\{1,\ldots,k\}}$ is like a depth-first search from left to right in the tree.
	We can encode all information contained in the subtree $L^\prime$ by the two sequences $\left(a_1,\ldots,a_{2k-2}\right)\in \{u,d\}^{2k-2}$ and $\left(m_1,\ldots, m_{k-1}\right)\in \N^{k-1}$. The first sequence $\left(a_1,\ldots,a_{2k-2}\right)$ encodes whether the process $Y_i$ goes one level up or down at a certain point. Here $a_i=u$ if the process goes one level up after $Y_i$, i.e.,  if $Y_{i+1}$ is the parent of $Y_i$. Otherwise we set $a_i=d$, i.e.,  if $Y_{i+1}$ is a child of $Y_i$.
	The sequence $\left(m_1,\ldots, m_{k-1}\right)$ encodes the direction of the process, where the $i$-th coordinate gives the direction when the walk goes down for the $i$-th time. This happens when it touches the vertex $Y_{i+1}^\prime$ for the first time. So if $v$ is the parent of $Y_{i+1}^\prime$, then $Y_{i+1}^\prime = \left(v \ m_i\right)$. \\
	
	For fixed $\overset{\to}{a} = \left(a_1,\ldots,a_{2k-2}\right) \in \{u,d\}^{2k-2}$, we want to upper bound the expected number of subtrees containing the origin with exactly this up-and-down structure. Assume that the exploration process $Y_i$ visits exactly $l$ children of some vertex $Y_j^\prime$. Then the expected number of ways to choose these $l$ children among the children of $Y_j^\prime$ in an increasing way is given by 
	\begin{align*}
	\sum_{m_1 \in \N  } q( m_1)  \sum_{\substack{m_2 \in \N: \\ m_2 > m_1} } q( m_2) \ \dots \sum_{\substack{m_l \in \N: \\ m_l > m_{l-1}} } q( m_l) \leq \mu^l   .
	\end{align*}
	We have this choice for all vertices $Y_j^\prime$ in the tree. The sum over the number of children of all the vertices is $k-1$, as every vertex, except the origin $\emptyset$, is the child of exactly one vertex.
	Thus the expected number of trees with a specified up-and-down structure can be bounded from above by 
	\begin{align*}
	\sum_{m_1 \in \N} \dots \sum_{m_{k-1} \in \N} \prod_{i=1}^{k-1} q( m_i) = \mu^{k-1}  .
	\end{align*}
	Up to now, we considered a fixed up-and-down-structure. However, there are at most $\left|\{u,d\}^{2k-2}\right|= 2^{2k-2}$ possible up-and-down structures $(a_1,\ldots,a_{2k-2})$ (In fact there are significantly less combinations, as one has additional constraints like $a_1=d$). So in total, we get
	\begin{align*}
	\E\left[|\mathcal{T}_k|\right] \leq \sum_{\overset{\to}{a} \in \{u,d\}^{2k-2}} \mu^{k-1} \leq \left(2^{2k-2}\right) \mu^{k-1} \leq 4^k \mu^k.
	\end{align*} 
\end{proof}
We now want to use the above lemma about the Galton-Watson tree in order to get results about the average degree of connected subsets of the long-range percolation graph. For this, we define the average degree of some set finite $Z\subset\Z^d$ by
\begin{align*}
\overline{\deg}(Z) \coloneqq \frac{1}{|Z|} \sum_{v\in Z} \deg(v) \text .
\end{align*}
One elementary inequality we will use in the following controls the exponential moments of certain random variables. Assume that $\left(U_{i}\right)_{i\in \N}$ are independent Bernoulli random variables and $U= \sum_{i=1}^{\infty} U_i$. Then
\begin{align}\label{eq:exp moment bernoulli}
\E\left[e^{U}\right] & = \E\left[e^{\sum_{i\in \N} U_i}\right] = \prod_{i\in \N} \E\left[e^{U_i}\right]  \leq \prod_{i\in \N} \left(1+e\E\left[U_i\right]\right) \leq \prod_{i\in \N} e^{e\E\left[U_i\right]} = e^{e\E\left[U\right]} \ 
\end{align}
and this already implies, by Markov's inequality, that for any $C>0$
\begin{align}\label{eq:exp moment bernoulli prob bound}
\p\left( U > C \E \left[U\right] \right) = \p\left( e^U > e^{C \E \left[U\right]} \right) \leq \E \left[ e^U\right] e^{-C \E \left[U\right]} \overset{\eqref{eq:exp moment bernoulli}}{\leq} e^{(e-C)\E\left[U\right]} .
\end{align}

\begin{lemma}\label{lem:connnectedsetsinLRP}
	Let $\mathcal{CS}_k= \mathcal{CS}_k\left(\Z^d\right)$ be all connected subsets of the long-range percolation graph with vertex set $\Z^d$, which are of size $k$ and contain the origin $\mz$. We write $\mu_\beta$ for $\E_\beta \left[\deg(\mz)\right]$. Then for all $\beta>0$ 
	\begin{align*}
	\p_\beta \left(\exists Z \in \mathcal{CS}_k : \overline{\deg}(Z) \geq 20 \mu_\beta\right) \leq e^{-4k\mu_\beta}.
	\end{align*}
\end{lemma}

\begin{proof}
	Consider percolation on the tree $L = \bigcup_{n=0}^\infty L_n$, where $L_n = \left(\Z^d \setminus \{\mz\}\right)^n$, the edge set is given by $S=\left\{\{v,(v \ m)\}: v \in L, m\in \Z^d \setminus \{\mz\} \right\}$ and an edge of the form $\{v,(v \ m)\}$ is open with probability $p\left(\beta, \{\mz,m\}\right)$. A total ordering on $\Z^d \setminus \{\mz\}$ is given by considering an arbitrary deterministic bijection with $\N$. From Lemma \ref{lem:galtonwatson_consets}, we know that the expected number of  connected sets of size $k$ in $L$ is bounded by $4^k \mu_\beta^k$.	We want to project a finite tree $T\subset L$ of size $k$ down to $\Z^d$. Remember the notation $\left(Y_i^\prime\right)_{i\in\{1,\ldots,k\}}$ for the depth-first search from left to right in the tree. The information contained in the structure of the tree can be represented by the vectors $\overset{\to}{a}=(a_1,\ldots,a_{2k-2}) \in \{u,d\}^{2k-2}$ and  $\overset{\to}{m}=(m_1,\ldots,m_{k-1}) \in \left(\Z^d\setminus\{\mz\}\right)^{k-1}$. We now define a  subgraph $\left(Z(T),E(T)\right)$ of the integer lattice and an exploration process $\left(X^\prime_i\right)_{i \in \{1,\ldots, k\}}$ as follows:
	\begin{enumerate}\addtocounter{enumi}{-1}
		\item Start with $X_1^\prime = \mz, E_1(T)= \emptyset$.
		\item For $i=2,\ldots,k:$ \\
		Let $j<i$ be such that $Y_i^\prime = (Y_j^\prime \  m)$ for some $m \in \Z^d\setminus\{\mz\}$. Set $X_i^\prime = X_j^\prime + m$ and $E_i(T) = E_{i-1}(T) \cup \left\{\{X_i^\prime, X_j^\prime\}\right\}$.
		\item $Z(T)= \bigcup_{i=1}^k \{X_i^\prime\}$ and $E(T)=E_{k}(T)$.
	\end{enumerate}
	
		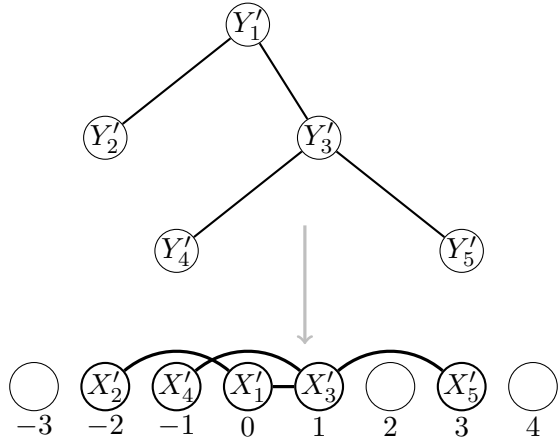
\begin{figure}
		
		\centering
		\begin{tikzpicture}[scale = 1.5]
		\vertex[ ] (01) at (0/1.6,3) {$Y_1^\prime$};
		\vertex[ ] (11) at (-2/1.6,2) {$Y_2^\prime$};
		\vertex[ ] (12) at (1/1.6,2) {$Y_3^\prime$};
		\vertex[ ] (21) at (-1/1.6,1) {$Y_4^\prime$};
		\vertex[ ] (22) at (3/1.6,1) {$Y_5^\prime$};
		
		\path[thick]
		(01) edge (11) (01) edge (12) (12) edge (21) (12) edge (22);
		
		\vertex[ ] (-3) at (-3/1.6,-0.2) {$\ \ \text{  }\ \text{  } $   };
		\vertex[ thick ] (-2) at (-2/1.6,-0.2) {$X_2^\prime$};
		\vertex[ thick ] (-1) at (-1/1.6,-0.2) {$X_4^\prime$};	
		\vertex[ thick ] (0) at (0/1.6,-0.2) {$X_1^\prime$};
		\vertex[ thick ] (1) at (1/1.6,-0.2) {$X_3^\prime$};
		\vertex[ ] (2) at (2/1.6,-0.2) {$\ \ \text{  }\ \text{  } $   };
		\vertex[ thick ] (3) at (3/1.6,-0.2) {$X_5^\prime$};
		\vertex[ ] (4) at (4/1.6,-0.2) {$\ \ \text{  }\ \text{  } $   };
		
		\vertex[ draw =none] (-3a) at (-3/1.6,-0.55) {$-3$   };
		\vertex[ draw=none ] (-2a) at (-2/1.6,-0.55) {$-2$};
		\vertex[draw=none ] (-1a) at (-1/1.6,-0.55) {$-1$};	
		\vertex[draw=none ] (0a) at (0/1.6,-0.55) {$0$};
		\vertex[draw=none ] (1a) at (1/1.6,-0.55) {$1$};
		\vertex[draw=none ] (2a) at (2/1.6,-0.55) {$2$   };
		\vertex[draw=none ] (3a) at (3/1.6,-0.55) {$3$};
		\vertex[draw=none ] (4a) at (4/1.6,-0.55) {$4$   };
		
		\path[very thick]
		(-2) edge[bend left = 40] (0) (-1) edge[bend left= 40] (1) (0) edge (1) (1) edge[bend left = 40] (3);

		\vertex[draw =none ] (a) at (0.5,1.3) {   };
		\vertex[draw =none ] (b) at (0.5,0.1) {   };
		
		\path[very thick]
		(a) edge[->, color = lightgray] (b);

		\end{tikzpicture}
		
		\parbox{11cm}{\vspace{0.4cm}{\caption{A tree $T$ with 5 vertices, $(a_1,\ldots,a_{8})=(d,u,d,d,u,d,u,u)$, and $(m_1,\ldots,m_4)=(-2,1,-2,2)$, and its projection on $\Z$. The vertices with thick boundary $\{-2,-1,0,1,3\} \subset \Z$ are the set $Z(T)$ and the thick edges between them are the set $E(T)$. Note that $(Z(T),E(T))$ really is a tree for this example.} \label{fig:explo} } }

	\end{figure}
	
	See Figure \ref{fig:explo} for an example of this projection. The graph $\left(Z(T),E(T)\right)$ is clearly connected, but it is not necessarily a tree, as there can be $i\neq j$ with $X_i^\prime = X_j^\prime$, in which case there exists a loop containing $X_i^\prime$. We call both the graph $\left(Z(T), E(T)\right)$ and the tree $T$ {\sl admissible} if $\left(Z(T),E(T)\right)$ is a tree. We also write $\mathcal{TA}_k$ for the set of admissible trees $T\subset(L,S)$ of size $k$. For a tree $T\subset (L, S)$ of size $k$, the condition $T \in \mathcal{TA}_k$ is equivalent to $|Z(T)|=k$, as every connected graph with $k$ vertices and $k-1$ edges is a tree. Assume that the graph $\left(Z(T),E(T)\right)$ is admissible. Then the probability that all edges exist in the random graph equals $\prod_{i=1}^{k-1}p(\beta, \{\mz,m_i\})$, which is exactly the probability that all edges of the tree $T$ exist inside $(L,S)$. Every connected set $Z \subset \Z^d$ has a spanning tree. Thus there exists a tree $T\subset L$ with $Z=Z(T)$ such that all edges in $E(T)$ exist. This and the result of Lemma \ref{lem:galtonwatson_consets} imply that
	\begin{align}\label{eq:connectedsetbound}
	\notag \E_\beta \left[ \left| \mathcal{CS}_k(\Z^d) \right| \right] &
	\leq \sum_{T \in \mathcal{TA}_k} \p_{\beta }\left( \text{all edges in } E(T) \text{ exist} \right) 
	= \sum_{T \in \mathcal{TA}_k} \p_{\beta }\left( T \in \mathcal{T}_k \right) \\
	& \leq \E_{\beta }\left[ \left| \mathcal{T}_k \right| \right] \leq 4^k \mu_\beta^k.
	\end{align}
	For an admissible tree $T$, the degree of each vertex $v\in Z(T)$ is the sum of an {\sl inside degree} and an {\sl outside degree}, which we will now define. The {\sl inside degree} $\deg_{Z(T)}(v)$  of a vertex $v \in Z(T)$ is defined by
	\begin{equation*}
	\deg_{Z(T)}(v) = \sum_{u\in Z(T)} \mathbbm{1}_{\left\{\{v,u\} \in E(T)\right\}}
	\end{equation*}
	which is just the number of edges in $E(T)$ containing $v$. Note that for a given admissible tree $T$, the inside degree is purely deterministic and does not depend on the environment. Also note that, by the handshaking lemma,
	\begin{align}\label{eq:inside_handshake}
	\sum_{v \in Z(T)}\deg_{Z(T)}(v) = 2 |E(T)| = 2\left(|Z(T)|-1\right) \text ,
	\end{align}
	where the last equality holds as $(Z(T),E(T))$ is a tree.
	Now let us turn to the {\sl outside degree} $\deg_{Z(T)^C}(v)$ of a vertex $v\in Z(T)$, which we define by
	\begin{equation*}
	\deg_{Z(T)^C}(v) = \sum_{\substack{u \in \Z^d \setminus \{v\}: \\ \{u,v\} \notin E(T)}} \omega(\{v,u\}) \text .
	\end{equation*}
	The outside degree depends on the random environment $\omega$ and is a non-constant random variable, contrary to $\deg_{Z(T)}(v)$. 
	Now we want to get bounds on the random variable $\sum_{v \in Z(T)} \deg_{Z(T)^C}(v)$. Note that $\{u,v\} \notin E(T)$ does not imply that $u\notin Z(T)$, but only that $u$ and $v$ are not neighbors in the graph induced by $T$. The random variable $\sum_{v \in Z(T)} \deg_{Z(T)^C}(v)$ is not the sum of independent Bernoulli random variables, as we might count some edges twice. But as one can count every edge at most twice in this sum, one has the bound
	\begin{align}\label{eq:leavesbound3}
	\frac{1}{2} \sum_{v \in Z(T)} \deg_{Z(T)^C}(v) \leq  \sum_{ \substack{\{u,v\} \notin E(T) : \\  \{u,v\} \cap Z(T) \neq \emptyset } } \omega\left(\{u,v\}\right)
	\end{align}
	where the expression on the right-hand side is a sum of independent Bernoulli random variables with expectation at most $|Z(T)|\mu_\beta$. So for each admissible tree $T$ we always have
	\begin{align*}
	\sum_{v\in Z(T)} \deg(v) 
	& =
	\sum_{v\in Z(T)} \deg_{Z(T)}(v) + \sum_{v\in Z(T)} \deg_{Z(T)^C}(v) 
	\leq 2 |Z(T)| + 2  \sum_{ \substack{\{u,v\} \notin E(T) : \\  \{u,v\} \cap Z(T) \neq \emptyset } } \omega\left(\{u,v\}\right) \text .
	\end{align*}
	We use the notation 
	\begin{align*}
	U = U(T) \coloneqq \sum_{ \substack{\{u,v\} \notin E(T) : \\  \{u,v\} \cap Z(T) \neq \emptyset } } \omega\left(\{u,v\}\right) \text .
	\end{align*}
	For a given finite admissible tree $T$, we have that
	\begin{align}\label{eq:deviationbound}
	\notag \p_\beta & \left( \overline{\deg}(Z(T)) \geq 20 \mu_\beta\right) = 
	\p_\beta \left( \sum_{v\in Z(T)} \deg(v)  \geq  20 |Z(T)| \mu_\beta\right)
	\leq \p_\beta\left( 2  U \geq 18|T|\mu_\beta\right)\\
	& = \p_\beta\left(  U \geq 9 |T|\mu_\beta\right) \overset{\eqref{eq:exp moment bernoulli prob bound}}{\leq}
	\E\left[e^U\right] e^{-9|T|\mu_\beta} \leq e^{e|T| \mu_\beta}  e^{-9|T|\mu_\beta} \leq  e^{-6|T|\mu_\beta}.
	\end{align}
	So far we only got this bound for a fixed admissible tree $T\subset (L,S)$. Remember that every connected set $Z\in \mathcal{CS}_k$ has a spanning tree and there exists a tree $T\subset (L,S)$ so that $\left(Z(T),E(T)\right)$ is exactly this spanning tree. Again, we use the notation $\mathcal{TA}_k$ for the set of admissible trees $T\subset (L,S)$ of size $k$. With the observation from before we get that
	\begin{align*}
	\p_\beta & \left(\exists Z \in \mathcal{CS}_k : \overline{\deg}(Z) \geq 20 \mu_\beta\right) 
	\leq \sum_{T \in \mathcal{TA}_k} \p_\beta\left(  \overline{\deg}(Z(T)) \geq 20 \mu_\beta , \text{ all edges in $E(T)$ exist}  \right)\\
	& \leq \sum_{T \in \mathcal{TA}_k} \p_\beta\left(U(T) \geq 9 k \mu_\beta, \text{ all edges in $E(T)$ exist}  \right)\\
	& = \sum_{T \in \mathcal{TA}_k} \p_\beta\left( U(T) \geq 9k\mu_\beta \right)  \p_\beta \left( \text{all edges in $E(T)$ exist}  \right)\\
	& \overset{ \eqref{eq:deviationbound}}{\leq } e^{-6 k \mu_\beta}  \sum_{T \in \mathcal{TA}_k}  \p_\beta \left( \text{all edges in $E(T)$ exist}  \right)
	\\
	& \overset{\eqref{eq:connectedsetbound}}{\leq} 
	e^{-6k\mu_\beta} 4^k \mu_\beta^k
	\leq  e^{-6k\mu_\beta} e^{2k} e^{\mu_\beta k} \leq e^{-4k\mu_\beta}
	\end{align*}
	where we used that $\mu_\beta \geq 2$ in the last inequality. This holds for long-range percolation with our parameters, as each vertex is always connected to its nearest neighbors. The final inequality is exactly the result that we wanted to show and thus we finish the proof.
\end{proof}

\section{Distances in $V_{\mz}^n$}

In this section, we give several bounds on the distribution of the graph distances between points, respectively sets, inside of certain boxes. In Section \ref{sec:graph distances of far away}, we determine several different properties of the function $(x,y)\mapsto \E_\beta \left[ D_{V_\mz^n} (x,y) \right]$. It is intuitively clear that the expression is large when $x,y$ also have a big Euclidean distance, for example when $x=\mz$ and $y=(n-1)\mo$. This intuition is made rigorous in Lemma \ref{lem:endisfar}. In Section \ref{sec:2nd moment}, we upper bound the second moment of random variables of the form $D_{V_\mz^n} (x,y)$. Then, in Section \ref{sec:distances between points and boxes} we use these results in order to bound the distance between certain points and sets in the long-range percolation graph.

\subsection{Graph distances of far away points}\label{sec:graph distances of far away}

From the definition of $\Lambda(n,\beta)$ in \eqref{eq:Lambda} it is not clear which pair $u,v$ maximizes the expected distance and how the expected graph distances can be compared for different graphs $V_\mz^n$ and $V_\mz^{n^\prime}$. In Proposition \ref{propo:scaling}, we construct a coupling between the long-range percolation graph on $V_{\mz}^n$ for different $n$. In Lemma \ref{lem:endisfar}, we show that, up to a constant factor, the maximum in the definition of $\Lambda(n,\beta)$ gets attained by the pair $\{\mz, (n-1)e_1\}$ or $\{\mz,(n-1)\mo\}$.

\begin{proposition}\label{propo:scaling}
	Let $\beta\geq 0$ and $n^\prime, n \in \N_{>0}$ with $n^\prime \leq n$. For $u,v \in V_{\mz}^n$ define $u^\prime \coloneqq \lfloor \frac{n^\prime }{n} u \rfloor, v^\prime \coloneqq \lfloor \frac{n^\prime }{n} v \rfloor$, where the rounding operation is componentwise. There exists a coupling of the random graphs with vertex sets $V_{\mz}^n$ and $V_{\mz}^{n^\prime}$ such that both are distributed according to $\p_\beta$ and
	\begin{equation}\label{eq:scaling}
		 D_{V_{\mz}^{n^\prime}}(u^\prime ,v^\prime) \leq  3  D_{V_{\mz}^n}(u, v)
	\end{equation}
	for all $u,v \in V_{\mz}^n$. The same holds true when one considers the graph $\Z^d$ instead of $V_{\mz}^n$ and this also implies that
	\begin{equation}\label{eq:scaling dia}
		\dia\big(V_0^{n^\prime}\big) \leq 3 \dia\big(V_0^{n}\big) \text .
	\end{equation}
\end{proposition}

\begin{proof}
	We prove the statement via a coupling with the underlying continuous model. As the claim clearly holds for $\beta=0$ or for $u=v$, we can assume $\beta>0 $, and $ u \neq v$ from here on. Let $\tilde{\mathcal{E}}$ be a Poisson point process on $\R^d\times \R^d$ with intensity $\frac{\beta}{2\|t-s\|^{2d}}$ and define $\e = \left\{(t,s)\in \R^d \times \R^d : (s,t)\in \tilde{\e}\right\} \cup \tilde{\e}$. Remember that this point process has a scaling invariance, namely that for a constant $\alpha >0$ the set $\alpha\e$ has exactly the same distribution as $\e$. We now define a random graph $G=(V,E)$: For $u,v \in V_{\mz}^n \eqqcolon V$ we place an edge between $u$ and $v$ if and only if $(u+\cC) \times (v + \cC) \cap n \e \neq \emptyset$. We have already seen in Section \ref{subsec:cts} about the continuous model that this creates a sample of independent long-range percolation where the connection probability between the vertices $u$ and $v$ is given by $p(\beta,|v-u|)=1-e^{-\int_{u+\cC} \int_{v+ \cC} \frac{\beta}{\|t-s\|^{2d}}\md t \md s}$. 
	We can do the same procedure for $V^\prime \coloneqq V_{\mz}^{n^\prime}$ and $n^\prime \e$ to get a random graph $G^\prime=(V^\prime, E^\prime)$. Formally, we place an edge between two vertices $u^\prime, v^\prime \in V^\prime$  if and only if $(u^\prime+\cC) \times (v^\prime + \cC) \cap n^\prime \e \neq \emptyset$. We now claim that for any two vertices $u,v \in V$ with $u\neq v$ and $u^\prime, v^\prime$ defined as above one has $D_{G^\prime}(u^\prime , v^\prime) \leq 2 D_{G}(u,v)+1$, which already implies \eqref{eq:scaling}. Assume that $(x_0=u,x_1,\ldots,x_l=v)$ is the shortest path between $u$ and $v$ in $G$, where $l=D_G(u,v)$. Then for all $i=1,\ldots,l$ there are points
	\begin{align*}
		\left(y(i,0), y(i,1)\right) \in (x_{i-1} + \cC) \times (x_i + \cC) \cap n \e \text .
	\end{align*}
	In particular one has 
	\begin{equation*}
	\|y(i-1,1) - y(i,0)\|_\infty< 1
	\end{equation*}
	for all $i=2,\ldots,l$, $\|y(1,0)-u\|_\infty < 1$, and $\|y(l,1) - v\|_\infty < 1$. For all $i=1,\ldots, l$ and $j \in \{0,1\}$ define $y^\prime(i,j)  = \frac{n^\prime}{n} y(i,j)$, which implies $\left( y^\prime(i,0) , y^\prime(i,1) \right) \in n^\prime \mathcal{E}$. With this definition one clearly has 
	\begin{equation*}
	\|y^\prime(i-1, 1)  - y^\prime(i, 0) \|_\infty < 1
	\end{equation*}
	for all $i=2,\ldots,l$, $\| y^\prime(1,0) - \frac{n^\prime}{n} u \|_\infty<1$, and $\| y^\prime (l,1) - \frac{n^\prime}{n} v \|_\infty <1$. So in $G^\prime$ we can use the path from $u^\prime$ to $v^\prime $ that uses all the edges $\left\{\lfloor y^\prime (i,0) \rfloor, \lfloor y^\prime (i,1) \rfloor \right\}$ and in the case where $\lfloor y^\prime (i-1,1) \rfloor \neq \lfloor y^\prime (i,0) \rfloor$ holds, respectively the analogous statement for $u^\prime$ or $v^\prime$ holds, we can use the nearest neighbor edge between those vertices, which exists as $\|y^\prime (i-1,1) - y^\prime (i,0) \|_\infty < 1$. So for each vertex that is touched by the shortest path between $u$ and $v$ in $G$ one needs to make at most one additional step for the path between $u^\prime$ and $v^\prime$ in $G^\prime$, which implies that $D_{G^\prime}(u^\prime,v^\prime) \leq 2D_G(u,v) + 1$. If one does not restrict to the sets $V=V_{\mz}^n$ and $V^\prime = V_{\mz}^{n^\prime}$, but works on the graph with vertex set $\Z^d$ instead, the same proof works.
\end{proof}

\begin{lemma}\label{lem:endisfar}
	For all $\beta\geq 0$, $n \in \N_{>0}$, and $u,v \in V_{\mz}^n$, we have
	\begin{equation}\label{eq:endisfar}
	\E_\beta\left[D_{V_{\mz}^{n}}(u,v)\right] \leq 6d \E_\beta\left[D_{V_{\mz}^{n}} (\mz,(n-1) e_1)\right] 
	\end{equation}
	and 
	\begin{equation}\label{eq:endisfar2}
	\E_\beta\left[D_{V_{\mz}^{n}}(\mz, (n-1)e_1)\right] \leq 6 \E_\beta \left[D_{V_{\mz}^{n}} (\mz,(n-1) \mo) \right] \text .
	\end{equation}
\end{lemma}
This lemma already has two interesting implications, that we want to discuss before going to the proof.

\begin{remark}\label{remark:Lambda}
	Combining \eqref{eq:endisfar} and \eqref{eq:endisfar2} already implies that for $\Lambda(n,\beta) = \max_{u,v \in V_{\mz}^{n}} \E_\beta\left[D_{V_\mz^n} (u,v) \right]+1$ one has
	\begin{align*}
		& \E_\beta\left[D_{V_{\mz}^{n}} (\mz,(n-1) e_1)\right] + 1 \leq \Lambda(n,\beta) \leq 6d \E_\beta\left[D_{V_{\mz}^{n}} (\mz,(n-1) e_1)\right] + 1 \ \text{ and } \\
		& \E_\beta \left[D_{V_{\mz}^{n}} (\mz,(n-1) \mo) \right] + 1 \leq \Lambda(n,\beta) 
		\leq 36 d \E_\beta \left[D_{V_{\mz}^{n}} (\mz,(n-1) \mo) \right] + 1 \text .
	\end{align*}
\end{remark}

\begin{remark}\label{remark:Lambda grows}
	For all bounded sets $K\subset \R_{\geq 0}$ there exists a constant $\dxp^\star >0$ such that for all $\beta \in K$ and all $M,N$ large enough one has
	\begin{align*}
		\Lambda(MN,\beta)\geq M^{\dxp^\star} \Lambda(N,\beta) \text .
	\end{align*}
\end{remark}

\begin{proof}
	\Cref{remark:Lambda} together with \eqref{eq:lowerboundtoshow1} already show the existence of such an $\dxp^\star$ along a subsequence of numbers of the form $M=m^k$. Proposition \ref{propo:scaling} shows the result for all large enough $M$.
\end{proof}

\begin{proof}[Proof of Lemma \ref{lem:endisfar}]
	Using the triangle inequality and linearity of expectation we get for all $u,v \in V_\mz^n$ that
	\begin{align*}
	\E_\beta \left[ D_{V_{\mz}^n}(u,v) \right] \leq \E_\beta \left[ D_{V_{\mz}^n}(u,\mz) \right] + \E_\beta \left[ D_{V_{\mz}^n}(\mz,v) \right] 
	\end{align*}
	and thus, in order to prove \eqref{eq:endisfar}, it suffices to show that 
	\begin{align}\label{eq:endisfar reduction}
		\E_\beta \left[ D_{V_{\mz}^n}(\mz,v) \right]  \leq 3d \E_\beta\left[D_{V_{\mz}^{n}} (\mz,(n-1) e_1)\right] 
	\end{align}
	for all $v \in V_{\mz}^n$. By symmetry, we can assume that $p_1(v) \geq p_2(v)\geq \ldots \geq p_d(v)$. For $k\in \{0,\ldots,d\}$, we define the vector $v(k) \in V_{\mz}^n$ by 
	\begin{equation*}
		v(k)=\sum_{i=1}^k p_i(v) e_i ,
	\end{equation*}
	i.e., the first $k$ coordinates of $v(k)$ equal the corresponding coordinates of $v$ and all other coordinates are $0$.
	By the triangle inequality and linearity of expectation we clearly have
	\begin{equation*}
		\E_\beta \left[ D_{V_{\mz}^n}(\mz,v) \right]  \leq \E_{\beta } \left[ \sum_{i=0}^{d-1} D_{V_{\mz}^n} (v(i),v(i+1)) \right]
		=
		\sum_{i=0}^{d-1} \E_{\beta } \left[  D_{V_{\mz}^n} (v(i),v(i+1)) \right]  \text .
	\end{equation*}
	 So in order to show \eqref{eq:endisfar reduction}, it suffices to show that 
	 \begin{equation}\label{eq:toshow scaling}
	 \E_\beta  \left[ D_{V_{\mz}^n} (v(i),v(i+1)) \right]  \leq 3 \E_\beta\left[D_{V_{\mz}^{n}} (\mz,(n-1) e_1)\right] 
	 \end{equation}
	 for all $i \in \{0,\ldots, d-1\}$. For each  such index $i$, the cube 
	\begin{align*}
		\mathcal{B}_i = \prod_{j=1}^{i} \left\{p_j(v)-p_{i+1}(v), \ldots, p_j(v)\right\} \times \{0,\ldots, p_{i+1}(v)\}^{d-i}
	\end{align*}
	is contained in the cube $V_{\mz}^n$. These cubes $\left(\mathcal{B}_i\right)_{i \in \{0,\ldots, d-1\}}$ are chosen in such a way that the cube $\mathcal{B}_i \subset V_\mz^n$, both points $v(i)$ and $v(i+1)$ are corners of the cube $\mathcal{B}_i$ and the line-segment connecting $v(i)$ to $v(i+1)$ is an edge of the cube. This property will then allow us to use the symmetry of the model, together with Proposition \ref{propo:scaling}. 
		\begin{figure}
		\begin{tikzpicture}[scale=0.75]
		\fill [superlightgray] (0,0) rectangle (7,7);
		\fill[pattern=north west lines, pattern color=gray] (3,0) rectangle (6,3);
		\foreach \x in {-1,0,1,2,3,4,5,6,7} {
			\foreach \y in {-1,0,1,2,3,4,5,6,7} {
				\vertex[fill,color=wellgray,minimum size=4 pt ]  at (\y,\x) {};
			};
		};
		\vertex[draw=none] (A) at (-1.2,0) {};
		\vertex[draw=none] (B) at (7.5,0) {};
		\path[very thick , ->] (A) edge (B);
		\vertex[draw=none] (C) at (0,-1.2) {};
		\vertex[draw=none] (D) at (0, 7.5) {};
		\path[very thick , ->] (C) edge (D);
		\vertex[fill, color=black, minimum size = 5.5pt, label=$v(0) \ \ \ \ \ \ $] (v0) at (0, 0) {};
		\vertex[fill, color=black, minimum size = 5.5pt, label=$\ \ \ \ \ \ v(1)   $] (v1) at (6, 0) {};
		\vertex[fill, color=black, minimum size = 5.5pt, label=$v(2)$] (v2) at (6, 3) {};
		\end{tikzpicture}
		\centering
		\parbox{11cm}{\caption{Let $v=(6,3)\in V_\mz^8$. The points in the gray area are the set $V_\mz^8$. The points in the hatched area are $\mathcal{B}_1$.} \label{Fig:B}}
	\end{figure}
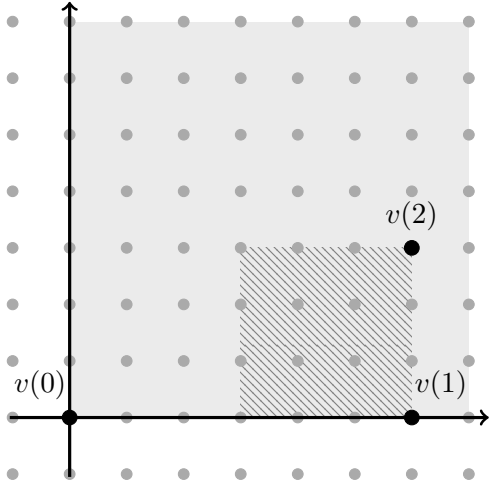
	See figure \ref{Fig:B} for an example. Allowing the geodesic to use less edges clearly increases the distance between two points, which implies $D_{V_{\mz}^n} (v(i),v(i+1)) \leq D_{\mathcal{B}_i} (v(i),v(i+1))$ as $\mathcal{B}_i \subset V_{\mz}^n$. As the model is invariant under changing the coordinates and under the action $e_i \mapsto -e_i$ we already get for all $i\in \{0,\ldots,d-1\}$
	\begin{align*}
		\E_\beta \left[ D_{\mathcal{B}_i} (v(i),v(i+1)) \right] = \E_\beta \left[ D_{V_{\mz}^{p_{i+1}(v)+1}} (\mz, p_{i+1}(v) e_1) \right] \leq 3 \E_\beta \left[ D_{V_{\mz}^{n}} (\mz, (n-1) e_1) \right]\text,
	\end{align*}
	where we used Proposition \ref{propo:scaling} for the last inequality. This shows \eqref{eq:toshow scaling} and thus finishes the proof of \eqref{eq:endisfar}.
	Now let us go to the proof of \eqref{eq:endisfar2}. Define $y\in \Z^d$ by $p_1(y)=1 , p_i(y) = -1$ for $i\geq 2$ and define the cube $\mathcal{B}$ by $\mathcal{B}=\{n-1,\ldots, 2n-2\} \times \{0,\ldots,n-1\}^{d-1}$. By the triangle inequality we have
	\begin{equation}\label{eq:sum gleich}
		D_{V_{\mz}^{2n-1}} (\mz, (2n-2)e_1) \leq D_{V_{\mz}^{n}} (\mz, (n-1)\mo) + 
		D_{\mathcal{B}} ((n-1)\mo, (2n-2)e_1) \text .
	\end{equation}
	Observe that $(2n-2)e_1 = (n-1)\mo + (n-1)y$. The pairs of vertices $\mz$ and $(n-1)\mo$ lie on opposite corners of the cube $V_\mz^n$. The vertices $(n-1)\mo$ and $(2n-2)e_1$ also lie on opposite corners of the cube $\cB$. The two cubes $V_\mz^n$ and $\cB$ differ by a translation only; in particular, they have the same side length. As the long-range percolation model is invariant under translation and reflection of any coordinate the two terms in the sum \eqref{eq:sum gleich} have the same distribution which implies that
	\begin{align*}
		\E_\beta \left[ D_{V_{\mz}^{2n-1}} (\mz, (2n-2)e_1) \right]
		\leq 2 \E_\beta \left[ D_{V_{\mz}^{n}} (\mz, (n-1)\mo) \right] \text .
	\end{align*}
	Using Proposition \ref{propo:scaling}, we finally get
	\begin{align*}
		\E_\beta \left[D_{V_{\mz}^n}(\mz,(n-1)e_1)\right] \leq 
		3 \E_\beta \left[ D_{V_{\mz}^{2n-1}} (\mz, (2n-2)e_1) \right]\leq
		6 \E_\beta \left[ D_{V_{\mz}^{n}} (\mz, (n-1)\mo) \right]
	\end{align*}
	which shows \eqref{eq:endisfar2}.	
\end{proof}

\subsection{The second moment bound}\label{sec:2nd moment}

The next lemma relates the second moment of the distances to their first moment. We use a technique that has already been used in \cite{ding2013distances} before in a slightly different form for dimension $d=1$ only. As we need the result in a uniform dependence on $\beta$ in our companion paper \cite{baeumler2022behaviour}, we directly prove the uniform statement here. The uniformity does not cause any complications for $d\geq 2$, but it causes minor technical difficulties for $d=1$. So we give the proof for $d=1$ in the \hyperlink{target:d=1}{appendix}. The situation for $d\geq 2$ is easier, as there are no cut points, in the sense that for every $u,v \in V_{\mz}^n$ there exist two disjoint paths between $u$ and $v$. For $d=1$, and in particular for $\beta < 1$, such a statement will typically not be true.

\begin{lemma}\label{lem:secondmomentbound}
	For all $\beta \geq 0$, there exists a constant $C_\beta<\infty$ such that for all $n \in \N$, all $\eps \in \left[0,1\right]$ and all $u,v \in V_{\mz}^n$
	\begin{align}\label{eq:second moment bound}
		\E_{\beta+\eps} \left[ D_{V_{\mz}^n} (u,v)^2 \right] \leq C_\beta \Lambda(n,\beta+\eps)^2.
	\end{align}
\end{lemma}

\begin{proof}[Proof of Lemma \ref{lem:secondmomentbound} for $d\geq 2$]
	Fix $\beta \geq 0$. We will prove that for all $\eps \in \left[0,1\right]$, all $m,n \in \N$, and all $u,v \in V_{\mz}^{mn}$
	\begin{align}\label{eq:second moment bound to show}
		\E_{\beta+\eps} \left[D_{V_{\mz}^{mn}} (u,v)^2 \right] \leq
		300 m^4 \Lambda(n,\beta+\eps)^2 +
		300 \max_{w,z  \in V_{\mz}^n } \E_{\beta+\eps}\left[D_{V_{\mz}^n}(w,z)^2\right]\text .
	\end{align}
	Iterating over this inequality one gets for some large enough $N$ that
	\begin{align}\label{eq:iterated bound second moment}
	& \notag \max_{u,v \in V_{\mz}^{m^k N}} \E_{\beta+\eps} \left[ D_{V_{\mz}^{m^k N}} (u,v)^2 \right]\\
	& \notag \leq 300 m^4 \sum_{i=0}^{k} 300^{i} \Lambda(m^{k-i}N,\beta+\eps)^2 + 300^k \max_{u,v \in V_{\mz}^{N}} \E_{\beta+\eps} \left[ D_{V_{\mz}^{N}} (u,v)^2 \right] \\
	& \leq 300 m^4 \sum_{i=0}^{k} 300^{i} \Lambda(m^{k-i}N,\beta+\eps)^2 + 300^k N^2.
	\end{align}
	for all $k \in \N$.	By \Cref{remark:Lambda grows} there exists $\dxp^\star = \dxp^\star(\beta) > 0$ such that for all $\eps \in \left[0,1\right]$, and all $m,n \in \N$ large enough one has
	\begin{align*}
		& \Lambda(mn,\beta+\eps) = \max_{u,v \in V_{\mz}^{mn}} \E_{\beta+\eps} \left[D_{V_{\mz}^{mn}} (u,v) \right] + 1 \\
		&
		\geq 
		m^{\dxp^\star}  \left(\max_{u,v \in V_{\mz}^{n}} \E_{\beta+\eps} \left[D_{V_{\mz}^{n}} (u,v) \right] + 1\right)
		= m^{\dxp^\star}  \Lambda(n,\beta+\eps) \text .
	\end{align*}
	Take $m$ large enough so that also $300m^{-\dxp^\star} < \frac{1}{2}$ is satisfied. Inserting this into $\eqref{eq:iterated bound second moment}$ gives
	\begin{align*}
		\max_{u,v \in V_{\mz}^{m^kN}} & \E_{\beta+\eps} \left[ D_{V_{\mz}^{m^kN}} (u,v)^2 \right] \leq  300 m^4 \sum_{i=0}^{k} 300^{i} \Lambda(m^{k-i}N,\beta+\eps)^2 + 300^k N^2\\
		& \leq 300 m^4 \sum_{i=0}^{k} 300^{i} m^{-\dxp^\star i} \Lambda(m^{k}N,\beta+\eps)^2 + N^2 \Lambda(m^k N,\beta+\eps)\\
		& \leq \left(600 m^4 + N^2\right) \Lambda(m^kN,\beta+\eps)^2
	\end{align*}
	for large enough $N$. This shows \eqref{eq:second moment bound} along the subsequence $N, mN, m^2N, \ldots$ . For general $n \in \N$, the desired result follows from Proposition \ref{propo:scaling}. So we are left with showing \eqref{eq:second moment bound to show}. For this, we use an elementary observation, that was already used in \cite{ding2013distances}. Assume that $X_1, \ldots, X_{\tilde{m}}$ are independent non-negative random variables and let $\tau = \arg \max_{i \in \{1,\ldots, \tilde{m}\}} \left(X_i\right)$. Then
	\begin{align}\label{eq:secondlargest bounds}
			\E\left[\left(\max_{i \neq \tau} X_i\right)^2\right] & \leq \E \left[ \sum_{i=1}^{\tilde{m}}  X_i \left(\sum_{j \neq i} X_j\right) \right] = \sum_{i=1}^{\tilde{m}} \sum_{j \neq i} \E \left[ X_i \right] \E\left[ X_j \right] \leq \tilde{m}^2 \max_{i} \E\left[X_i\right]^2  .
	\end{align}
	We still need to show inequality \eqref{eq:second moment bound to show}, i.e.,  that
	\begin{align*}
	\E_{\beta+\eps} \left[D_{V_{\mz}^{mn}} (u,v)^2 \right] \leq
	300 m^4 \Lambda(n,\beta+\eps)^2 +
	300 \max_{w,z  \in V_{\mz}^n } \E_{\beta+\eps}\left[D_{V_{\mz}^n}(w,z)^2\right] \text .
	\end{align*}
	Let $u,v \in V_{\mz}^{mn}$, say with $u \in V_x^n , v \in V_y^n$, where $x,y \in V_{\mz}^m$. Inequality \eqref{eq:second moment bound to show} clearly holds in the case where $x=y$. For the case $x \neq y$, let $x_0=x, x_1, \ldots, x_l = y$ and  $x^\prime_0=x, x^\prime_1, \ldots, x^\prime_{l^{\prime}} = y$ be two completely disjoint and deterministic paths between $x$ and $y$ inside $V_{\mz}^m$ that are of length at most $m+1$ and use only nearest-neighbor edges, i.e.,  $\|x_i-x_{i-1}\|_\infty =1$ and $\|x_i^\prime - x_{i-1}^\prime\|_\infty =1$ for all suitable indices $i$. By completely disjoint we mean that $\left\{x_1,\ldots, x_{l-1}\right\} \cap \left\{x_1^\prime ,\ldots, x_{l^\prime-1}^\prime\right\} = \emptyset$; the starting point $x=x_0=x_0^\prime$ and the end point $y=x_l=x_{l^\prime}^\prime$ already need to agree by the construction. Now we iteratively define sequences $(L_i,R_i)_{i=0}^l$ and $(L_i^\prime,R_i^\prime)_{i=0}^{l^\prime}$ as follows:
	\begin{enumerate}\addtocounter{enumi}{-1}
		\item Set $L_0 = u, R_l=v$.
		\item For $i=1,\ldots, l$, choose $R_{i-1} \in V_{x_{i-1}}^n$ and $L_i \in V_{x_i}^n$ such that ${\|R_{i-1}-L_i\|_\infty = 1}$.
	\end{enumerate}
	Analogously, we define $(L_i^\prime,R_i^\prime)_{i=0}^{l^\prime}$ by
	\begin{enumerate}\addtocounter{enumi}{-1}
		\item Set $L^\prime_0 = u, R^\prime_{l^\prime}=v$.
		\item For $i=1,\ldots, l^\prime$, choose $R_{i-1}^\prime \in V_{x_{i-1}^\prime}^n$ and $L_i^\prime \in V_{x_i^\prime}^n$ such that ${\|R^\prime_{i-1}-L^\prime_i\|_\infty = 1}$.
	\end{enumerate}
	
	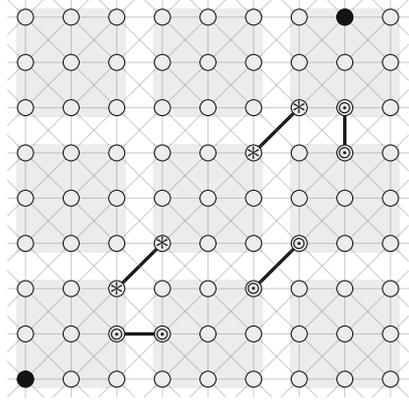
\begin{figure}
		\[\begin{tikzpicture}[scale = 0.6]
			
			\draw[black, very thick] (2,1) -- (3,1);
			\draw[black, very thick] (5,2) -- (6,3);
			\draw[black, very thick] (7,5) -- (7,6);
			\draw[black, very thick] (2,2) -- (3,3);
			\draw[black, very thick] (5,5) -- (6,6);
			
			\draw[gray, thin, opacity=0.4] (-0.4,-0.4) grid (8.5,8.5);
			
			\foreach \n in {0,...,8}
			\draw[gray, thin, opacity=0.4] (0-0.4,\n+0.4) -- (\n+0.4,0-0.4);
			
			\foreach \n in {1,...,8}
			\draw[gray, thin, opacity=0.4] (\n-0.4,8+0.4) -- (8+0.4,\n-0.4);
			
			\foreach \n in {0,...,8}
			\draw[gray, thin, opacity=0.4] (0-0.4,8-\n-0.4) -- (\n+0.4,8+0.4);
			
			\foreach \n in {1,...,8}
			\draw[gray, thin, opacity=0.4] (\n-0.4,0-0.4) -- (8+0.4,8-\n+0.4);
			
			\foreach \x in {0,...,8} 
			\foreach \y in {0,...,8} 
			\vertex[fill=white]  at (\x,\y) {};
			
			\vertex[fill=black]  at (0,0) {};
			\vertex[fill=black]  at (7,8) {};
			
			\vertex[draw=none]  at (2,2) {$*$};
			\vertex[draw=none]  at (3,3) {$*$};
			\vertex[draw=none]  at (5,5) {$*$};
			\vertex[draw=none]  at (6,6) {$*$};

			\vertex[draw=none]  at (2,1-0.02) {$\circ$}; \vertex[draw=none]  at (2,1-0.02) {$\cdot$};
			\vertex[draw=none]  at (3,1-0.02) {$\circ$}; \vertex[draw=none]  at (3,1-0.02) {$\cdot$};
			\vertex[draw=none]  at (5,2-0.02) {$\circ$}; \vertex[draw=none]  at (5,2-0.02) {$\cdot$};
			\vertex[draw=none]  at (6,3-0.02) {$\circ$}; \vertex[draw=none]  at (6,3-0.02) {$\cdot$};
			\vertex[draw=none]  at (7,5-0.02) {$\circ$}; \vertex[draw=none]  at (7,5-0.02) {$\cdot$};
			\vertex[draw=none]  at (7,6-0.02) {$\circ$}; \vertex[draw=none]  at (7,6-0.02) {$\cdot$};

			\foreach \a in {0,...,2} 
			\foreach \b in {0,...,2}
			\filldraw[gray,fill opacity=0.15, draw=none] (3*\a-0.2,3*\b-0.2) rectangle (3*\a+2.2,3*\b+2.2);
			
		\end{tikzpicture}\]
		\centering
		\parbox{12cm}{\caption{The graph $V_\mz^9$ for $d=2$. The division into $3\times 3$-boxes is marked in gray. In this picture, $u,v$ are the black vertices, $R_0,L_1,R_1,L_2$ are the starred vertices and the nearest-neighbor edges $\{R_0, L_1\}$ and $\{R_1,L_2\}$ are black and thick. The vertices $R_0^\prime,L_1^\prime,R_1^\prime,L_2^\prime,R_2^\prime,L_3^\prime$are the vertices with the extra circle and dot in the inside; The nearest-neighbor edges \{$R_0^\prime,L_1^\prime\}, \{R_1^\prime,L_2^\prime\}, \text{ and } \{R_2^\prime,L_3^\prime\}$ are black and thick.} \label{fig:2ndmom}}
		
	\end{figure}

	The choice of these algorithms in step $1$ is typically not unique. If there are several possibilities, we always choose the vertices with some deterministic rule that does not depend on the environment. An example of such a construction is given in Figure \ref{fig:2ndmom}. The idea behind this construction is, that there exists a path that goes from $L_0$ to $R_0$ to $L_1$ to $R_1$ ... to $L_l$ to $R_l$. Furthermore, there also exists a path that goes from $L_0^\prime$ to $R_0^\prime$ to $L_1^\prime$ ... to $L_{l^\prime}^\prime$ to $R_{l^\prime}^\prime$. We then compare the length of these two paths.
	
	 By construction we have $L_i, R_i \in V_{x_i}^n$ and $L_i^\prime, R_i^\prime \in V_{x_i^\prime}^n$ for all $i \in \{0,\ldots, l\}$, respectively $i \in \{0,\ldots, l^\prime\}$. Define 
	\begin{align*}
		X_i = D_{V_{x_i}^n}(L_i, R_i) & \text{ for } i\in \{1,\ldots, l-1\} \text{ and }\\
		X_i^\prime = D_{V_{x_i^\prime}^n}(L_i^\prime , R_i^\prime ) & \text{ for } i\in \{1,\ldots, l^\prime-1\} \text .
	\end{align*}
	These are at most $l-1+l^\prime -1 \leq2m$ random variables and they are independent, as the boxes $V_{x_i^\prime}^n$ and $V_{x_i}^n$ are disjoint. We order the random variables $\left\{X_i : i \in \{1,\ldots,l-1\}\right\} \cup \left\{X_i^\prime : i \in \{1,\ldots,l^\prime-1\}\right\}$ in a descending way and call them $Y_1,Y_2,\ldots, Y_{l+l^\prime -2 }$. The idea in finding a short path between $u$ and $v$ is now to avoid the box where the maximum of the $Y_i$-s is attained. Assume that the maximum of them is one of the $X_i$-s, i.e., $X_i=Y_1$ for some $i\in \{1,\ldots,l-1\}$. Then we consider the path that goes from $L_0^\prime =u$ to $R_0^\prime$ and from there to $L_1^\prime$, and from there we go successively to $R_{l^\prime}^\prime =v$. Otherwise, we have $X_i^\prime = Y_1$ for some $i \in \{1,\ldots,l^\prime-1\}$. In this situation, we consider the path that goes from $L_0 =u$ to $R_0$, from there to $L_1$, and successively to $R_{l} =v$. In both cases we have constructed a path between $u$ and $v$. The length of this path is an upper bound on the chemical distance between $u$ and $v$ and thus we get
	\begin{align}\label{eq:elementary counting}
		D_{V_{\mz}^{mn}} (u,v)  \leq   D_{V_x^n}(L_0, R_0) +  D_{V_x^n}(L_0^\prime , R_0^\prime)  + D_{V_y^n}(L_l, R_l) +  D_{V_y^n}(L_{l^\prime}^\prime , R_{l^\prime }^\prime) + m Y_2 + (m+1) \text ,
	\end{align}
	where the summand $(m+1)$ arises as one still needs to go from $R_i$ to $L_{i+1}$ for all $i\in \{0,\ldots, l-1\}$, or from $R_i^\prime$ to $L_{i+1}^\prime$ for all $i\in \{0,\ldots, l^\prime-1\}$. But by assumption one has $l, l^\prime \leq m+1$, so one needs at most $m+1$ additional steps.
	From \eqref{eq:secondlargest bounds} we know that 
	\begin{align}\label{eq:2nd moment bound Y_2}
		\E_{\beta + \eps} \left[Y_2^2\right] \leq 4m^2 \max_{w,z \in V_\mz^n}  \E_{\beta + \eps} \left[D_{V_\mz^n}(w,z)  \right]^2.
	\end{align}
	For the distance between $L_0 $ and $R_0$ one clearly has
	\begin{align*}
		\E_{\beta + \eps } \left[ D_{V_x^n}(L_0, R_0)^2  \right] 
		\leq 
		\max_{w,z \in V_{\mz}^{n}} \E_{\beta+\eps} \left[D_{V_{\mz}^n}(w,z)^2\right]
	\end{align*}
	and the same statements hold for $D_{V_x^n}(L_0^\prime , R_0^\prime) , D_{V_y^n}(L_l, R_l), $ and $ D_{V_y^n}(L_{l^\prime}^\prime , R_{l^\prime }^\prime)$. Using the elementary inequality $\left(\sum_{i=1}^{6} a_i\right)^2 \leq 36 \sum_{i=1}^{6} a_i^2$ that holds for any six numbers $a_1,\ldots, a_6 \in \R$ for the term in \eqref{eq:elementary counting}, we get that
	\begin{align*}
		& \E_{\beta+\eps} \left[D_{V_{\mz}^{mn}} (u,v)^2 \right] \\
		& \leq   36 \E_{\beta+\eps} \left[ D_{V_x^n}(L_0, R_0)^2 +  D_{V_x^n}(L_0^\prime , R_0^\prime)^2  + D_{V_y^n}(L_l, R_l)^2 +  D_{V_y^n}(L_{l^\prime}^\prime , R_{l^\prime }^\prime)^2 + m^2 Y_2^2 + (m+1)^2 \right]\\
		& \leq 4\cdot 36 \max_{w,z \in V_{\mz}^{n}} \E_{\beta+\eps} \left[D_{V_{\mz}^n}(w,z)^2\right] + 36 m^2 \E_{\beta+\eps} \left[Y_2^2\right] + 36(m+1)^2\\
		& \overset{\eqref{eq:2nd moment bound Y_2}}{\leq} 300 \max_{w,z \in V_{\mz}^{n}} \E_{\beta+\eps} \left[D_{V_{\mz}^n}(w,z)^2\right] + 300 m^4 \Lambda(n,\beta+\eps)^2
	\end{align*}
	which shows \eqref{eq:second moment bound to show} and thus finishes the proof.
\end{proof}

\begin{corollary}\label{coro:allmoments}
	Iterating this technique one can show that for all $k\in \N$ of the form $k=2^l$ and for all $\beta > 0$ there exists a constant $C_\beta<\infty$ such that for all $n \in \N$, and all $u,v \in V_{\mz}^n$
	\begin{align}\label{eq:all moment bound}
	\E_{\beta} \left[ D_{V_{\mz}^n} (u,v)^k \right] \leq C_\beta \Lambda(n,\beta)^k.
	\end{align}
	Then, one can extend this bound to all $k\in \R_{\geq 0}$ with Hölder's inequality.
\end{corollary}

\begin{proof}[Proof of \Cref{coro:allmoments} for $d\geq 2$]
	For $r> 0$, define the quantity 
	\begin{align*}
		\Lambda^r(\beta,n) \coloneqq \max_{x,y \in V_\mz^n} \E_\beta \left[ D_{V_\mz^n} (x,y)^r \right] 
	\end{align*}
	and assume that $\Lambda^r(\beta,n) \leq C \Lambda(\beta,n)^r$ for some constant $C$ and all $n\in \N$. Using the same notation as in \eqref{eq:elementary counting} above we get that for any $u,v \in V_\mz^{mn}$, say with $u\in V_x^n$ and $y\in V_y^n$,
	\begin{align*}
	D_{V_{\mz}^{mn}} (u,v)  \leq   D_{V_x^n}(L_0, R_0) +  D_{V_x^n}(L_0^\prime , R_0^\prime)  + D_{V_y^n}(L_l, R_l) +  D_{V_y^n}(L_{l^\prime}^\prime , R_{l^\prime }^\prime) + m Y_2 + (m+1) \text ,
	\end{align*}
	and thus we also get that
	\begin{align*}
	D_{V_{\mz}^{mn}} (u,v)^{2r}  \leq 6^{2r}  \Big( & D_{V_x^n}(L_0, R_0)^{2r} +  D_{V_x^n}(L_0^\prime , R_0^\prime)^{2r} \\
	& + D_{V_y^n}(L_l, R_l)^{2r} +  D_{V_y^n}(L_{l^\prime}^\prime , R_{l^\prime }^\prime)^{2r} +( m Y_2)^{2r} + (m+1)^{2r} \Big) \text .
	\end{align*}
	We have that
	\begin{align*}
		\E_{\beta } \left[ ( m Y_2)^{2r}  \right] = m^{2r} \E_\beta \left[ \left(Y_2^r\right)^2 \right] \leq m^{2r}  \max_{w,z  \in V_{\mz}^n } \E_\beta \left[ D_{V_\mz^n}(w,z)^r \right]^2
		\leq m^{2r} C^2 \Lambda(\beta,n)^r
	\end{align*}
	and from here the same proof as in \Cref{lem:secondmomentbound} shows that $\Lambda^{2r}(\beta,n) \leq C(r) \Lambda(\beta,n)^{2r}$ for some constant $C(r)< \infty$. Inductively, we thus get that for all $r=2^k$, with $k\in \N$ one has $\Lambda^r(\beta,n) \leq C(r) \Lambda(\beta,n)^r$. Whenever $r>0$ is not of the form $r=2^k$ for some $k\in \N$, let $k$ be large enough so that $r < 2^k$. Then we get that
	\begin{align*}
		\Lambda^r(\beta,n) = \max_{x,y \in V_\mz^n} \E_\beta \left[ D_{V_\mz^n} (x,y)^r \right] 
		\leq
		\max_{x,y \in V_\mz^n} \E_\beta \left[ D_{V_\mz^n} (x,y)^{2^k} \right]^{\frac{r}{2^k}} \leq C \Lambda(\beta,n)^r
	\end{align*}
	for some constant $C$.
\end{proof}

\subsection{Graph distances between points and boxes}\label{sec:distances between points and boxes}

So far, we only considered distances between two different points in a box. In this section, we investigate the distance between certain points and boxes. For $n\in \N$ and $0< \iota < \frac{1}{2}$ we define the boxes $L_\iota^n \coloneqq \left[0,\iota n\right]^d$ and $R_\iota^n \coloneqq \left[n-1-\iota n, n-1\right]^d$. These are boxes that lie in opposite corners of the cube $V_{\mz}^n$, where $L_\iota^n$ lies in the corner containing $\mz$ and $R_\iota^n$ lies in the corner containing $\mo$. The next lemma deals with the graph distance of these two boxes. A similar statement of Lemma \ref{lem:set_to_set} for the continuous model and $d=1$, was already proven in \cite{ding2013distances}. We follow the same strategy for the proof of this lemma. Again, we prove it uniformly for $\beta$ in some compact intervals, as we will need this uniform statement in \cite{baeumler2022behaviour}. The uniformity does not make any complications in this proof here.

\begin{lemma}\label{lem:set_to_set}
	For all $\beta \geq 0$, there exists an $\iota >0$ such that uniformly over all $\eps \in \left[0,1\right]$ and $n\in \N$
	\begin{align}\label{eq:set_to_set expectation}
		\E_{\beta+\eps} \left[ D_{V_{\mz}^n}\left( L_\iota^n , R_\iota^n \right) \right] \geq \frac{1}{2} \E_{\beta + \eps} \left[ D_{V_{\mz}^n}\left( \mz , (n-1) \mo \right) \right] \text , 
	\end{align}
	and there exists $c^\star >0$ such that uniformly over all $\eps \in \left[0,1\right]$ and $n\in \N$
	\begin{align}\label{eq:set_to_set probability}
		\p_{\beta + \eps} \left(  D_{V_{\mz}^n}\left( L_\iota^n , R_\iota^n \right) \geq \frac{1}{4} \E_{\beta+\eps} \left[  D_{V_{\mz}^n}\left( \mz , (n-1)\mo \right) \right] \right) \geq c^\star  .
	\end{align}
\end{lemma}
\begin{proof}
	The statement clearly holds for small $n$, so we focus on $n \in \N$ large enough from here on. 
	Let $x \in L_\iota^n$ and $y \in R_\iota^n$ be the minimizers of $D_{V_{\mz}^n}\left( L_\iota^n , R_\iota^n \right)$, i.e., $D_{V_{\mz}^n}\left( L_\iota^n , R_\iota^n \right) = D_{V_{\mz}^n}\left( x,y \right)$. If the minimizers are not unique, pick two minimizers with a deterministic rule not depending on the environment. The choice of $x,y$, and the distance $D_{V_{\mz}^n}\left( L_\iota^n , R_\iota^n \right)$ depend only on edges with at least one endpoint in $V_{\mz}^n \setminus \left( L_\iota^n \cup R_\iota^n\right)$. The distances $D_{L_\iota^n} (\mz,x)$, respectively $D_{R_\iota^n} (y,(n-1)\mo)$, depend only on edges with both endpoints in $L_\iota^n$, respectively $R_\iota^n$. Thus we get that
	\begin{align*}
		\E_{\beta+\eps} \left[ D_{V_{\mz}^n} (\mz,(n-1)\mo) \right] 
		& \leq  
		\E_{\beta+\eps}\left[ D_{R_\iota^n} (\mz,x) \right] 
		+ \E_{\beta+\eps}\left[ D_{V_{\mz}^n} (L_\iota^n, R_\iota^n) \right]
		+ \E_{\beta+\eps}\left[ D_{R_\iota^n} (y,(n-1)\mo) \right]\\
		& \leq 2 \Lambda(\lfloor \iota n \rfloor, \beta+\eps ) + \E_{\beta+\eps}\left[ D_{V_{\mz}^n} (L_\iota^n, R_\iota^n) \right] .
	\end{align*}
	For $\iota$ small enough and $n$ large enough, we get uniformly over $\eps \in \left[0,1\right]$ that 
	\begin{align*}
		\Lambda( n, \beta+\eps ) \geq \left(\frac{1}{\iota}\right)^{\dxp^\prime}   \Lambda(\lfloor \iota n \rfloor, \beta+\eps ) 
	\end{align*}
	for some $\dxp^\prime > 0$ by Remark \ref{remark:Lambda grows}. So by Lemma \ref{lem:endisfar}, respectively Remark \ref{remark:Lambda}, we can choose $\iota$ small enough so that	uniformly over $n \in \N$ large enough and $\eps \in \left[0,1\right]$
	\begin{align*}
		2 \Lambda(\lfloor \iota n \rfloor, \beta+\eps ) \leq \frac{1}{2} 	\E_{\beta+\eps} \left[ D_{V_{\mz}^n} (\mz,(n-1)\mo) \right] \text ,
	\end{align*}
	and this implies that
	\begin{align*}
		& \E_{\beta+\eps}\left[ D_{V_{\mz}^n} (L_\iota^n, R_\iota^n) \right] \geq \E_{\beta+\eps} \left[ D_{V_{\mz}^n} (\mz,(n-1)\mo) \right] - 2 \Lambda(\lfloor \iota n \rfloor , \beta + \eps ) \\
		& \geq \frac{1}{2} \E_{\beta+\eps} \left[ D_{V_{\mz}^n} (\mz,(n-1)\mo) \right] 
	\end{align*}
	which proves \eqref{eq:set_to_set expectation}.
	For such an $\iota$, define $\mathcal{A} = \left\{ D_{V_{\mz}^n}\left( L_\iota^n , R_\iota^n \right) \geq \frac{1}{4} \E_{\beta+\eps} \left[  D_{V_{\mz}^n}\left( \mz, (n-1) \mo \right) \right]\right\}$. By the Cauchy-Schwarz inequality we have
	\begin{align*}
		& \E_{\beta+\eps} \left[ D_{V_{\mz}^n} (\mz,(n-1)\mo) \right]
		\leq 2
		 \E_{\beta+\eps} \left[ D_{V_{\mz}^n} (L_\iota^n, R_\iota^n )   \right]\\
		 &
		 = 2 \E_{\beta+\eps} \left[ D_{V_{\mz}^n} (L_\iota^n, R_\iota^n )  \mathbbm{1}_{\mathcal{A}^C} \right]
		+ 2 \E_{\beta+\eps} \left[ D_{V_{\mz}^n} (L_\iota^n, R_\iota^n )  \mathbbm{1}_{\mathcal{A}} \right]\\
		& \leq \frac{1}{2} \E_{\beta+\eps} \left[ D_{V_{\mz}^n} (\mz,(n-1)\mo) \right] + 2  \E_{\beta+\eps} \left[ D_{V_{\mz}^n} (\mz,(n-1)\mo)^2 \right]^{1/2} \sqrt{\p_{\beta + \eps}\left( \mathcal{A} \right)}\\
		& \leq \frac{1}{2} \E_{\beta+\eps} \left[ D_{V_{\mz}^n} (\mz,(n-1)\mo) \right] + C^\prime \E_{\beta+\eps} \left[ D_{V_{\mz}^n} (\mz,(n-1)\mo) \right] \sqrt{\p_{\beta + \eps}\left( \mathcal{A} \right)} \text ,
	\end{align*}
	where the last inequality holds for some $C^\prime < \infty$, by Lemma \ref{lem:endisfar} and Lemma \ref{lem:secondmomentbound}. Solving the previous line of inequalities for $\p_{\beta + \eps}\left( \mathcal{A} \right)$ shows \eqref{eq:set_to_set probability}.
\end{proof}

\begin{lemma}\label{lem:point to pointonsphere}
	For all $\beta \geq 0$ and all dimensions $d$, there exists a constant $c_1>0$ such that uniformly over all $n\in \N$ and all $x \in \mathcal{S}_n$
	\begin{align}\label{eq:point to pointonsphere}
		\E_\beta \left[ D_{B_n(\mz)} (\mz,x) \right] \geq c_1 \E_\beta \left[ D_{V_{\mz}^n}(\mz,(n-1)\mo ) \right]
	\end{align}
	and the constant $c_1$ can be chosen in such a way so that it only depends on the dimension $d$ and the value $\iota > 0$ in \eqref{eq:set_to_set expectation}.
\end{lemma}
\begin{proof}
	Let $v \in \mathcal{S}_n$ be one of the minimizers of $y \mapsto \E_\beta \left[D_{B_n(\mz)}(\mz,y)\right]$ among all vertices $y \in \mathcal{S}_n$. By reflection symmetry, we can assume that all coordinates of $v$ are non-negative. With the notation $e_0=e_d$ we define the vectors $v_0,\ldots,v_{d-1}$ by
	\begin{align*}
		\langle e_j , v_i \rangle = \langle e_{i+j \text{ mod } d} , v \rangle
	\end{align*}
	which are just versions of the vector $v$, where we cyclically permuted the coordinates. By invariance under changes of coordinates, we have
	\begin{align*}
		\E_\beta \left[ D_{B_n(\mz)}(\mz,v) \right] = \E_\beta \left[ D_{B_n(\mz)}(\mz,v_i) \right]
	\end{align*}
	for all $i \in \{0,\ldots, d-1 \}$. Define the vertices $u_0,\ldots, u_{d}$ by $u_j = \sum_{i=1}^{j} v_i$. By our construction we have $u_0 = \mz$ and $u_{d} = \sum_{i=1}^{d} v_i = N \mo$ for some integer $N \geq n$. The balls $B_n(u_i)$ are all contained in the cube $\Upsilon=\{-n,\ldots, N+n\}^d$ for all $i\in \{0,\ldots, d\}$. Thus we have
	\begin{align*}
		\E_\beta\left[ D_\Upsilon (\mz, N \mo ) \right] 
		\leq \sum_{i=1}^{d-1} \E_\beta\left[ D_\Upsilon ( u_{i-1} , u_{i-1} + v_i \mo ) \right] 
		\leq d \E_\beta \left[ D_{B_n(\mz)}(\mz,v) \right] \text ,
	\end{align*}
	and by translation invariance we also have for the cube $\Upsilon_1 = \{0,\ldots, 2n + N\}^d$
	\begin{align*}
		\E_\beta\left[ D_{\Upsilon_1 } (n \mo , (n+N) \mo ) \right] 
		\leq d \ \E_\beta \left[ D_{B_n(\mz)}(\mz,v) \right] \text .
	\end{align*}
	Using the triangle inequality, we see that for all $k\in \N$  the expected distance between $n\mo$ and $(n+kN)\mo$ inside the cube $\Upsilon_k = \{0,\ldots, 2n + kN\}^d $ is upper bounded by 
	\begin{align*}
		\E_\beta\left[ D_{\Upsilon_k } (n \mo , (n+kN) \mo ) \right] 
		\leq k \cdot d \ \E_\beta \left[ D_{B_n(\mz)}(\mz,v) \right] \text . 
	\end{align*}
	But Proposition \ref{propo:scaling} also gives that for $s=\frac{n}{kN + 2n}$ and $w_1 = \lfloor s n\mo \rfloor, w_2 = \lfloor s (n+kN)\mo \rfloor$
	\begin{align*}
		\E_\beta \left[ D_{V_{\mz}^n} (w_1,w_2) \right] \leq 3 k \cdot d \ \E_\beta \left[ D_{B_n(\mz)}(\mz,v) \right] \text .
	\end{align*}
	As $N\geq n$, for each fixed $\iota >0$ we can choose $k$ large enough so that $w_1 \in L_\iota^n$ and $w_2 \in R_\iota^n$ and thus $\E_\beta \left[ D_{V_{\mz}^n} (w_1,w_2) \right]
	\geq \E_\beta \left[ D_{V_{\mz}^n} (L_\iota^n, R_\iota^n) \right]$. Then we get by the lower bound on the expected distance between the boxes $L_\iota^n$ and $R_\iota^n$ \eqref{eq:set_to_set expectation} that for such a $k$
	\begin{align*}
		\E_\beta \left[ D_{B_n(\mz)}(\mz,v) \right] &
		\geq \frac{1}{3kd} \E_\beta \left[ D_{V_{\mz}^n} ( w_1 , w_2 ) \right]
		\geq
		\frac{1}{3kd}
		\E_\beta \left[ D_{V_{\mz}^n} (L_\iota^n, R_\iota^n) \right]\\
		&
		\overset{\eqref{eq:set_to_set expectation}}{\geq } \frac{1}{6kd } \E_\beta \left[ D_{V_{\mz}^n} (\mz , (n-1) \mo) \right] 
	\end{align*}
	which finishes the proof, as $v \in \mathcal{S}_n$ was assumed to minimize the expected distance $\E_\beta\left[D_{B_n(\mz)}(\mz,y)\right]$ among all vertices $y \in \mathcal{S}_n$.
\end{proof}

\begin{lemma}
	For all dimensions $d$ and all $\beta >0$, there exists an $\eta \in \left(0, \frac{1}{2}\right)$ such that  uniformly over all $n\in \N$ and all $x \in \mathcal{S}_n$
	\begin{align}\label{eq:box to boxonsphere expectation}
	\E_\beta \left[ D_{B_n(\mz)} \left(B_{\eta n}(\mz), B_{\eta n}(x) \right) \right] \geq \frac{c_1}{2} \Lambda(n,\beta)
	\end{align}
	where $c_1$ is the constant from \eqref{eq:point to pointonsphere} and there exists a constant $c_2$ such that
	\begin{align}\label{eq:box to boxonsphere probability}
	\p_\beta \left(  D_{B_n(\mz )} \left(B_{\eta n}(\mz), B_{\eta n}(x) \right)  \geq \frac{c_1}{4} \Lambda(n,\beta) \right)  \geq c_2 \text .
	\end{align}
	Furthermore, for each $\beta \geq 0$ there exist constants $c_3 > 0$ such that
	\begin{align}\label{eq:somequantiles of boxtobox}
	\p_\beta \left( D\left( B_n(\mz ) , B_{2n}(\mz )^C \right) \geq \frac{c_1}{4} \Lambda(n,\beta) \right) \geq c_3
	\end{align}
	uniformly over all $n \in \N$.
\end{lemma}

Note that in the above lemma, for $x\in \cS_n$ the box $B_{\eta n}(x)$ is not completely contained inside $B_n(\mz)$, but from the definition of $D_{B_n(\mz)} \left(\cdot, \cdot \right)$, we only consider the part that intersects $B_n(\mz)$.

\begin{proof}
	Given the results of \Cref{lem:point to pointonsphere}, the proof of \eqref{eq:box to boxonsphere expectation} and \eqref{eq:box to boxonsphere probability} works in the same way as the proof of Lemma \ref{lem:set_to_set} and we omit it. Regarding the statement of \eqref{eq:somequantiles of boxtobox}, we will first prove that for $\eta >0$ small enough
	\begin{align}\label{eq:somequantiles to show}
		\p_\beta \left(  D \left(B_{\eta n}(\mz), B_{ n}(\mz )^C \right)  \geq \frac{c_1}{4} \Lambda(n,\beta) \right) \geq c_4
	\end{align}
	for some constant $c_4>0$ and uniformly over all $n\in \N$. For this, we use the {\sl FKG inequality}, see \cite[Section 1.3]{heydenreich2017progress} or \cite{fortuin1971correlation,harris1960lower} for the original papers. We can cover the set $\bigcup_{x\in \mathcal{S}_n} B_{\eta n} (x)$ with uniformly (in $n$) finitely many sets of the form $B_{\eta n} (x)$. For example, we have
	\begin{align*}
		\bigcup_{\substack{x\in \mathcal{S}_n : \\ \langle x, e_1 \rangle = n }} B_{\eta n} (x) \subset \bigcup_{x \in F}  B_{\eta n} (x) 
	\end{align*}
	where $F=\left\{ne_1 + \sum_{i=2}^{d} k_i e_i  :  k_i \in \left\{ -  \big\lceil \frac{n}{\lceil \eta n \rceil} \big\rceil,\ldots, \big\lceil \frac{n}{\lceil \eta n \rceil} \big\rceil\right\} \text{ for all } i \in \{2,\ldots, d\} \right\}$, and all other faces of the set  $\bigcup_{x\in \mathcal{S}_n} B_{\eta n} (x)$ can be covered in a similar way. Suppose that $A^\prime_n \subset \mathcal{S}_n$ is a sequence of finite sets with $\sup_n |A_n^\prime|\eqqcolon A^\prime < \infty$ such that
	\begin{align*}
		\bigcup_{x\in \mathcal{S}_n} B_{\eta n} (x) = \bigcup_{x \in A_n^\prime}  B_{\eta n} (x) 
	\end{align*}
	for all $n\in \N$. So in particular we have
	\begin{align*}
		& \left\{D_{B_n(\mz)} \left(B_{\eta n}(\mz), B_{\eta n}(x) \right)  \geq \frac{c_1}{4} \Lambda(n,\beta) \text{ for all } x \in \mathcal{S}_n \right\}\\
		= 
		& \left\{D_{B_n(\mz)} \left(B_{\eta n}(\mz), B_{\eta n}(x) \right)  \geq \frac{c_1}{4} \Lambda(n,\beta) \text{ for all } x \in A_n^\prime \right\} \text .
	\end{align*}
	The events $\left\{D_{B_n(\mz)} \left(B_{\eta n}(\mz), B_{\eta n}(x) \right)  \geq \frac{c_1}{4} \Lambda(n,\beta)\right\}$ are decreasing for all $x \in \mathcal{S}_n$ in the sense that they are stable under the deletion of edges. Thus the FKG inequality and \eqref{eq:box to boxonsphere probability} imply that that
	\begin{align*}
		&\p_\beta \left( D_{B_n(\mz)} \left(B_{\eta n}(\mz), B_{\eta n}(x) \right)  \geq \frac{c_1}{4} \Lambda(n,\beta) \text{ for all } x \in \mathcal{S}_n \right)\\
		= & \p_\beta \left( D_{B_n(\mz)} \left(B_{\eta n}(\mz), B_{\eta n}(x) \right)  \geq \frac{c_1}{4} \Lambda(n,\beta) \text{ for all } x \in A_n^\prime \right) 
		\geq c_2^{|A_n^\prime|} \geq c_2^{A^\prime}.
	\end{align*}
	Assume that there is no direct edge from $\left[-(n-\eta n) , (n-\eta n) \right]^d$ to $\Z^d \setminus \left[-n,n\right]^d$. This has a uniform positive probability in $n$ and is also a decreasing event. Then any path from $B_{\eta n}(\mz)$ to $B_n(\mz)^C$ goes through at least one box $B_{\eta n }(x) \cap B_n(\mz)$ for some $x \in \mathcal{S}_n$. So with another application of the FKG inequality we get that
	\begin{align*}
		\p_\beta \left( D\left(B_{\eta n}(\mz) , B_n(\mz)^C\right) \geq \frac{c_1}{4} \Lambda(n,\beta) \right) \geq c_5
	\end{align*}
	for some $c_5>0$ and uniformly over all $n \in \N$.
	Next, let $A_n \subset B_n(\mz)$ be a sequence of sets such that $\bigcup_{x \in A_n} B_{\eta n} (x) = B_n(\mz)$ and $\sup_n |A_n| \eqqcolon \bar{A} < \infty$. Then $D\left(B_n(\mz),B_{2n}(\mz)^C\right) < \frac{c_1}{4} \Lambda(n,\beta)$ already implies that there exists a point $x \in A_n$ such that $ D\left(B_{\eta n}(x) , B_n(x)^C\right) < \frac{c_1}{4} \Lambda(n,\beta)$. By another application of the FKG inequality we have
	\begin{align*}
		& \p_\beta \left( D\left(B_n(\mz),B_{2n}(\mz)^C\right) \geq \frac{c_1}{4} \Lambda(n,\beta) \right) \\
		&
		\geq \p_\beta \left( D\left(B_{\eta n}(x) , B_n(x)^C\right) \geq \frac{c_1}{4} \Lambda(n,\beta) \text{ for all } x \in A_n \right) \geq c_5^{|A_n|} \geq c_5^{\bar{A}}
	\end{align*}
	which proves \eqref{eq:somequantiles of boxtobox}.
\end{proof}

\begin{lemma}\label{lem:all quantiles of pointtobox}
	For all $\beta \geq 0$ and all $\eps>0$, there exist $0<c<C<\infty$ such that
	\begin{align}\label{eq:all quantiles of pointtobox}
		\p_\beta \left( c\Lambda(n,\beta) \leq D \left( \mz,B_n(\mz)^C \right) \leq C\Lambda(n,\beta)\right) > 1-\eps
	\end{align}
	for all $n \in \N$.
\end{lemma}

Similar statements for one dimension and the continuous model were already proven in \cite{ding2013distances}. We follow a similar strategy here.

\begin{proof}
	By the inequality $D\left(\mz,B_n(\mz)^C\right)\leq D_{V_\mz^{n+2}}\left(\mz ,(n+1)\mo \right)$ we get that
	\begin{align*}
		\E_\beta\left[D\left(\mz,B_n(\mz)^C\right)\right] \leq \Lambda(n+2,\beta) \leq \Lambda(n,\beta) + 2 \text.
	\end{align*}
	Using Markov's inequality we see that
	\begin{align*}
		\p_\beta \left( D(\mz,B_n(\mz)^C) > C\Lambda(n,\beta)\right) \leq \frac{\Lambda(n,\beta) + 2}{C \Lambda(n,\beta)} \text,
	\end{align*}
	and thus the probability $\p_\beta \left( D(\mz,B_n(\mz)^C) \leq C\Lambda(n,\beta)\right)$ can be made arbitrarily close to $1$ by taking $C$ large enough. We will also refer to this case as the upper bound. The probability of the lower bound $\p_\beta \left( c \Lambda(n,\beta)\right) \leq D(\mz,B_n(\mz)^C)$ can be made arbitrarily close to $1$ for small $n$ by taking $c$ small enough. So we will always focus on $n$ large enough from here on. Fix $K, N \in \N_{>1}$ such that the function $i \mapsto \Lambda \left(K^{2i}N,\beta\right)$ is increasing in $i$. This is possible by Remark \ref{remark:Lambda grows}. We now consider boxes of the form $B_{K^{2i}N}(\mz)$. The probability of a direct edge from $B_{K^{2(i-1)} N}(\mz)$ to $B_{K^{2i}N}(\mz)^C$ equals the probability of a direct edge between $\mz$ and $B_{K^2}(\mz)^C$, and is by \eqref{eq:point to distance greater k} bounded by $\beta 50^d K^{-2}$. So the probability that there is some $i \in \{1,\ldots, K\}$ for which there is a direct edge from $B_{K^{2(i-1)} N}(\mz)$ to $B_{K^{2i}N}(\mz)^C$ is bounded by $\beta 50^d K^{-1}$. We denote the complement of this event by $\cA$. Conditioned on the event $\cA$, where there exists no edge between $B_{K^{2(i-1)} N}(\mz)$ and $B_{K^{2i}N}(\mz)^C$ for all $i\in \{1,\ldots, K\}$, each path from $B_{N}(\mz)$ to $B_{K^{2K}N}(\mz)^C$ needs to cross all the distances from $B_{K^{2(i-1)} N}(\mz)$ to $B_{2K^{2(i-1)}N}(\mz)^C$. For odd $i$, these distances are independent. As the function  $i \mapsto \Lambda \left(K^{2i}N,\beta\right)$ is increasing in $i$, conditioned on the event $\cA$, we have the bound
	\begin{align*}
		& \p_\beta  \left( D\left(\mz, B_{K^{2K} N}(\mz)^C\right) < \frac{c_1}{4} \Lambda(N, \beta)  \Big| \cA \right)\\
		&
		\leq \p_\beta \left( D\left(B_{K^{2(i-1)} N}(\mz), B_{2K^{2(i-1)} N}(\mz)^C\right) < \frac{c_1}{4} \Lambda\left( K^{2(i-1)} N, \beta\right)  \forall i \in \{1,\ldots, K\} \text{ odd}  \Big| \cA \right) \\
		& = \prod_{\substack{i=1:\\i \text{ odd} }}^{K} \p_\beta \left( D\left(B_{K^{2(i-1)} N}(\mz), B_{2K^{2(i-1)} N}(\mz)^C\right) < \frac{c_1}{4} \Lambda\left( K^{2(i-1)} N, \beta\right)   \Big| \cA \right) \\
		& \leq \prod_{\substack{i=1:\\i \text{ odd} }}^{K} \p_\beta \left( D\left(B_{K^{2(i-1)} N}(\mz), B_{2K^{2(i-1)} N}(\mz)^C\right) < \frac{c_1}{4} \Lambda\left( K^{2(i-1)} N, \beta\right)  \right)  \leq (1-c_3)^{\lfloor\frac{K}{2}\rfloor} \text,
	\end{align*}
	where the second last inequality holds because of FKG, as events of the form $\left\{D(\cdot,\cdot) < x\right\}$ are increasing and $\cA$ is decreasing, and 
	where $c_3$ is the constant from \eqref{eq:somequantiles of boxtobox}. Thus we have that
	\begin{align*}
		& \p_\beta  \left( D\left(\mz, B_{K^{2K} N}(\mz)^C\right) < \frac{c_1}{4} \Lambda(N, \beta)  \right)\\
		&
		\leq \p_\beta  \left( D\left(\mz, B_{K^{2K} N}(\mz)^C\right) < \frac{c_1}{4} \Lambda(N, \beta)  \Big| \cA \right) + \p_\beta \left( \cA^C  \right)
		\leq (1-c_3)^{\lfloor\frac{K}{2}\rfloor} + \beta 50^d K^{-1}
	\end{align*}
	and this quantity can be made arbitrary small by suitable choice of $K$. To finish the proof, remember that $\Lambda(N, \beta)$ and $\Lambda\left( K^{2K} N, \beta \right) $ are off by a multiplicative factor of at most $K^{2K}$, as 
	\begin{equation*}
		\Lambda\left(  N, \beta \right) \leq \Lambda\left( K^{2K} N, \beta \right) \overset{\eqref{eq:Lambda}}{\leq} \Lambda\left( K^{2K} , \beta \right) \Lambda\left(  N, \beta \right)
		\leq  K^{2K}  \Lambda\left(  N, \beta \right) \text. 
	\end{equation*}
	Thus we have
	\begin{align*}
		\p_\beta  &\left( D\left(\mz, B_{K^{2K} N}(\mz)^C\right) < \frac{c_1}{4 K^{2K}} \Lambda(K^{2K} N, \beta)  \right) \leq \p_\beta  \left( D\left(\mz, B_{K^{2K} N}(\mz)^C\right) < \frac{c_1}{4 } \Lambda( N, \beta)  \right) 
		\\
		&
		\leq (1-c_3)^{\lfloor\frac{K}{2}\rfloor} + \beta 50^d K^{-1}.
	\end{align*}
	Now, for fixed $\eps>0$, take $K$ large enough so that $(1-c_3)^{\lfloor\frac{K}{2}\rfloor} + \beta 50^d K^{-1} < \eps$. For $n \in \N$ large enough with $n > K^{2K}$, let $N$ be the largest integer for which $K^{2K}N \leq n$. We know that $K^{2K}N \leq n \leq K^{2K}2N$ and this also yields, by Proposition \ref{propo:scaling}, that
	\begin{align*}
		\Lambda(n,\beta) \leq 3 \Lambda \left( K^{2K}2 N , \beta  \right)
		\leq 6 \Lambda \left( K^{2K} N , \beta  \right)
	\end{align*}
	which already implies
	\begin{align*}
	 & \p_\beta  \left( D\left(\mz, B_{n}(\mz)^C\right) < \frac{c_1}{24 K^{2K}} \Lambda(n, \beta)  \right)
	\leq
	\p_\beta  \left( D\left(\mz, B_{n}(\mz)^C\right) < \frac{c_1}{4 K^{2K}} \Lambda(K^{2K}N, \beta)  \right)\\
	&
	\leq
	\p_\beta  \left( D\left(\mz, B_{K^{2K}N}(\mz)^C\right) < \frac{c_1}{4 K^{2K}} \Lambda(K^{2K}N, \beta)  \right) \leq (1-c_3)^{\lfloor\frac{K}{2}\rfloor} + \beta 50^d K^{-1} \leq \eps \text .
	\end{align*}
\end{proof}

The previous lemma tells us that for fixed $\beta > 0$ all quantiles of $D\left(\mz,B_n(\mz)^C\right)$ are of order $\Lambda(n,\beta)$. We want to prove a similar statement for the quantiles of $D\left(B_n(\mz),B_{2n}(\mz)^C\right)$. However, an analogous statement can not be true, as there is a uniform positive probability of a direct edge between $B_n(\mz)$ and $B_{2n}(\mz)^C$. But if we condition on the event that there is no such direct edge, the statement still holds.

\begin{lemma}\label{lem:all quantiles of boxtobox}
	Let $\mathcal{L}$ be the event that there is no direct edge between $B_n(\mz)$ and $B_{2n}(\mz)^C$.	For all $\beta > 0$ and all $\eps>0$, there exist $0<c<C<\infty$ such that
	\begin{align}\label{eq:all quantiles of boxtobox}
	\p_\beta \left( c\Lambda(n,\beta) \leq D(B_n(\mz),B_{2n}(\mz)^C) \leq C\Lambda(n,\beta) \ \big| \ \mathcal{L} \right) > 1-\eps
	\end{align}
	for all $n \in \N$.
\end{lemma}
\begin{proof}
	We know that
	\begin{align*}
		& \E_\beta \left[ D(B_n(\mz),B_{2n}(\mz)^C) \right] \leq \E_\beta \left[ D (n\mo,(2n+1)\mo) \right]\\
		&
		\leq \E_\beta \left[ D_{V_\mo^{n}}(n\mo,(2n-1)\mo) \right] + 2
		\leq \Lambda(n,\beta) + 1 \text ,
	\end{align*}
	and thus the probability $\p_\beta \left(  D(B_n(\mz),B_{2n}(\mz)^C) \leq C\Lambda(n,\beta) \ \big| \ \mathcal{L} \right)$ can be made arbitrarily close to $1$ by taking $C$ large enough, by Markov's inequality.
	For the lower bound, we first consider integers of the form $N_k=M^kN_0$, where we fix $M \in \N$ first. Let $M$ be the smallest natural number such that $M \geq 100$ and $\Lambda(M,\beta) \leq \frac{M}{10}$. The inequality $\Lambda(M,\beta) \leq \frac{M}{10}$ holds for large enough $M$, as $\Lambda(M,\beta)$ asymptotically grows like a power of $M$ that is strictly less than one, see Section \ref{sec:submultiplicativity}. As $\beta$ is fixed for the rest of the proof, we simply write $\Lambda(n)$ for $\Lambda(n,\beta)$.
	We write $C_n$ for the annulus $B_{2n}(\mz)\setminus B_n(\mz)$. Let $\mathcal{A}_\delta$ denote the event that for all vertices $x \in C_{MN}$ for which there exists an edge $e=\{x,y\}$ with $\|x-y\|\infty \geq N$ one has $D\left(x,C_{MN}^C; \omega^{e-}\right) \geq \delta \Lambda(MN)$. We will now show that the probability of the event $\mathcal{A}_\delta$ converges to $1$ as $\delta \to 0$. Remember that $\Lambda(MN)$ and $\Lambda(N)$ differ by a factor of at most $M$. Let us first consider the event that for some $\delta_1 > 0$ there exists a vertex incident to a long edge in one of the boundary regions of thickness $\delta_1 N$ of $C_{MN}$, Formally, for $\delta_1 \in \left(0,\tfrac{1}{2}\right)$, we define the boundary region $\partial^{\delta_1} C_{MN}$ of $C_{MN}$ by
	\begin{equation*}
		\partial^{\delta_1} C_{MN} = \left\{B_{MN+\delta_1 N}(\mz)\setminus B_{MN}(\mz)\right\} \cup \left\{B_{2MN}(\mz) \setminus B_{2MN-\delta_1 N}(\mz)\right\} \text .
	\end{equation*}
	The set $\partial^{\delta_1} C_{MN}$ has a size of at most $4d \delta_1 N \left(5MN\right)^{d-1}$, as for $x\in \partial^{\delta_1} C_{MN}$ one needs to fix one of the coordinates within the interval $\left(MN, MN+\delta_1 N\right]$, respectively in the interval $\left( 2MN-\delta_1 N , 2MN  \right]$, or one of the reflections of these intervals, and then has at most $4MN+1$ possibilities for each of the remaining $d-1$ coordinates. Combining this gives
	\begin{align*}
		\left| \partial^{\delta_1} C_{MN} \right|
		\leq (4d \delta_1 N) (4MN+1)^{d-1} \leq 4d \delta_1 N (5MN)^{d-1}.
	\end{align*}
	The probability that a vertex is incident to some edge of length $\geq N$ is proportional to $\frac{\beta}{N^d}$ as shown in \eqref{eq:point to distance greater k}. So together with \eqref{eq:point to distance greater k} we get that
	\begin{align*}
		\p_\beta \left( \exists x \in \partial^{\delta_1} C_{MN} , y \in B_{N-1}(x)^C : x \sim y  \right) 
		&
		\leq 4d \delta_1 N (5MN)^{d-1} \p_\beta \left(\mz \sim \cS_{\geq N}\right)
		\\
		&
		\leq \delta_1 \cdot 4d (5MN)^d \beta 50^d N^{-d} \leq \delta_1 \cdot \beta \left(10^3M\right)^d.
	\end{align*}
	Furthermore, the expected number of points $x \in C_{MN}$ which are incident to a long edge is bounded by
	\begin{align}\label{eq:gleich bound expectation}
		&\notag \E_\beta \left[ \left|\{ x \in C_{MN} : x \sim B_{N-1}(x)^C \} \right| \right] 
		\leq
		\sum_{x \in C_{MN}} \sum_{y \in B_{N-1}(x)^C} \p_\beta (x\sim y)\\
		&
		\leq
		|C_{MN}| \sum_{y \in B_{N-1}(\mz)^C} \p_\beta (\mz \sim y)
		\overset{\eqref{eq:point to distance greater k helpful bound}}{\leq} (5MN)^d \beta 50^d N^{-d} \leq\beta (250M)^d.
	\end{align}
	where the second last inequality holds as $|C_{MN}| \leq (4MN+1)^d \leq (5MN)^d$, and because the sum $\sum_{y \in B_{N-1}(\mz)^C} \p_\beta (\mz \sim y)$ can be upper bounded by $\beta 50^d N^{-d}$ in the exact same way as in \eqref{eq:point to distance greater k helpful bound}.
	As the existence of an edge $\{x,y\}$ with $\|x-y\|_\infty\geq N$ and the distance $D\left( x,B_{\delta_1 N}(x)^C; \omega^{\{x,y\}-} \right)$ are independent random variables, we get with a union bound that
	\begin{align*}
		& \p_\beta\left( \exists x \in C_{MN} , y \in B_{N-1}(x)^C : x \sim y , D\left( x,B_{\delta_1 N}(x)^C; \omega^{\{x,y\}-} \right) < \delta \Lambda(MN) \right)\\
		& \leq \sum_{x \in C_{MN}}  \sum_{y \in B_{N-1}(x)^C} \p_\beta (x\sim y) 
		\p_\beta \left( D\left( x, B_{\delta_1 N}(x)^C \right) < \delta \Lambda(MN) \right)\\
		& \leq \beta (250M)^d \p_\beta \left( D\left( \mz , B_{\delta_1 N}(\mz)^C \right) < \delta \Lambda(MN) \right)
	\end{align*}
	where we used \eqref{eq:gleich bound expectation} for the last inequality. 
	Thus we also get that
	\begin{align}
		\notag &\p_\beta \left( \mathcal{A}_\delta^C \right) = \p_\beta 
		\left( \exists x \in C_{MN}, y \in B_{N-1}(x)^C : x \sim y , D\left( x, C_{MN}^C ; \omega^{\{x,y\}-} \right) < \delta \Lambda(MN)  \right)\\
		\notag & \leq \p_\beta \left( \exists x \in \partial^{\delta_1} C_{MN} , y \in B_{N-1}(x)^C : x \sim y  \right)\\
		\notag & \ \
		+ \p_\beta\left( \exists x \in C_{MN} , y \in B_{N-1}(x)^C : x \sim y , D\left( x,B_{\delta_1 N}(x)^C; \omega^{\{x,y\}-} \right) < \delta \Lambda(MN) \right)\\
		\label{eq:probability of long bad edge}& \leq \delta_1 \beta \left(10^3M\right)^d + \beta (250M)^d \p_\beta \left( D\left( \mz, B_{\delta_1 N}(\mz)^C \right) < \delta \Lambda(MN) \right)
	\end{align}
	and this converges to $0$ as $\delta \to 0$, for an appropriate choice of $\delta_1(\delta)$, by Lemma \ref{lem:all quantiles of pointtobox} uniformly over $N \in \N$. We write $f(\delta)$ for the supremum of $\p_\beta\left(\mathcal{A}_\delta^C\right)$ over all $N\in \N$ and for $A,B \subset V$, we write $D^\star\left(A, B\right)$ for the indirect distance between $A$ and $B$, i.e.,  the length of the shortest path between $A$ and $B$ that does not use a direct edge between $A$ and $B$. Now assume that $D^\star\left(B_{MN}(\mz), B_{2MN}(\mz)^C\right) < \delta \Lambda(MN)$. We now consider the path between $B_{MN}(\mz) $ and $ B_{2MN}(\mz)^C$ that achieves this distance. Either this path uses some long edge (of length greater than $N-1$), or it only jumps from one block of the form $V_v^N$ to directly neighboring blocks. The probability that there exists a point $x \in C_{MN}$ and a long edge $e$ incident to it such that $D\left(x,C_{MN}^C; \omega^{e-}\right) < \delta \Lambda(MN)$ is relatively small by \eqref{eq:probability of long bad edge}. Any path that does not use long edges can only do jumps between neighboring blocks of the form $V_v^N$. Say the geodesic between $B_{MN}(\mz)$ and $B_{2MN}(\mz)^C$ uses the blocks $\left(V_{v_i^\prime}^N\right)_{i=0}^{L^\prime}$. Consider the loop-erased trace of this walk on the blocks, i.e., say that the path uses the blocks $\left(V_{v_i}^N\right)_{i=0}^{L} \subset C_{MN}$ with $\|v_i-v_{i-1}\|_\infty=1$ and never returns to $V_{v_i}^N$ after going to $V_{v_{i+1}}^N$. There need to be at least $\frac{M}{3}$ transitions between blocks $V_{v_{j_i}}^N$ and $V_{v_{j_{i+1}}}^N$ with $\|v_{j_{i}}-v_{j_{i+1}}\|_\infty = 2$ and $0 \leq j_{1} < j_{2} < \ldots < j_{M/3} \leq L$ as the path needs to walk a distance in the infinity-norm of at least $MN$. So in particular we have
	\begin{align*}
		\sum_{i=1}^{\lceil M/3 \rceil} & D^\star \left( V_{v_{j_{i}}}^N, V_{v_{j_{i+1}}}^N \right) \leq D^\star \left(B_{MN}(\mz), B_{2MN}(\mz)^C \right) \\
		& < \delta \Lambda(MN,\beta) \leq \delta \Lambda(M,\beta) \Lambda(N,\beta) \leq \frac{M}{10} \delta \Lambda(N,\beta)
	\end{align*}
	where we used our assumption on $M$ for the last step. So in particular there need to be at least two transitions between $V_{v_{j_{i}}}^N$ and $V_{v_{j_{i+1}}}^N$ that satisfy $D^\star \left( V_{v_{j_{i}}}^N, V_{v_{j_{i+1}}}^N \right) < \delta \Lambda(N,\beta)$. In fact, there need to be some linear number in $M$ many such transitions, but two are sufficient for our purposes here. All these transitions need to be disjoint, as shortest paths never use the same edge twice. Thus we get by the {\sl BK inequality} (see \cite[Section 1.3]{heydenreich2017progress} or \cite{van1985inequalities,reimer2000proof}) that
	\begin{align*}
		\p_\beta \left(\sum_{i=1}^{\lceil M/3 \rceil} D^\star \left( V_{v_{j_{i}}}^N, V_{v_{j_{i+1}}}^N \right) \leq \frac{M}{10} \delta \Lambda(N) \right) \leq M^2 \left( \min_i \p_\beta \left( D^\star \left( V_{v_{j_{i}}}^N, V_{v_{j_{i+1}}}^N \right) < \delta \Lambda(N) \right)\right)^2.
	\end{align*}
	For each combination of vectors $v_{j_{i}},v_{j_{i+1}}$ with $\|v_{j_{i}}-v_{j_{i+1}}\|_\infty=2$, we can translate and rotate the boxes $V_{v_{j_{i}}}^N$ and $V_{v_{j_{i+1}}}^N$ to boxes $T\left(V_{v_{j_{i}}}^N\right)$ and $T\left(V_{v_{j_{i+1}}}^N\right)$ in such a way that  $T\left(V_{v_{j_{i}}}^N\right) \subset B_N(\mz)$ and $T\left(V_{v_{j_{i+1}}}^N\right) \subset B_{2N}(\mz)^C$. By translational and rotational invariance of our long-range percolation model, this already implies that
	\begin{align*}
		& \min_i \p_\beta \left( D^\star \left( V_{v_{j_{i}}}^N, V_{v_{j_{i+1}}}^N \right) < \delta \Lambda(N) \right)  \leq \p_\beta\left( D^\star \left(B_N(\mz),B_{2N}(\mz)^C\right) < \delta \Lambda(N) \right) \text .
	\end{align*}
	There are at most $\left((5M)^d\right)!$ choices for possible choice of vertices $v_0,v_1,\ldots, v_L$, as there are at most $(5M)^d$ possibilities for $v_0$ and $(5M)^d-1$ possibilities for $v_1$ and so on. Overall we see that the probability that there exists an indirect path between $B_{MN}(\mz)$ and $B_{2MN}(\mz)^C$ of length $\delta \Lambda(MN)$, which jumps between neighboring blocks of the form $V_{v}^N$ only, is bounded by 
	\begin{align*}
		\left(5M^d\right)!M^2 \p_\beta \left(D^\star \left(B_N(\mz),B_{2N}(\mz)^C\right) < \delta \Lambda(N)\right)^2 .
	\end{align*}
	We write $S$ for the constant $\left(5M^d\right)!M^2$. Thus we get that
	\begin{align*}
		\p_\beta & \left(D^\star \left(B_{MN}(\mz),B_{2MN}(\mz)^C < \delta \Lambda(MN)\right)\right) \\
		& \leq
		S \p_\beta \left( D^\star \left(B_N(\mz),B_{2N}(\mz)^C\right) < \delta \Lambda(N)\right)^2 + \p_\beta \left(\mathcal{A}_\delta^C\right) \text .
	\end{align*}
	We define the sequence $(a_n)_{n\in \N}$ by 
	\begin{equation*}
		a_0 = \p_\beta \left(D^\star \left(B_N(\mz), B_{2N}(\mz)^C\right) < \delta \Lambda(N)\right)
	\end{equation*}
	and $a_{n+1} = S a_n^2 + f(\delta)$. Inductively it follows that
	\begin{align*}
			\p_\beta \left( D^\star \left(B_{M^kN}(\mz),B_{2M^kN}(\mz)^C\right) < \delta \Lambda(M^kN)\right) \leq a_k 
	\end{align*}
	for all $k\in \N$.
	For $f(\delta)< \frac{1}{4S}$, the equation $a= S a^2 +f(\delta)$ has the two solutions
	\begin{align*}
		a_- = \frac{1 - \sqrt{1-4Sf(\delta)}}{2S} \text{ and } a_+ = \frac{1+ \sqrt{1-4Sf(\delta)}}{2S} > \frac{1}{2S} \text .
	\end{align*}
	For $a_0 \in \left[0, a_+\right)$, and thus in particular for $a_0 \in \left[0, \frac{1}{2S}\right]$, the sequence $a_n$  converges to $a_- = \frac{1-\sqrt{1-4Sf(\delta)}}{2S} \approx f(\delta)$ and thus we get
	\begin{align*}
		\limsup_{k\to \infty} \p_\beta \left( D^\star \left(B_{M^kN}(\mz),B_{2M^kN}(\mz)^C\right) < \delta \Lambda(M^kN)\right) \leq a_- \ .
	\end{align*}
	For fixed $N\in \N$, the requirement
	\begin{align*}
		a_0 = \p_\beta \left( D^\star \left(B_N(\mz),B_{2N}(\mz)^C\right) < \delta \Lambda(N)\right) \leq \frac{1}{2S}
	\end{align*}
	is satisfied for small enough $\delta >0$ and this shows \eqref{eq:all quantiles of boxtobox} along the subsequence $N_k=M^kN$. To get the statement for all integer numbers, one can use Proposition \ref{propo:scaling} and the fact that $\Lambda(n) \leq \Lambda(mn) \leq m \Lambda(n)$ for all integers $m,n$.
\end{proof}

With the same technique as above one can also prove that the indirect distance between $V_\mz^n$ and the set $B_n\left(V_\mz^n\right)^C= \left\{x\in \Z^d : D_\infty(x,V_\mz^n) > n\right\} = \bigcup_{u \in \Z^d : \| u \|_\infty \geq 2} V_u^n$ scales like $\Lambda(n,\beta)$.

\begin{corollary}\label{coro:quantiles of boxes}
	For all $\beta \geq 0$ and  $\eps >0$ there exist $0<c_\eps < C_\eps <\infty$ such that
	\begin{align*}
		\p_\beta \left( c_\eps \Lambda(n,\beta) \leq D^\star \left( V_\mz^n , B_n\left(V_\mz^n\right)^C \right) \leq C_\eps \Lambda(n,\beta) \right)   \geq 1-\eps \ .
	\end{align*}
\end{corollary}

\section{The proof of Theorem \ref{theo:exponent}}\label{sec:prooftheo1}

We first give an outline of the proof of Theorem \ref{theo:exponent}. In Lemma \ref{lem:all quantiles of pointtobox}, we showed that $D \left(\mz,B_n\left(\mz\right)^C\right) \approx_P \Lambda(n,\beta)$, and Lemma \ref{lem:endisfar} shows that $\Lambda(n,\beta) \approx \E_\beta\left[D_{V_\mz^n}(\mz,(n-1)e_1)\right]$, meaning that the ratio of these two expressions is uniformly bounded from below and above by constants $0<c<C<\infty$. In Lemma \ref{lem:supermultiplicativity} below we prove supermultiplicativity of $\Lambda(n,\beta)$. Together with the submultiplicativity proven in Lemma \ref{lem:submultiplicativity} this shows that for each $\beta \geq 0$ there exists $c_\beta>0$ such that $c_\beta\Lambda(m,\beta)\Lambda(n,\beta) \leq \Lambda(mn,\beta) \leq \Lambda(m,\beta)\Lambda(n,\beta)$. We define $a_k=\log\left(\Lambda(2^k,\beta)\right)$. The sequence is subadditive and thus
\begin{align*}
\dxp(\beta) = \lim_{k\to \infty} \frac{\log(\Lambda(2^k,\beta))}{\log(2^k)} = \lim_{k\to \infty} \frac{a_k}{\log(2) k} = \inf_{k \in \N} \frac{a_k}{\log(2) k}
\end{align*}
exists, where the last inequality holds because of Fekete's Lemma. On the other hand, the sequence $b_k = \log(c_\beta \Lambda(2^k,\beta))$ satisfies
\begin{align*}
b_{k+l} = \log(c_\beta \Lambda(2^{k+l},\beta)) \geq \log(c_\beta \Lambda(2^{k},\beta)c_\beta \Lambda(2^{l},\beta)) = b_k+b_l
\end{align*}
and thus
\begin{align*}
\dxp(\beta)=\lim_{k\to \infty} \frac{ \log(c_\beta\Lambda(2^k,\beta))}{\log(2^k)} = \lim_{k\to \infty} \frac{b_k}{\log(2) k} = \sup_{k \in \N} \frac{b_k}{\log(2) k} \text .
\end{align*}
This already implies that
\begin{align*}
2^{k\dxp(\beta)} \leq \Lambda(2^k,\beta) \leq c_\beta^{-1} 2^{k\dxp(\beta)}
\end{align*}
for all $k \in \N$. These two inequalities can be extended from points of the form $2^k$ to all integers with \Cref{propo:scaling}. So there exists a constant $0<C_\beta^\prime <\infty$ such that for all $n\in \N$
\begin{align*}
\frac{1}{C_\beta^\prime } n^{\dxp(\beta)} \leq \Lambda(n,\beta) \leq C_\beta^\prime n^{\dxp(\beta)}.
\end{align*}
which shows \eqref{eq:theo:exponent1}. So we still need to prove supermultiplicativity of $\Lambda(\cdot,\beta)$ in order to prove the first item in Theorem \ref{theo:exponent}. The second item of Theorem \ref{theo:exponent}, i.e.,  the bounds on the diameter of cubes \eqref{eq:theo:exponent2}, we show in Section \ref{sec:dia}.

\subsection{Distances between certain points}\label{sec:distance of long edge points}

In this chapter, we examine the typical behavior of distances between points that are connected to long edges. In \Cref{lem:linearspacing2}, we consider the infinity distance between such points. Using a coupling argument with the continuous model, we compare the situation to the situations occurring in  \Cref{lem:separation both far} and \Cref{lem:separation one far}. Then, in  \Cref{lem:graphspacing} we translate these bounds on the infinity distance into bounds on the typical graph distance between points that are incident to long edges.

Fix the three blocks $V_u^n, V_w^n$ and $V_{\mz}^n$ with $\|u\|_\infty \geq 2$. The next lemma deals with the infinity distance between points $x,y \in V_{\mz}^n$ with $x \sim V_u^n, y\sim V_w^n$. 
\begin{lemma}\label{lem:linearspacing2}
	For all $\frac{1}{n} <\eps \leq \frac{1}{4}$ and $u,w \in \Z^d \setminus \{0\}$ with $\|u\|_\infty \geq 2$ one has
	\begin{align*}
		\p_\beta \big( \exists x,y \in V_{\mz}^n : \|x-y\|_\infty \leq \eps n, x \sim V_u^n, y \sim V_w^n \ \big|  \ V_\mz^n \sim V_u^n, V_\mz^n \sim V_w^n \big) \leq C_d^\prime \eps^{1/2} \lceil\beta \rceil^2
	\end{align*}
	where $C_d^\prime$ is a constant that depends only on the dimension $d$.
\end{lemma}
\begin{proof}
	Let $\e$ be the symmetrized Poisson process constructed in subSection \ref{subsec:cts} about the continuous model, i.e.,  $\tilde{\e}$ is a Poisson process on $\R^d\times \R^d$ with intensity $\frac{\beta}{2\|t-s\|^{2d}}$ and $\e$ is defined by $\e \coloneqq \left\{(s,t)\in \R^d \times \R^d : (t,s) \in \tilde{\e}\right\} \cup \tilde{\e}$. Now we place an edge between $x,y \in \Z^d$ if and only if 
	\begin{equation*}
		(x+\cC) \times (y+\cC) \cap n\e \neq \emptyset
	\end{equation*}
	and call this graph $G=(V,E)$. The distribution of the resulting graph is identical to $\p_\beta$ by the dilation invariance of $\e$. We can do the same procedure for $\lfloor \frac{1}{2\eps}-1 \rfloor \e$, i.e., place an edge between $x^\prime,y^\prime \in \Z^d$ if and only if 
	\begin{equation*}
		(x^\prime+\cC) \times (y^\prime+\cC) \cap \Big\lfloor \frac{1}{2\eps} - 1 \Big\rfloor \e \neq \emptyset
	\end{equation*}
	and call the resulting graph $G^\prime =( V^\prime, E^\prime)$. Now assume that in the graph $G$ there exist $x,y \in V_{\mz}^n$ with $\|x-y\|_\infty \leq \eps n$ such that $x\sim V_u^n$ and $y\sim V_w^n$ in $G$. Then there exist 
	\begin{equation*}
		x_c \in x+\cC, u_c \in nu+\left[0,n\right)^d , y_c\in y+\cC, w_c \in nw+\left[0,n\right)^d
	\end{equation*}
	with $(x_c,u_c), (y_c, w_c) \in n\e$. We also have $\|x_c-y_c\|_\infty \leq \eps n+1 \leq 2\eps n$. Now we rescale the process from size $n$ to size $\lfloor \frac{1}{2\eps} - 1 \rfloor$. For $(\tilde{x_c},\tilde{u_c}) =  \frac{\lfloor \frac{1}{2\eps}-1\rfloor}{n} (x_c , u_c) $ and $(\tilde{y_c},\tilde{w_c}) =  \frac{\lfloor \frac{1}{2\eps}-1\rfloor}{n} (y_c , w_c) $ we have 
	\begin{align*}
	& (\tilde{x_c},\tilde{u_c}) \in \left(\left(\frac{\lfloor \frac{1}{2\eps}-1 \rfloor}{n} x + \left[0, \frac{\lfloor \frac{1}{2\eps}-1 \rfloor}{n} \right)^d \right) 
	\times 
	\left(\Big\lfloor\frac{1}{2\eps}-1 \Big\rfloor u + \left[0,\Big\lfloor\frac{1}{2\eps}-1\Big\rfloor\right)^d \right)\right)
	\cap \Big\lfloor  \frac{1}{2\eps} - 1 \Big\rfloor\mathcal{E} \text , \\
	&
	(\tilde{y_c},\tilde{w_c}) \in \left(\left(\frac{\lfloor \frac{1}{2\eps}-1 \rfloor}{n} y +  \left[0, \frac{\lfloor \frac{1}{2\eps}-1 \rfloor}{n} \right)^d \right) 
	\times 
	\left(\Big\lfloor\frac{1}{2\eps}-1\Big\rfloor w + \left[0,\Big\lfloor\frac{1}{2\eps}-1\Big\rfloor\right)^d  \right)\right)
	\cap \Big\lfloor \frac{1}{2\eps} - 1 \Big\rfloor \mathcal{E} \text .
	\end{align*}
	From the rescaling we also have $\|\tilde{x_c}-\tilde{y_c}\|_\infty \leq 2\eps \lfloor \frac{1}{2\eps}-1\rfloor < 1 $. So in particular there are vertices $x^\prime, y^\prime \in \left\{0, \ldots ,\lfloor \frac{1}{2\eps}-1\rfloor - 1\right\}^d$ with $x^\prime \sim V_u^{\lfloor \frac{1}{2\eps}-1\rfloor}, y^\prime \sim V_w^{\lfloor \frac{1}{2\eps}-1\rfloor}$ in $G^\prime$, and $\|x^\prime -y^\prime \|_\infty \leq 1$. Write $N= \lfloor \frac{1}{2 \eps} - 1 \rfloor$. From Lemmas \ref{lem:separation both far} and \ref{lem:separation one far} we get
	\begin{align*}
	\p_\beta & \Big( \exists x , y \in V_{\mz}^{N} : \|x-y\|_\infty\leq 1 ,  x\sim V_u^{N}, y \sim V_w^{N} \ 
	\Big| \ V_{\mz}^{N}\sim V_u^{N}, V_{\mz}^{N} \sim V_w^{N}\Big) \\
	& \leq \frac{C_d \lceil \beta \rceil^2}{N^{1/2}}
	= \frac{C_d \lceil \beta \rceil^2}{\lfloor \frac{1}{2\eps}-1\rfloor^{1/2}}
	\leq C_d^\prime \eps^{1/2} \lceil \beta \rceil^2
	\end{align*}
	for some constants $C_d, C_d^\prime <\infty$. With the coupling argument from before we thus also get
	\begin{align*}
	\p_\beta & \left(\exists x,y \in V_{\mz}^n : \|x-y\|_\infty \leq \eps n , x \sim V_u^{n}, y \sim V_w^{n} \  \big| \ V_{\mz}^{n} \sim V_u^{n} , V_{\mz}^{n} \sim V_w^{n}  \right) \\
\leq 	\p_\beta & \Big( \exists x , y \in V_{\mz}^{N} : \|x-y\|_\infty\leq 1 , x\sim V_u^{N}, y \sim V_w^{N}  \ 
	\Big| \ V_{\mz}^{N}\sim V_u^{N}, V_{\mz}^{N} \sim V_w^{N}\Big)  \leq C_d^\prime \eps^{1/2} \lceil \beta \rceil^2 .
	\end{align*}
\end{proof}

\begin{lemma}\label{lem:graphspacing}
	For all dimensions $d$ and all $\beta \geq 0$, there exists a function $g_1(\eps)$ with $g_1(\eps) \underset{\eps \to 0}{\longrightarrow} 1$ such that for all $u,w \in \Z^d \setminus \{\mz\}$ with $\|u\|_\infty \geq 2$ and all large enough $n\geq n(\eps)$
	\begin{align*}
		\p_\beta \big(  D_{V_{\mz}^n}(x,y) > \eps \Lambda(n,\beta) \text{ for all } x,y \in V_{\mz}^n \text{ with } x\sim V_u^n, y\sim V_w^n \ \big| \ V_u^n \sim V_\mz^n \sim V_w^n \big) \geq g_1(\eps) \text .
	\end{align*}
\end{lemma}
\begin{proof}
	We write $\p_\beta^{u,w}\left(\cdot\right)$ for the conditional  probability measure $\p_\beta \big( \cdot \ \big| \ V_u^n \sim V_\mz^n \sim V_w^n \big)$. As $\beta$ is fixed throughout the rest of the proof, we write $\Lambda(n)$ for $\Lambda(n,\beta)$. We define the event 
	\begin{align*}
		& \mathcal{A}(K,\eps_1,\eps)  = 
		\left\{\|x-y\|_\infty > \eps_1n \text{ for all } x,y \in V_{\mz}^n \text{ with } x\sim V_u^n, y\sim V_w^n \right\}\\
		&
		\cap
		\left\{D_{V_{\mz}^n}\left(x,B_{\eps_1 n}(x)^C\right) > \eps \Lambda(n) \text{ for all } x \in V_{\mz}^n \text{ with } x\sim V_u^n \right\}
		\cap
		\left\{ \left| \{x\in V_{\mz}^n : x\sim V_u^n \} \right| \leq K \right\}
	\end{align*}
	and observe that 
	\begin{equation*}
		\left\{  D_{V_{\mz}^n}(x,y) > \eps \Lambda(n) \text{ for all } x,y \in V_{\mz}^n \text{ with } x\sim V_u^n, y\sim V_w^n \right\} \supset  \mathcal{A}(K,\eps_1,\eps) \text .
	\end{equation*}
	Thus it suffices to show that $\p_\beta^{u,w}\left( \mathcal{A}(K,\eps_1,\eps) \right)$ converges to $1$ as $\eps \to 0$ for an appropriate choice of $K = K(\eps), \eps_1 = \eps_1(\eps)$. Respectively, that $\p_\beta^{u,w}\left( \mathcal{A}(K,\eps_1,\eps)^C \right)$ converges to $0$. We have that
	\begin{align*}
	& \mathcal{A}(K,\eps_1,\eps)^C  = \left\{ \left| \{x\in V_{\mz}^n : x\sim V_u^n \} \right| > K \right\}\\
	&
	\cup
	\left\{\|x-y\|_\infty < \eps_1n \text{ for some } x,y \in V_{\mz}^n \text{ with } x\sim V_u^n, y\sim V_w^n \right\}\\
	&
	\cup
	\Big(\left\{D_{V_{\mz}^n}\left(x,B_{\eps_1 n}(x)^C\right) \leq \eps \Lambda(n) \text{ for some } x \in V_{\mz}^n \text{ with } x\sim V_u^n \right\} \cap \left\{ \left| \{x\in V_{\mz}^n : x\sim V_u^n \} \right| \leq K  \right\} \Big)
	\end{align*}
	and thus we get with Lemma \ref{lem:linearspacing2} that
	\begin{align*}
		& \p_\beta^{u,w}\left(\mathcal{A}(K,\eps_1,\eps)^C\right)  \leq  \p_\beta^{u,w} \left(  \left| \{x\in V_{\mz}^n : x\sim V_u^n \} \right| > K \right) + C_d^\prime \eps_1^{1/2} \lceil \beta \rceil^2 \\
		&
		+
		\p_\beta^{u,w} \left(\left\{ D_{V_{\mz}^n}\left(x,B_{\eps_1 n}(x)^C\right) \leq \eps \Lambda(n) \text{ for some } x \in V_{\mz}^n \text{ with } x\sim V_u^n \right\} \cap \left\{ \left| \{x\in V_{\mz}^n : x\sim V_u^n \} \right| \leq K  \right\} \right)  \\
		& \leq \p_\beta^{u,w} \left(  \left| \{x\in V_{\mz}^n : x\sim V_u^n \} \right| > K \right) + C_d^\prime \eps_1^{1/2} \lceil \beta \rceil^2 + K \p_\beta \left( D \left(\mz,B_{\eps_1 n}(\mz)^C\right) \leq \eps \Lambda(n) \right) \text .
	\end{align*}
	The expression  $\p_\beta^{u,w} \left(  \left| \{x\in V_{\mz}^n : x\sim V_u^n \} \right| > K \right)$ converges to $0$ for $K\to \infty$, by Markov's inequality, as one has the bound
	\begin{align*}
		\E_\beta & \left[  \left| \{x\in V_{\mz}^n , z \in V_u^n : x\sim z \} \right| \right] = \sum_{x\in V_{\mz}^n} \sum_{ z \in V_u^n} \p_\beta \left( x \sim z \right)
		\overset{\eqref{eq:connectionupper bound}}{\leq}
		\sum_{x\in V_{\mz}^n} \sum_{ z \in V_u^n} \frac{\beta}{\left(\|x-z\|_\infty -1 \right)^{2d}}\\
		&
		\leq \sum_{x\in V_{\mz}^n} \sum_{ z \in V_u^n} \frac{\beta}{\left((\|u\|_\infty -1)n \right)^{2d}}
		\leq \frac{\beta 2^{2d}}{\|u\|_\infty^{2d}} \leq \beta 2^{2d}.
	\end{align*}
	We need an upper bound on this quantity for the conditional measure $\p_\beta^{u,w}$.
	Lemma \ref{lem:number of edges conditioned} then gives that
	\begin{align*}
		\E_\beta^{u,w} \left[ \left| \{x\in V_{\mz}^n : x\sim V_u^n \} \right| \right] \leq  \E_\beta^{u,w} \left[  \left| \{x\in V_{\mz}^n , z \in V_u^n : x\sim z \} \right| \right] \leq \beta 2^{2d} +1
	\end{align*}
	and this upper bound does not depend on $n$ or $u$. Using Lemma \ref{lem:all quantiles of pointtobox}, we see that for fixed $\eps_1>0$ the term $\p_\beta \left( D \left(\mz,B_{\eps_1 n}(\mz)^C\right) \leq \eps \Lambda(n) \right)$ converges to $0$ as $\eps \to 0$ and thus we can take $K=K(\eps)$ and $\eps_1=\eps_1(\eps)$ that converge to $+\infty$, respectively $0$, slow enough such that $K \p_\beta \left( D  \left(\mz,B_{\eps_1 n}(\mz)^C\right) \leq \eps \Lambda(n) \right)$ also converges to $0$ for $\eps \to 0$.
\end{proof}

We want a similar function for the indirect distance between boxes. Such a function exists by Corollary \ref{coro:quantiles of boxes}.

\begin{definition}\label{lem:tightness_box to outside boxes}
	Let $g_2(\eps)$ be a function with $g_2(\eps) \underset{\eps \to 0}{\longrightarrow} 1$ such that the indirect distance $D^\star$ between the sets $V_\mz^n$ and $ B_n\left(V_\mz^n\right)^C $ satisfies
	\begin{align*}
	\p_\beta \left(  D^\star \left(V_\mz^n, B_n\left(V_\mz^n\right)^C \right) > \eps\Lambda(n,\beta) \right) \geq g_2(\eps)
	\end{align*}
	for all  $n \geq n(\eps)$ large enough.
\end{definition}
\noindent
Consider long-range percolation on $\Z^d$.
We split the long-range percolation graph into blocks of the form $V_v^n$, where $v\in \Z^d$. For each $v\in \Z^d$, we contract the block $V_v^n \subset \Z^d$ into one vertex $r(v)$. We call the graph that results from contracting all these blocks $G^\prime = (V^\prime, E^\prime)$.
For $r(v)\in G^\prime$, we define the neighborhood $\cN\left(r(v)\right)$ by
\begin{align*}
\cN\left(r(v)\right) = \left\{r(u) \in G^\prime : \|v-u\|_\infty \leq 1 \right\} \text ,
\end{align*}
and we define the neighborhood-degree of $r(v)$ by 
\begin{align*}
\deg^{\cN}(r(v)) = \sum_{r(u) \in \cN\left(r(v)\right)} \deg (r(u)) \text .
\end{align*}
We also define these quantities in the same way when we start with long-range percolation on the graph $V_\mz^{mn}$, and contract the box $V_v^n$ for all $v \in V_\mz^m$. The next lemma concerns the indirect distance between two sets, conditioned on the graph $G^\prime$.

\begin{lemma}\label{lem:goodeventD}
	Let $\mathcal{W}(\eps)$ be the event
	\begin{align*}
	\mathcal{W}(\eps) \coloneqq \left\{ D^\star \left( V_v^{n} , \bigcup_{u \in \Z^d : \| u - v\|_\infty \geq 2} V_u^{n} \right) > \eps \Lambda(n,\beta) \right\} \text .
	\end{align*}
	For all large enough $n\geq n(\eps)$ one has
	\begin{equation*}
	\p_{\beta} \left( \mathcal{W}(\eps)^C \ | \ G^\prime \right) \leq 3^d \deg^{\cN} \left(r(\mz)\right)  \left(1-g_1(\eps)\right) + \left(1-g_2(\eps)\right).
	\end{equation*}
\end{lemma}

\begin{proof}
	By translation invariance we can assume $v=\mz$. We define the set $T = V_\mz^n \cup \bigcup_{u \in \Z^d : \| u \|_\infty \geq 2} V_u^{n}$, and we define the events $\mathcal{W}_1(\eps)$ and $\mathcal{W}_2(\eps)$ by
	\begin{align*}
	\mathcal{W}_1(\eps) & = \Big\{ \exists a , b ,x,y\in \Z^d \text{ with } \|a\|_\infty=1 , \|a-b\|_\infty \geq 2, x \in V_a^n, y\in V_b^n : \\
	& \hspace{30mm} e=\{x,y\} \text{ open}, D(x,T;\omega^{e-}) \leq \eps \Lambda(n), D(y,T;\omega^{e-}) \leq \eps \Lambda(n) \Big\}
	\end{align*}
	and
	\begin{align*}
		\mathcal{W}_2(\eps) & = \Big\{ \text{There is an open path } P \text{ of length at most } \eps \Lambda(n) \text{ from } V_\mz^n \text{ to } \bigcup_{u \in \Z^d : \| u \|_\infty \geq 2} V_u^{n} : \\
		&
		\hspace{21mm} \forall \{x,y\}\in P 
		\text{ there exist } a,b \in \Z^d \text{ with } x \in V_a^n, y \in V_b^n , \|a-b\|_\infty \leq 1 \Big\} \text .
	\end{align*}
	We will now show that $\mathcal{W}(\eps)^C \subset \mathcal{W}_1(\eps) \cup \mathcal{W}_2(\eps)$. Assuming that $\mathcal{W}(\eps)^C$ holds, there exists an open path $P$ from $ V_\mz^n $ to $\bigcup_{u \in \Z^d : \| u \|_\infty \geq 2} V_u^n$ with length $\leq\eps \Lambda(n)$, and this path does not use a direct edge between these two sets. The path $P$ can either be of the form as described in the event $\mathcal{W}_2(\eps)$, or it contains an edge $e=\{x,y\}$ such that $x\in V_a^{n}, y \in V_b^n$ with $\|a\|_\infty =1, \|a-b\|_\infty \geq 2$. Let us assume that this path $P$ is not of the form as described in the event $\mathcal{W}_2(\eps)$. As the length of the path is at most $\eps \Lambda(n)$, the distance from the endpoints $x,y$ of such an edge to the set $T$ is at most $\eps \Lambda(n)$, even when the edge $\{x,y\}$ is removed. This holds, as the path $P$ starts at $V_\mz^n$, then uses the edge $e$, and then arrives in the set $\bigcup_{u \in \Z^d : \| u \|_\infty \geq 2} V_u^{n}$. Also note that $y\in T$ is possible, in which case the distance between $y$ and $T$ equals $0$.
	However, we see that $\mathcal{W}_1(\eps)$ holds. So we showed that $\mathcal{W}(\eps)^C \cap \mathcal{W}_2(\eps)^C$ implies $\mathcal{W}_1(\eps)$, and thus we also showed that $\mathcal{W}(\eps)^C \subset \mathcal{W}_1(\eps) \cup \mathcal{W}_2(\eps)$.
	
	The event $\mathcal{W}_2(\eps)$ implies that $D^\star \left(V_\mz^n, B_n(V_\mz^n)^C\right) \leq \eps \Lambda(n,\beta)$. Furthermore, the event $\mathcal{W}_2(\eps)$ is independent of $G^\prime$, which implies that
	\begin{equation*}
		\p_\beta \left(\mathcal{W}_2(\eps) | G^\prime\right) = \p_\beta \left(\mathcal{W}_2(\eps) \right) \leq \p_\beta \left( D^\star \left(V_\mz^n, B_n(V_\mz^n)^C\right) \leq \eps \Lambda(n,\beta) \right) \leq 1-g_2(\eps).
	\end{equation*}
	Suppose that $\|a\|_\infty = 1$ and $\|a-b\|_\infty \geq 2$, with $V_a^n \sim V_b^n$.
	Assume that there exists a path $P$ from $ V_\mz^n $ to $\bigcup_{u \in \Z^d : \| u \|_\infty \geq 2} V_u^n$ with length $\leq\eps \Lambda(n)$, that uses an edge $e=\{x,y\}$ with $x \in V_a^n, y\in V_b^n$. The path needs to get to $x$, and it enters the box $V_a^n$ from some box $V_e^n$ with $\|e\|_\infty \leq 1$. Say that the path enters the box $V_a^n$ through the vertex $z\in V_a^n$ with $z\sim V_e^n$. The chemical distance between $x$ and $z$ can be at most $\eps \Lambda(n,\beta)$. There are $3^d$ such vectors $e$, so the probability that there exists such a path is bounded by $3^d (1-g_1(\eps))$, as $\|a-b\|_\infty \geq 2$.
	With a union bound we get that
	\begin{align*}
		\p_\beta \left(\mathcal{W}_1(\eps) | G^\prime\right) & \leq \sum_{a:\|a\|_\infty = 1} \ \sum_{b : \|a-b\|_\infty\geq 2,V_a^n\sim V_b^n} \p_\beta \left( \exists x \in V_a^n, y\in V_b^n :  x\sim y, D(x,T;\omega^{e-}) < \eps \Lambda(n) \right)\\
		&
		\leq \sum_{a:\|a\|_\infty = 1} \deg(r(a)) 3^d (1-g_1(\eps))
		\leq
		\deg^{\cN}(r(\mz)) 3^d(1-g_1(\eps))
	\end{align*}
	and thus we finally get that
	\begin{align*}
		\p_\beta \left( \mathcal{W}(\eps)^C \right) \leq
		\p_\beta \left( \mathcal{W}_1(\eps) \right)
		+
		\p_\beta \left( \mathcal{W}_2(\eps) \right) 
		\leq 
		\deg^{\cN}(r(\mz)) 3^d(1-g_1(\eps)) + (1-g_2(\eps)) \text .
	\end{align*}
\end{proof}

\subsection{Supermultiplicativity of $\Lambda(n,\beta)$}

In this section, we prove the supermultiplicativity of $\Lambda(n,\beta)$. Our main tools for this are the results of the previous section and Lemma \ref{lem:connnectedsetsinLRP}. We also use the same notation as in Lemma \ref{lem:connnectedsetsinLRP}, i.e., $\mu_\beta=\E_\beta \left[\deg(\mz)\right]$ and $\overline{\deg}(Z) = \frac{1}{|Z|} \sum_{v \in Z}\deg(v)$.

\begin{lemma}\label{lem:supermultiplicativity}
	For all $\beta>0$, there exists a constant $c>0$ such that for all $m,n \in \N$
	\begin{align}\label{eq:supermultiplicativity}
		\Lambda(mn,\beta )\geq c \Lambda(n,\beta) \Lambda(m,\beta)\text .
	\end{align}
\end{lemma}

\begin{proof}
	Inequality \eqref{eq:supermultiplicativity} holds for all small $m $ or $n \in \N$ for some $c>0$, so it suffices to consider $m$ and $n$ large enough. We split the graph $V_{\mz}^{mn}$ into blocks of the form $V_v^n$, where $v\in V_{\mz}^m$. For each $v\in V_{\mz}^m$, we contract the block $V_v^n \subset V_{\mz}^{mn}$ into one vertex. We call the graph that results from contracting all these blocks $G^\prime = (V^\prime, E^\prime)$. The graph $G^\prime$ has the same distribution as long-range percolation on $V_{\mz}^m$ under the measure $\p_\beta$. By $r(v)$, we denote the vertex in $G^\prime$ that results from contracting the box $V_v^n$. We also define an analogy of the infinity-distance on $G^\prime$ by $\|r(u) - r(v)\|_\infty= \|u-v\|_\infty$. Our goal is to bound the expected distance between the vertices $\mz$ and $(mn-1)e_1$ from below, conditioned on the graph $G^\prime$. For this, we consider all loop-erased walks $P^\prime = \left(r(v_0), r(v_1),\ldots, r(v_k) \right)$ between $r(\mz)$ and $r((m-1)e_1)$ in $G^\prime$. In the following we always work on a certain event $\mathcal{H}_t$, which is defined by
	\begin{align*}
	\mathcal{H}_t = \bigcap_{k\geq t} \left\{	\left| \mathcal{CS}_k \left( G^\prime \right) \right| \leq 10^k \mu_\beta^k \right\} \cap 
	\bigcap_{k\geq t} \left\{ \overline{\deg}(Z) \leq 20 \mu_\beta, \ \forall Z \in \mathcal{CS}_k \left( G^\prime \right) \right\} \text .
	\end{align*}
	Note that, by Lemma \ref{lem:connnectedsetsinLRP}, \eqref{eq:connectedsetbound}, and Markov's inequality one has
	\begin{align}\label{eq:Htbound}
	\notag \p_\beta \left( \mathcal{H}_t^C \right) & \leq \sum_{k=t}^{\infty} \p_\beta \left( 	\left| \mathcal{CS}_k \left( G^\prime \right) \right| > 10^k \mu_\beta^k \right) + 
	\sum_{k=t}^{\infty} \p_\beta \left( \exists  Z \in \mathcal{CS}_k \left( G^\prime \right) : \overline{\deg}(Z) > 20 \mu_\beta  \right)\\
	& \leq \sum_{k=t}^{\infty} 0.4^k +  \sum_{k=t}^{\infty} e^{-4\mu_\beta k} \leq  \sum_{k=t}^{\infty} 0.5^k = 2\cdot 2^{-t}.
	\end{align}
	Let $P^\prime=\left(r(v_0),\ldots,r(v_k)\right)$ be a self-avoiding path in $G^\prime$ starting at the origin vertex, i.e., $v_0=\mz$. Assume that $k$ is large enough (which will be specified later) and let $\eps$ be small enough such that
	\begin{align}\label{eq:supermult eps }
		\left(	(27^d 50 \mu_\beta)^2 (1-g_1(\eps)) + 2 (1-g_2(\eps))\right)^{\frac{1}{30^d 200 \mu_\beta}}
		\leq
		\frac{1}{40 \mu_\beta} \text .
	\end{align}
	We will see later on, where this condition on $\eps$ comes from. We will now describe what it means for a block $V_{v_i}^n$ to be {\sl separated}; we will also say that the vertex $r(v_i) \in G^\prime$ is separated in this case. Intuitively, a block being separated ensures that a path in the original model that passes through this block needs to walk a distance of at least $\eps \Lambda(n,\beta)$. Formally, let $P$ be a path in the original graph $V_\mz^{mn}$ between $\mz$ and $(mn-1)e_1$, such that this path goes through the blocks corresponding to $r(u_0), r(u_1),\ldots,r(u_K)$ in this order. Let $P^\prime = \left(r(v_0),\ldots,r(v_k)\right)$ be the loop-erasure of the path $\left(r(u_0), r(u_1),\ldots,r(u_K)\right)$. So in particular, $P^\prime$ is self-avoiding.
	Suppose that $\|v_{i} - v_{i+1} \|_\infty \geq 2$. Then we call the block $r(v_i)$ separated if
	\begin{align*}
		D_{V_{v_{i}}^n}(x,y) \geq \eps \Lambda(n,\beta) \text{ for all } x,y \in V_{v_{i}}^n \text{ with } x \sim V_{v_{i+1}}^n, y \sim V_w^n, w \notin \{v_{i}, v_{i+1}\}.
	\end{align*}
	If $\|v_{i} - v_{i+1} \|_\infty = 1$, we call the block $r(v_i)$ separated if
	\begin{align*}
		D^\star_{V_\mz^{mn}} \left( V_{v_{i}}^n , \bigcup_{r(w) \in G^\prime : \|w-v_{i} \|_\infty \geq 2 } V_w^n\right) \geq \eps \Lambda(n,\beta)  \text .
	\end{align*}
	Next, we want to upper bound the probability that a block $r(v_i)$ is not separated, given the graph $G^\prime$. Assume that $\|v_{i}- v_{i+1} \|_\infty \geq 2$. Conditioned on the graph $G^\prime$, the probability that $r(v_i)$ is not separated is bounded by $\deg\left(r(v_i) \right)^2 \left(1-g_1(\eps)\right)$ for large enough $n$, by Lemma \ref{lem:graphspacing} and a union bound over all pairs of neighbors of $r(v_i)$. Assume that $\|v_{i}- v_{i+1} \|_\infty = 1$. Given the graph $G^\prime$, we have that
	\begin{align*}\label{eq:ungleichung}
	\p  & \left( D^{\star}_{V_\mz^{mn}}\Big(V_{v_{i}}^n , \bigcup_{r(v) : \|v-v_{i}\|_\infty \geq 2} V_v^n \Big) < \eps \Lambda(n,\beta) \ \big| \ G^\prime \right)\\
	&
	\leq 3^d
	\deg^{\cN}\left(r(v_{i})\right) (1-g_1(\eps)) + (1-g_2(\eps))
	\end{align*}
	for all large enough $n$, by \Cref{lem:goodeventD}. No matter whether $\|v_i-v_{i+1}\|_\infty = 1$ or $\|v_i-v_{i+1}\|_\infty > 1$, in both cases we have that
	\begin{align*}
		\p_\beta \left( r(v_i) \text{ not separated } | \ G^\prime \right) \leq 3^d
		\deg^{\cN}\left(r(v_{i})\right)^2 (1-g_1(\eps)) + (1-g_2(\eps)) \text .
	\end{align*}
	We define the set 
	\begin{equation*}
		\tilde{R}_k= \bigcup_{i=0}^{k-1}  \cN\left(r(v_i)\right) \text .
	\end{equation*}
	The set $\tilde{R}_k$ is a connected set in $G^\prime$, containing the origin $r(v_0)$, and its size is bounded from above and below by
	\begin{equation*}
		k \leq |\tilde{R}_k|\leq 3^d k  \text.
	\end{equation*}
	Assuming that the event $\cH_{k}$ holds, we get that the average degree of the set $\tilde{R}_k$ is bounded by $20 \mu_\beta k$. A vertex $r(v)$ can be included in several sets $\cN(r(v_i))$ for different $i$, but in at most $3^d$ many. So in particular we have
	\begin{align*}
		\sum_{i=0}^{k-1} \deg^{\cN} \left(r(v_i)\right) \leq 3^d |\tilde{R}_k| 20 \mu_\beta k \leq 9^d 20 \mu_\beta
	\end{align*}
	and thus there can be at most $\tfrac{k}{2}$ many indices $i \in \{0,\ldots, k-1\}$ with $\deg^{\cN} \left(r(v_i)\right) > 9^d 50 \mu_\beta$. We now define a set of special indices $\mathcal{IND}(P^\prime) \subset \{1,\ldots,k-1\}$ via the algorithm below. For abbreviation, we will mostly just write $\mathcal{IND}$ for $\mathcal{IND}(P^\prime)$, but one should remember that the indices really depend on the chosen path $P^\prime$.
	\begin{enumerate}\addtocounter{enumi}{-1}
		\item Start with $\mathcal{IND}_0 = \emptyset$.
		\item For $i=1,\ldots,k-1$: \\ If $\deg^{\cN} \left(r(v_i)\right) \leq 9^d 50 \mu_\beta$ and $\cN \left(r(v_i)\right) \nsim \bigcup_{j \in \mathcal{IND}_{i-1}} \cN \left( r(v_j) \right)$, then define $\mathcal{IND}_{i} = \mathcal{IND}_{i-1} \cup \{i\}$. Otherwise set $\mathcal{IND}_{i} = \mathcal{IND}_{i-1}$.
		\item Set $\mathcal{IND} \coloneqq \mathcal{IND}_{k-1}$.
	\end{enumerate}
	So in particular we have that for an index $i\in \mathcal{IND}$ it always holds that
	\begin{align*}
		\p_\beta & \left( r(v_i) \text{ not separated } | \ G^\prime \right) \leq
		3^d \deg^{\cN}\left(r(v_{i})\right)^2 (1-g_1(\eps)) + (1-g_2(\eps)) \\
		& \leq 
		(27^d 50 \mu_\beta)^2 (1-g_1(\eps)) + (1-g_2(\eps)) \eqqcolon g^\prime(\eps)
	\end{align*}
	On the event $\cH_{k}$, there are at least $\frac{k}{2}-1$ many indices $i \in \{1,\ldots,k-1\}$ with $\deg^{\cN} \left(r(v_i)\right) \leq 9^d 50 \mu_\beta$. Suppose that $V_v^n$ is a block with $V_v^n \sim \bigcup_{r(w) \in \cN \left(r(v_i)\right)} V_w^n$. (Note that all boxes $V_v^n$ with $r(v) \in  \cN \left(r(v_i)\right)$ are by definition adjacent to $\bigcup_{r(w) \in \cN \left(r(v_i)\right)} V_w^n$.) When we include the index $i$ to the set $\mathcal{IND}$, we can block all the indices $j>i$ with $r(v) \in \cN \left(r(v_j)\right)$. But for fixed $v$, there can be at most $3^d$ indices $j>i$ with $r(v) \in \cN \left(r(v_j)\right)$. So including one index $i$ with $\deg^{\cN} \left(r(v_i)\right) \leq 9^d 50 \mu_\beta$ to the set $\mathcal{IND}$, can block at most $3^d 9^d 50 \mu_\beta$ other indices. Thus we get that on the event $\cH_{k}$ one has for large enough $k$ that
	\begin{align*}
		|\mathcal{IND}| \geq \frac{\frac{k}{2}-1}{27^d 50 \mu_\beta+1} \geq \frac{k}{30^d 100 \mu_\beta} \text .
	\end{align*}
	Whether a block $V_{v_i}^n$ is separated in the path $P^\prime$ depends only on the edges with at least one endpoint in $\cN \left(r(v_i)\right)$. So in particular for different indices $i\in \mathcal{IND}$, it is independent whether the underlying blocks $V_{v_i}^n$ are separated. Thus we get that
	\begin{align*}
		&\p_\beta \left( \left|\left\{i \in \mathcal{IND}(P^\prime) : r(v_i) \text{ is separated }\right\}\right| \leq \frac{k}{30^d 200 \mu_\beta} \ \Big| \ G^\prime  \right)\\
		&
		\leq 2^{|\mathcal{IND}(P^\prime)|}
		 \left(g^\prime(\eps)\right)^{\frac{|\mathcal{IND}(P^\prime)|}{2}}
		\leq 2^{k} \left(g^\prime(\eps)\right)^{\frac{k}{30^d 200 \mu_\beta}}
		\leq \left(20 \mu_\beta\right)^{-k}
	\end{align*}
	where the last inequality holds because of our assumption on $\eps$ \eqref{eq:supermult eps }. With another union bound we get that
	\begin{align*}
		&\p_\beta \left( \exists P^\prime \text{ in } G^\prime \text{ of length } k \text{ s.t.} \left|\left\{i \in \mathcal{IND}(P^\prime) : r(v_i) \text{ is separated }\right\}\right| \leq \frac{k}{30^d 200 \mu_\beta} \ \Big| \ \mathcal{H}_k  \right)\\
		&
		\leq (10 \mu_\beta)^k \left(20 \mu_\beta\right)^{-k} = 2^{-k},
	\end{align*}
	where we say $P^\prime \text{ in } G^\prime$ if the path $P^\prime$ starts at $r(\mz)$ and is contained in the graph $G^\prime$.
	Using that $\p_\beta \left( \mathcal{H}_k^C \right) \leq 2 \cdot 2^{-k}$, we thus get that
	\begin{align*}
		\p_\beta \left( \exists P^\prime \text{ in } G^\prime \text{ of length } k \text{ s.t. } \left|\left\{i \in \mathcal{IND}(P^\prime) : r(v_i) \text{ is separated }\right\}\right| \leq \frac{k}{30^d 200 \mu_\beta} \right)
		\leq
		3 \cdot 2^{-k}.
	\end{align*}
	For abbreviation, we define the event $\cG_k$ by
	\begin{align*}
		\cG_k^C = \left\{\exists P^\prime \text{ in } G^\prime \text{ of length } k \text{ s.t. } \left|\left\{i \in \mathcal{IND}(P^\prime) : r(v_i) \text{ is separated }\right\}\right| \leq \frac{k}{30^d 200 \mu_\beta}\right\} \text .
	\end{align*}
	Assuming that the events $\cG_k$ and $D_{G^\prime}\left(r(\mz),r((m-1)e_1)\right) = k$ both hold, we get that for large enough $k$ one has $D_{V_\mz^{mn}}(\mz,(mn-1)e_1) \geq \frac{k\eps \Lambda(n,\beta)}{30^d 200 \mu_\beta}$.
	So in total we get that for some large enough $k^\prime$
	\begin{align}\label{eq:finalsubmulti}
	\notag &\E_\beta \left[ D_{V_{\mz}^{mn}}(\mz,(mn-1)e_1) \right] \geq 
	\sum_{k=k^\prime}^{\infty} 
	\E_\beta \left[ D_{V_{\mz}^{mn}}(\mz,(mn-1)e_1) \mathbbm{1}_{\mathcal{G}_k} \mathbbm{1}_{\left\{ D_{G^\prime}(r(\mz),r((m-1)e_1)) = k \right\}} \right]\\
	& \geq  \frac{\eps \Lambda(n,\beta)}{30^d 200 \mu_\beta} \sum_{k=k^\prime}^{\infty} k \E_\beta \left[ \mathbbm{1}_{\mathcal{G}_k} \mathbbm{1}_{\left\{ D_{G^\prime}(r(\mz),r((m-1)e_1)) = k \right\}} \right] \text ,
	\end{align}
	and we can further bound the last sum by
	\begin{align*}
	& \sum_{k=k^\prime}^{\infty} k \E_\beta \left[ \mathbbm{1}_{\mathcal{G}_k} \mathbbm{1}_{\left\{ D_{G^\prime}(r(\mz),r((m-1)e_1)) = k \right\}} \right] \\
	& =  
	\sum_{k=k^\prime}^{\infty} k \E_\beta \left[ \mathbbm{1}_{\left\{ D_{G^\prime}(r(\mz),r((m-1)e_1)) = k \right\}} \right]
	-
	\sum_{k=k^\prime}^{\infty} k \E_\beta \left[ \mathbbm{1}_{\left\{\mathcal{G}_k^C \right\}} \mathbbm{1}_{\left\{ D_{G^\prime}(r(\mz),r((m-1)e_1)) = k \right\}} \right]  \\
	& \geq  \sum_{k=k^\prime}^{\infty} k \E_\beta \left[ \mathbbm{1}_{\left\{ D_{V_{\mz}^m}(\mz,(m-1)e_1) = k \right\}} \right]
	-
	\sum_{k=k^\prime}^{\infty} k \E_\beta \left[ \mathbbm{1}_{\left\{\mathcal{G}_k^C \right\}}  \right]  \\
	& \geq \sum_{k=1}^{\infty} k \E_\beta \left[ \mathbbm{1}_{\left\{ D_{V_{\mz}^m}(\mz,(m-1)e_1) = k \right\}} \right] 
	-
	\sum_{k=1}^{k^\prime-1} k \E_\beta \left[ \mathbbm{1}_{\left\{ D_{V_{\mz}^m}(\mz,(m-1)e_1) = k \right\}} \right]
	-
	3\sum_{k=k^\prime}^{\infty} k 2^{-k} \\
	&
	\geq \E_\beta\left[ D_{V_{\mz}^m}(\mz,(m-1)e_1) \right]
	-
	k^\prime 
	-
	6 \geq c^\prime \Lambda(m, \beta)
	\end{align*}
	for small enough $c^\prime >0$ and $m$ large enough. 
	Inserting this into \eqref{eq:finalsubmulti} finishes the proof.
\end{proof}

\subsection{The diameter of boxes}\label{sec:dia}

In this section, we prove the second item of Theorem \ref{theo:exponent}, i.e., that the diameter of the box $\left\{0,\ldots, n-1\right\}^d$ and its expectation both grow like $n^\dxp$.

\begin{lemma}\label{lem:diameter}
	For all $\beta\geq 0$ one has
	\begin{align*}
	n^{\dxp(\beta)} \approx_P \dia\left( \left\{0,\ldots, n-1\right\}^d \right)  \approx_P \E_\beta\left[ \dia\left( \left\{0,\ldots, n-1 \right\}^d \right) \right] \text .
	\end{align*}
\end{lemma}
\begin{proof}
	By Proposition \ref{propo:scaling}, it suffices to consider the case when $n=2^k$ for some $k\in \N$. We have 
	\begin{align*}
		 \dia\left( \left\{0,\ldots, n-1\right\}^d \right) \geq D_{V_\mz^n}\left(\mz, (n-1)e_1\right)
	\end{align*}
	and this already implies that for each $\eps >0$ there exist constants $c,c_\eps>0$ such that
	\begin{align*}
		\p_\beta\left( c_\eps n^{\dxp(\beta)} < \dia\left( \left\{0,\ldots, n-1\right\}^d \right)  \right) > 1-\eps
	\end{align*}
	and
	\begin{align*}
		c n^{\dxp(\beta)} \leq \E_\beta\left[ \dia\left( \left\{0,\ldots, n-1 \right\}^d \right) \right] 
	\end{align*}
	uniformly over $n$. For the upper bound, we make a dyadic decomposition of the box $V_\mz^n$. Similar ideas were also used in \cite{ding2013distances} for one dimension or in \cite{baeumler2022behaviour}. For a constant $S\geq 1$, we say that a box $V_y^{2^l} \subset V_\mz^{2^k}$ is {\sl $S$-good} if
	\begin{align*}
		D_{V_y^{2^l}} \left(2^l y,2^ly+(2^l-1)e\right) \leq S \left(\frac{3}{2}\right)^{(l-k)\dxp} 2^{k\dxp} 
	\end{align*}
	for all $e \in \{0,1\}^d$, where we simply write $\dxp$ for $\dxp(\beta)$ from here on.  We use the notation
	\begin{align*}
	\Omega_l^S = \bigcap_{y \in V_{\mz}^{2^{k-l}}} \left\{ V_y^{2^l} \text{ is $S$-good} \right\} \ \ \text{ and } \ \ 
	\Omega^S = \bigcap_{l=1}^k \Omega_l^S .
	\end{align*}
	On the event $\Omega^S$, we can bound the graph distance between $\mz$ and any $y \in V_\mz^{2^k}$ by considering a path that goes along the boxes in a dyadic decomposition. Let $y_0,\ldots, y_k \in \Z^d$ be such that $y \in V_{y_i}^{2^i}$ for all $i$. So in particular $y_0=y$ and $y_k=\mz$. We also have that $V_{y_0}^{2^0} \subset V_{y_1}^{2^1} \subset \ldots \subset V_{y_k}^{2^k}$ and thus also $2^{i-1} y_{i-1} \in V_{y_i}^{2^i}$ for all $i\geq 1$. This implies that $2^{i-1} y_{i-1} = 2^iy_i + 2^{i-1}e$ for some $e \in \{0,1\}^d$. As all the boxes inside $V_\mz^{2^k}$ were assumed to be $S$-good we have
	\begin{align*}
		& D_{V_{\mz}^{2^k}} \left( 2^{i} y_{i} , 2^{i-1} y_{i-1} \right) \leq D_{V_{\mz}^{2^k}} \left( 2^{i} y_{i} , 2^i y_i + (2^{i-1}-1)e \right) + 1\\
		& \leq 
		D_{V_{2 y_i}^{2^{i-1}}} \left( 2^{i-1} 2 y_{i} , 2^{i-1} 2 y_i + (2^{i-1}-1)e \right) + 1 \leq S \left(\frac{3}{2}\right)^{(i-1-k)\dxp} 2^{k\dxp}+1\text.
	\end{align*}
	Now we have by the triangle inequality
	\begin{align*}
		D_{V_0^{2^k}} \left(\mz, v \right) & \leq \sum_{l=1}^{k} \left( S \left(\frac{3}{2}\right)^{(l-1-k) \dxp}2^{k\dxp} + 1 \right) \leq  S 2^{k\dxp}  \sum_{l=1}^{k} \left(  \left(\frac{3}{2}\right)^{(l-1-k) \dxp} + \frac{1}{2^{k\dxp}} \right)\\
		& \leq C_\dxp S 2^{k\dxp}
	\end{align*}
	where the constant $C_\dxp$ depends only on $\dxp$. As $D(u,v) \leq D(\mz,u) + D(\mz,v)$ for all $u,v \in V_\mz^{2^k}$, the previous bound already implies that on the event $\Omega^S$ one has
	\begin{align}\label{eq:diabound on Omega}
		\dia \left( V_\mz^{2^k} \right) \leq 2 C_\dxp S 2^{k\dxp}
	\end{align}
	and thus it suffices to bound the probability of $\left(\Omega^S\right)^C$. We know from Corollary \ref{coro:allmoments} that the $r$-th moment of $D_{V_y^{2^l}} \left(2^l y,2^ly+(2^l-1)e\right)$ is of order $2^{rl\dxp}$, for all $r \geq 0$.
	So by a union bound and Markov's inequality we get that for every fixed box $V_y^{2^l}$
	\begin{align*}
		&\p_\beta\left( V_y^{2^l} \text{ is not $S$-good}\right) \leq \sum_{e \in \{0,1\}^d} \p_\beta \left( D_{V_y^{2^l}} \left(2^l y,2^ly+(2^l-1)e\right) > S \left(\frac{3}{2}\right)^{(l-k)\dxp} 2^{k\dxp} \right)
		\\
		& = \sum_{e \in \{0,1\}^d} \p_\beta \left( D_{V_\mz^{2^l}} \left(\mz,\mz+(2^l-1)e\right)^\frac{4d}{\dxp} > S^\frac{4d}{\dxp}  \left(\frac{3}{2}\right)^{(l-k)4d} 2^{4dk} \right)\\
		& \leq \sum_{e \in \{0,1\}^d} \frac{\E_\beta \left[ D_{V_\mz^{2^l}} \left(\mz,\mz+(2^l-1)e\right)^{\frac{4d}{\dxp}} \right]}{S^\frac{4d}{\dxp} \left(\frac{3}{2}\right)^{(l-k)4d} 2^{4dk}} 
		\leq \frac{C\cdot 2^{l\dxp\frac{4d}{\dxp}}}{ S^\frac{4d}{\dxp} \left(\frac{3}{2}\right)^{(l-k)4d} 2^{4dk} }
		=
		\frac{C\cdot 2^{4dl}}{ S^\frac{4d}{\dxp} \left(\frac{3}{2}\right)^{(l-k)4d} 2^{4dk} } \\
		& \leq \frac{C}{S^{\frac{4d}{\dxp}}} \left(\frac{2}{3}\right)^{(l-k)4d} 2^{(l-k)4d}
		=  \frac{C}{S^{\frac{4d}{\dxp}}} \left(\frac{4}{3}\right)^{(l-k)4d}
	\end{align*}
	for some constant $C<\infty$ depending only on $d$ and $\beta$.
	With another union bound we get that
	\begin{align*}
		\p_\beta \left(\left(\Omega_l^S\right)^C\right) & \leq 
		\sum_{y \in V_{\mz}^{2^{k-l}}} \p_\beta \left( V_y^{2^l} \text{ is not $S$-good} \right)
		\leq \sum_{y \in V_{\mz}^{2^{k-l}}}  \frac{C}{S^{\frac{4d}{\dxp}}} \left(\frac{4}{3}\right)^{(l-k)4d} \\
		& = \frac{C}{S^{\frac{4d}{\dxp}}} 2^{(k-l)d}   \left(\frac{4}{3}\right)^{(l-k)4d} = \frac{C}{S^{\frac{4d}{\dxp}}} \left(\frac{81}{128}\right)^{(k-l)d}
	\end{align*}
	which implies that
	\begin{equation}\label{eq:Omega bound}
		\p_\beta\left( \left( \Omega^S \right)^C \right) \leq \sum_{l=1}^{k} \frac{C}{S^{\frac{4d}{\dxp}}} \left(\frac{81}{128}\right)^{(k-l)d}
		\leq \frac{C^\prime}{S^{\frac{4d}{\dxp}}} 
	\end{equation}
	for some constant $C^\prime< \infty$. Together with \eqref{eq:diabound on Omega}, this proves that $\dia \left( V_\mz^{2^k} \right) \approx_P 2^{k\dxp}$. Inequality \eqref{eq:diabound on Omega} also implies that 
	\begin{align*}
		 \left\{\frac{\dia\left(V_\mz^{2^k}\right)}{2^{k\dxp}} > S\right\} \subset \left(\Omega^{\frac{S}{2C_\dxp}}\right)^C
	\end{align*}
	whenever $\frac{S}{2C_\dxp} > 1$, and this implies that for some finite $K \in \N$ and all $k \in \N$
	\begin{align*}
		\E_\beta \left[ \frac{\dia\left(V_\mz^{2^k}\right)}{2^{k\dxp}} \right] 
		& \leq K + \sum_{S=K}^{\infty} \p_\beta \left( \frac{\dia\left(V_\mz^{2^k}\right)}{2^{k\dxp}} > S \right)
		\leq K + \sum_{S=K}^{\infty} \p_\beta \left(  \left(\Omega^{\frac{S}{2C_\dxp}}\right)^C \right)\\
		&
		\leq K + \sum_{S=1}^{\infty} \frac{C^\prime \left(2C_\dxp\right)^\frac{4d}{\dxp} }{S^\frac{4d}{\dxp}} < \infty
	\end{align*}
	where the last term is finite as $\frac{4d}{\dxp} > 1$. This also shows that 
	\begin{align*}
		\E_\beta \left[ \dia\left(V_\mz^{2^k}\right) \right] = \mathcal{O} \left( 2^{k\dxp} \right)
	\end{align*}
	and thus we finish the proof of  Lemma \ref{lem:diameter}.
\end{proof}

\section{Tail behavior of the distances and diameter}\label{section:tail behavior}

\Cref{theo:exponent} shows that the random variables $\frac{D(\mz,u)}{\|u\|^{\dxp(\beta)}}$ are tight in $(0,\infty)$ under the measure $\p_\beta$. In this section, we give more precise estimates on the tail-behavior of the random variables $\frac{D(\mz,u)}{\|u\|^{\dxp(\beta)}}$. We describe this tail behavior via functions $f$ for which $\sup_{u\in \Z^d \setminus\{\mz\}}\E_\beta \left[ f\left(\frac{D(\mz,u)}{\|u\|^{\dxp(\beta)}}\right) \right]$ is finite or infinite. This result is also a useful tool in \cref{section:comparison}, and in the companion paper \cite{baeumler2022behaviour}.

\begin{theorem}\label{theo:tail behavior}
	For all $\eta<\frac{1}{1-\dxp(\beta)}$ one has
	\begin{align}\label{eq:stretched upper bound}
	\sup_{n \in \N } \E_\beta \left[ \exp \left(\left(\frac{\dia \left(\{0,\ldots,n\}^d\right)}{n^{\dxp(\beta)}}\right)^\eta\right) \right] < \infty \text .
	\end{align}	
\end{theorem}

For dimension $d=1$, the bound given by \eqref{eq:stretched upper bound} is sharp, as the following lemma shows.

\begin{lemma}
For all dimensions $d$ and all $\beta>0$, there exists a constant $t>0$ such that
	\begin{align}\label{eq:stretched lower bound}
	\sup_{u \in \Z^d\setminus \{\mz\}} \E_\beta \left[ \exp \left( t \left(\frac{D(\mz,u)}{\|u\|^{\dxp(\beta)}}\right)^\frac{d}{1-\dxp(\beta)}\right) \right] = \infty \text .
	\end{align}
\end{lemma}
\begin{proof}
	We define the event
	\begin{align*}
	\mathcal{D}_n = \bigcap_{v \in B_n(\mz)} \left\{v \nsim w \text{ for all $w\in \Z^d$ with }\|v-w\|_\infty \geq 2\right\}
	\end{align*}
	If the event $D_{\|u\|_\infty}$ occurs, the shortest path between $\mz$ and $u$ uses nearest-neighbor edges only and thus has a length of $\|u\|_\infty$.
	Using the FKG-inequality, we get that
	\begin{align*}
	\p_\beta \left(\mathcal{D}_n\right) \geq \p_\beta( \mz\nsim w \text{ for all $w\in \Z^d$ with }\|w\|_\infty \geq 2 )^{|B_n(\mz)|} \geq e^{-C n^d}
	\end{align*}
	for some constant $C < \infty$. Thus we see that 
	\begin{align*}
	\p_\beta \left( \frac{D(\mz,u)}{\|u\|_\infty^{\dxp(\beta)}} = \|u\|_\infty^{1-\dxp(\beta)} \right) \geq \p_\beta \left(\mathcal{D}_{n}\right) \geq \exp\left(-C \|u\|_\infty^d\right)
	\end{align*}
	and from here one can easily verify that \eqref{eq:stretched lower bound} holds for $t$ large enough.
\end{proof}

\begin{remark}
Conditioning on the event that there is no edge longer than $\eps m^{-\tfrac{1}{1-\dxp(\beta)}} n$ open in the box $B_n(\mz)$ for $\eps>0$ small enough, one can actually show that for all $u\in \Z^d$ with $\|u\|_\infty=n$ one has
	\begin{align*}
		\p_\beta \left( \frac{D(\mz,u)}{\|u\|^{\dxp(\beta)}} > m \right) \geq \exp\left(- C m^{-\frac{d}{1-\dxp(\beta)}}\right)
	\end{align*}
	for some constant $C \in \R_{>0}$, and all large enough $n$.
\end{remark}
For a sequence of positive random variables $(X_n)_{n\in\N}$ and some $\eta>0$, we have that
\begin{align*}
	&\E \left[ \exp \left(X_n^\eta\right) \right] = \int_{0}^{\infty} \p \left(  \exp \left(X_n^\eta\right) > s \right) \md s
	=
	\int_{0}^{\infty} \p \left(    X_n^\eta > \log(s) \right) \md s
	\\
	&
	= 1 + 
	\int_{1}^{\infty} \p \left(    X_n > \log(s)^{1/\eta} \right) \md s
	=
	1 + 
	\int_{0}^{\infty} \p \left(    X_n > s \right) \eta s^{\eta-1} \exp \left(s^\eta\right) \md s \text .
\end{align*}
So in particular, if there exist constants $0<c,C<\infty$ such that 
\begin{align}\label{eq:stretched prob upper bound}
	\p \left(X_n > s\right) \leq C \exp\left(-c s^{\bar{\eta}}\right) 
\end{align}
for all $n\in \N$, this implies that $\sup_{n \in \N} \E \left[ \exp \left(X_n^\eta\right) \right] < \infty$ for all $\eta \in (0,\bar{\eta})$. So in fact we will often show \eqref{eq:stretched prob upper bound} in the following, as this will already imply statements of the form $\sup_{n \in \N} \E \left[ \exp \left(X_n^\eta\right) \right] < \infty$, as in \eqref{eq:stretched upper bound}.
\Cref{theo:tail behavior} directly implies that that for all $\eta<\frac{1}{1-\dxp}$ one has
\begin{align}\label{eq:stretched upper bound2}
\sup_{n \in \N} \E_\beta \left[ \exp \left(\left(\frac{D_{V_\mz^n}(\mz,(n-1)e_1)}{n^\dxp}\right)^\eta\right) \right] < \infty \text ,
\end{align}
whereas \eqref{eq:stretched upper bound2} does not directly imply any statements about the diameter of boxes as in \eqref{eq:stretched upper bound}. However, a slightly weaker statement can be deduced from a slight modification of \eqref{eq:stretched upper bound2}, as the next lemma shows.

\begin{lemma}\label{lem:distance to dia}
	Suppose that
	\begin{align}\label{eq:distance some eta}
	\sup_{n \in \N} \E_\beta \left[ \exp \left(\left(\frac{D_{V_\mz^n}(\mz,(n-1)e)}{n^\dxp}\right)^\eta\right) \right] < \infty
	\end{align}
	for some $\eta>0$ and all $e\in \{0,1\}^d$. Then there exist constants $C,C_\dxp \in \R_{> 0}$ such that 
	\begin{align*}
	\p_\beta \left( \dia \left(V_\mz^{\bar{n}}\right) > S  C_\dxp n^\dxp \text{ for some } \bar{n} \in \{0,\ldots,n\} \right) 
	\leq C \exp\left(- S^{\eta}\right) \text ,
	\end{align*}
	which implies that
	\begin{align}\label{eq:diameter lower eta}
	\sup_{n \in \N} \E_\beta \left[ \exp \left( \left( \frac{\dia \left(\{0,\ldots,n-1\}^d\right)}{n^\dxp}\right)^{\bar{\eta}}\right) \right] < \infty
	\end{align}
	for all $\bar{\eta} \in (0,\eta)$.
\end{lemma}

\begin{proof}
	We do the proof for $n=2^k$ with $k\in \N$. The proof for general $n \in \N$ follows by \Cref{propo:scaling}. For $S\geq 1$ and $l \in \{0,\ldots,k\}$, define the events
	\begin{align*}
	\Omega^S_l = \bigcap_{y \in V_{\mz}^{2^{k-l}}} \bigcap_{e \in \{0,1\}^d} \left\{ D_{V_y^{2^l}}\left(2^l y , 2^l y + (2^l-1)e \right) \leq S \left(\frac{3}{2}\right)^{(l-k)\dxp} 2^{k\dxp}  \right\}
	\end{align*}
	and
	\begin{align*}
	\Omega^S = \bigcap_{l=0}^{k} \Omega_l^S.
	\end{align*}
	On the event $\Omega^S$, for all $\bar{n}\leq n$, and for any $y\in V_\mz^{\bar{n}}$, we can bound the graph distance between $\mz$ and $y$ by considering a dyadic path between them, and thus we get that on the event $\Omega^S$
	\begin{align*}
		D_{V_\mz^{\bar{n}}} \left(\mz, y \right) \leq  \sum_{l=0}^{k} S \left(\frac{3}{2}\right)^{(l-k)\dxp} 2^{k\dxp}  + k \text,
	\end{align*}
	and this already implies that
	\begin{align*}
	\dia \left(V_\mz^{\bar{n}}\right) \leq 2 \left( \sum_{l=0}^{k} S \left(\frac{3}{2}\right)^{(l-k)\dxp} 2^{k\dxp}  + k\right)
	\leq 
	C_\dxp S 2^{k\dxp}
	\end{align*}
	for some constant $C_\dxp < \infty$ and all $\bar{n}\leq n$. So in particular we see that the event $\left\{\dia \left(V_\mz^{\bar{n}}\right) > S  C_\dxp n^\dxp \text{ for some } \bar{n}\leq n \right\}$ implies that $\Omega_l^S$ does not hold for some $l \in \{0,\ldots, k\}$. So with a union bound we get that
	\begin{align}\label{eq:stretched exp 1}
	& \notag \p_\beta \left( \dia \left(V_\mz^{\bar{n}}\right) > S  C_\dxp n^\dxp \text{ for some } \bar{n}\leq n \right) \\
	&
	\leq
	\sum_{l=0}^{k} 2^{(k-l)d} \sum_{e\in \{0,1\}^d} \p_\beta \left( D_{V_\mz^{2^l}}\left(\mz , (2^l-1)e \right) > S \left(\frac{3}{2}\right)^{(l-k)\dxp} 2^{k\dxp}  \right).
	\end{align}
	By Markov's inequality and \eqref{eq:distance some eta} we have for any $e\in \{0,1\}^d$
	\begin{align*}
	& \p_\beta \left( D_{V_\mz^{2^l}}\left(\mz , (2^l-1)e \right) > S \left(\frac{3}{2}\right)^{(l-k)\dxp} 2^{k\dxp}  \right)
	=
	\p_\beta \left( D_{V_\mz^{2^l}}\left(\mz , (2^l-1)e \right) > S \left(\frac{4}{3}\right)^{(k-l)\dxp} 2^{l\dxp}  \right)\\
	&
	=
	\p_\beta \left( \exp \left( \left(\frac{D_{V_\mz^{2^l}}\left(\mz , (2^l-1)e \right)}{2^{l\dxp}} \right)^\eta \right) > \exp \left(  S^\eta \left(\frac{4}{3}\right)^{(k-l)\dxp \eta}   \right) \right)\\
	&
	\leq
	\E_\beta \left[ \exp \left( \left(\frac{D_{V_\mz^{2^l}}\left(\mz , (2^l-1)e \right)}{2^{l\dxp}} \right)^\eta \right)  \right] 
	\exp \left( - S^\eta \left(\frac{4}{3}\right)^{(k-l)\dxp \eta}   \right) 
	\leq C_\eta \exp \left( - S^\eta \left(\frac{4}{3}\right)^{(k-l)\dxp \eta}   \right) 
	\end{align*}
	for some constant $C_\eta < \infty$. Inserting this into \eqref{eq:stretched exp 1} shows that
	\begin{align*}
	\p_\beta \left( \dia \left(V_\mz^{\bar{n}}\right) > S  C_\dxp n^\dxp \text{ for some } \bar{n}\leq n \right) & \leq 2^d
	\sum_{l=0}^{k} 2^{(k-l)d} C_\eta \exp \left( - S^\eta \left(\frac{4}{3}\right)^{(k-l)\dxp \eta}   \right) \\
	&
	\leq C \exp\left(- S^{\eta}\right)
	\end{align*}
	for some constant $C< \infty$. By taking the constant $C$ large enough we can also guarantee that the above inequality holds for all $S>0$. This already implies that \eqref{eq:diameter lower eta} holds for all $\bar{\eta} \in (0,\eta)$.
\end{proof}

\begin{lemma}\label{lem:moments for e1}
	For all $\beta \geq 0$ and all $e \in \{0,1\}^d$ one has
	\begin{align}\label{eq:expomoments for 1}
	\sup_{n \in \N}
	\E_\beta \left[\exp \left( \left(\frac{D_{V_\mz^n}(\mz,(n-1)e)}{n^\dxp}\right)^{0.5} \right) \right] < \infty \text .
	\end{align}
\end{lemma}

\begin{proof}
	First, we will consider $e=e_1$ only.
	We define a process $(a_k(n))_{k\in \N}$.
	We start with $a_0(n)=0$ and define $a_k(n)$ inductively by
	\begin{align}\label{eq:ak definition}
	\notag & a_{k+1}(n) =  (a_k(n) + 2) \\
	& + \sup \left\{z \in \N_{>0} : D_{(a_k(n)+2)e_1+\{0,\ldots,z\}^d}\big((a_k(n)+2)e_1, (a_k(n)+2)e_1 + ze_1\big) \leq n^\dxp \right\}.
	\end{align}
	
	\begin{figure}
		\begin{tikzpicture}[scale=0.35]
			\fill[color=lightergray] (2,0) rectangle (6,4);
			\fill[color=lightergray] (8,0) rectangle (11,3);
			\fill[color=lightergray] (13,0) rectangle (17,4);
			\draw[thick] (0,0) rectangle (14,14);
			\foreach \x in {-1,0,1,2,3,4,5,6,7,8,9,10,11,12,13,14,15,16,17} {
				\foreach \y in {-1,0,1,2,3,4,5,6,7,8,9,10,11,12,13,14,15} {
					\vertex[fill,color=gray,minimum size=2 pt ]  at (\x,\y) {};
				};
			};
			\vertex[draw=none] (A) at (-1.2,0) {};
			\vertex[draw=none] (B) at (18,0) {};
			\path[very thick , ->] (A) edge (B);
			\vertex[draw=none] (C) at (0,-1.2) {};
			\vertex[draw=none] (D) at (0, 16) {};
			\path[very thick , ->] (C) edge (D);
			\vertex[fill, color=black, minimum size = 4pt, label=$a_0 \  \  \ $] (a0) at (0, 0) {};
			\vertex[fill, color=black, minimum size = 4pt, label=$a_1 \  \  \ $] (a1) at (6, 0) {};
			\vertex[fill, color=black, minimum size = 4pt, label=$a_2 \  \  \ $] (a2) at (11, 0) {};
			\vertex[fill, color=black, minimum size = 4pt, label=$a_3 \  \  \ $] (a3) at (17, 0) {};
		\end{tikzpicture}
		\centering
		\parbox{11cm}{\caption{An example of the process $a_0e_1,a_1e_1,a_2e_1,a_3e_1$ for the Graph $V_\mz^{15}$. We have $K=3$, as $a_2(15) \leq 14 < a_3(15)$.} \label{Fig:K}}
	\end{figure}
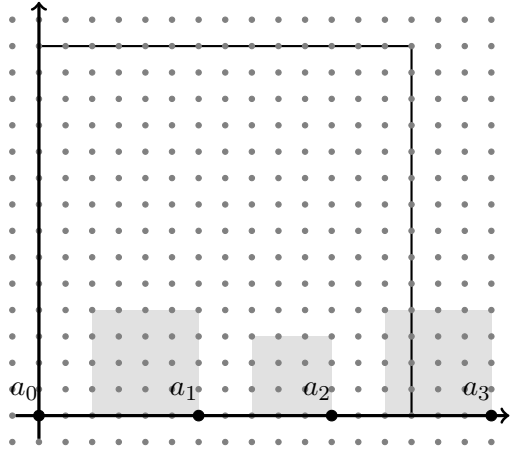

	The idea behind this definitions is that the graph distance between $a_k(n)e_1$ and $a_{k+1}(n)e_1$ is always upper bounded by $n^\dxp +2$; so in order to bound the graph distance between $\mz=a_0(n)e_1$ and $(n-1)e_1$, one can consider the path that goes from $a_0(n)e_1$ to $a_1(n)e_1$ to $a_2(n)e_1$, and inductively to $a_K(n)e_1$ for properly chosen $K$ such that $a_K(n)e_1$ is close to $(n-1)e_1$.
	Given the long-range percolation graph, the sequence $\left(a_k(n)\right)_{k\in \N}$ can be constructed as follows: for given $a_{k-1}(n)$, we walk along the $e_1$-axis in positive direction, starting at $(a_{k-1}(n)+2)e_1$. We do this until the graph distance between $(a_{k-1}(n)+2)e_1$ and $(a_{k-1}(n)+2+z)e_1$ exceeds a certain threshold $(n^\dxp)$, and then we go one step back, i.e.,  in negative $e_1$-direction, and then define this point as $a_{k}e_1$. See Figure \ref{Fig:K} for an illustration. This procedure only reveals information about edges with both endpoints in the slice $\{y\in \Z^d : a_{k-1}(n)+2 \leq \langle y, e_1 \rangle \leq a_{k}(n)+1\}$, so in particular the differences $(a_{k^\prime+1}(n)-a_{k^\prime}(n))$ are independent of $a_k(n)$ for $k^\prime \geq k$. By translation invariance, the differences $\left(a_{k+1}(n)-a_k(n)\right)_{k\in \N_0}$ are independent and identically distributed random variables.
	The graph distance between $a_k(n)e_1$ and $a_{k+1}(n)e_1$ is always bounded by $n^\dxp+2$, as we can go from $a_k(n) e_1$ to $(a_k(n)+2)e_1$ in two steps and from there to $a_{k+1}(n)e_1$ in at most $n^\dxp$ steps. Define
	\begin{equation*}
	K_n= \inf \left\{k \in \N : a_k(n) \geq n \right\} \ 
	\end{equation*}
	as the index of the first point $a_k(n)e_1$ that lies outside of $V_\mz^n$.
	Then one has 
	\begin{align*}
	D_{V_\mz^n}\left(\mz,(n-1)e_1\right) \leq K_n n^\dxp + 2K_n \leq 3K_n n^\dxp
	\end{align*}
	as one can walk through the path that goes from $\mz$ to $a_1(n)e_1$, from $a_1(n)e_1$ to $a_2(n) e_1$, and from there in the same manner inductively to $a_{K_n-1}(n)e_1$, and from there to $(n-1)e_1$. So our next goal is to show that $K_n$ is typically not too large. We use that for all $\beta \geq 0$ there exists an $\alpha > 0$ such that
	\begin{equation}\label{eq:incrementbound}
	\p_\beta \left(\frac{a_{k+1}(n)- a_{k}(n)}{n} \geq \alpha \right) \geq 0.5, 
	\end{equation}
	which we will prove in Lemma \ref{lem:incrementbound} below. We define the indices $k_0(n),k_1(n),\ldots$ by $k_0(n)=0$ and
	\begin{align*}
	k_{i+1}(n)= \inf\{k > k_i(n): \frac{a_{k}(n)-a_{k-1}(n)}{n} \geq \alpha \} .
	\end{align*}
	By construction we have $K_n \leq k_{\lceil 1/\alpha \rceil +1}(n)$. So in particular we have
	\begin{align*}
	\frac{D_{V_\mz^n}\left(\mz, (n-1)e_1\right)}{n^\dxp} \leq 3 K_n \leq 3  k_{\lceil 1/\alpha \rceil +1}(n) = 3 \sum_{i=0}^{\lceil \frac{1}{\alpha} \rceil} \left(k_{i+1}(n)-k_i(n)\right) \text.
	\end{align*}
	The differences $\left(k_{i+1}(n)-k_i(n)\right)_{i\geq 0}$ are independent random variables and are, by \eqref{eq:incrementbound}, dominated by $\text{Geometric}\left(\frac{1}{2}\right)$-distributed random variables. This already implies that 
	\begin{align}\label{eq:t exponential}
	&\notag \E_{\beta } \left[ \exp \left( t \frac{D_{V_\mz^n}\left(\mz, (n-1)e_1\right)}{n^\dxp} \right) \right]
	\leq
	\E_{\beta } \left[ \exp \left(  t3 \sum_{i=0}^{\lceil \frac{1}{\alpha} \rceil} k_{i+1}(n)-k_i(n) \right) \right]\\
	&
	=
	\prod_{i=0}^{\lceil \frac{1}{\alpha} \rceil}
	\E_{\beta } \left[ \exp \left(  t3  (k_{i+1}(n)-k_i(n))  \right) \right] \leq C < \infty
	\end{align}
	for some $t>0$ small enough and a uniform constant $C$ that does not depend on $n$, as the differences $k_{i+1}(n)-k_i(n)$ are dominated by a $\text{Geometric}\left(\frac{1}{2}\right)$-distributed random variable. This shows that \eqref{eq:expomoments for 1} holds for $e=e_1$. 
	To extend this proof to general $e\in \{0,1\}^d$, we use the same technique as in the proof of \Cref{lem:endisfar}. For $i\in \{0,\ldots, d\}$, we define $e(i)$ by
	\begin{equation*}
		e(i)=\sum_{j=1}^{i} p_j(e) e_i \text ,
	\end{equation*}
	and thus we get by the triangle inequality that
	\begin{equation*}
		D_{V_\mz^n}(\mz, (n-1)e) \leq \sum_{i=1}^d D_{V_\mz^n}((n-1)e(i-1), (n-1)e(i)) \text .
	\end{equation*}
	The random variables $D_{V_\mz^n}((n-1)e(i-1), (n-1)e(i))$ are either equal to $0$, when $e(i-1)$ and $e(i)$ coincide, or they have the same distribution as $D_{V_\mz^n}(\mz, (n-1)e_1)$, when $e(i-1)$ and $e(i)$ lie on adjacent corners of the cube $V_\mz^n$.
	H{\"o}lder's inequality implies that
	\begin{align*}
		& \E_\beta \left[ \exp \left( D_{V_\mz^n}(\mz, (n-1)e)^{0.5} \right) \right]
		\leq
		\E_\beta \left[ \exp \left( \sum_{i=1}^d D_{V_\mz^n}((n-1)e(i-1), (n-1)e(i))^{0.5}  \right) \right]\\
		&
		\leq
		\prod_{i=1}^d
		\E_\beta \left[ \exp \left(  d D_{V_\mz^n}((n-1)e(i-1), (n-1)e(i))^{0.5}  \right) \right]^{\frac{1}{d}}
		\leq
		\E_\beta \left[ \exp \left(  d D_{V_\mz^n}(\mz,  (n-1)e_1)^{0.5}  \right) \right]
	\end{align*}
	and the last term is finite uniformly over all $n\in \N$, which follows from \eqref{eq:t exponential}.
\end{proof}

\begin{lemma}\label{lem:incrementbound}
	For all $\beta >0$, there exists a constant $\alpha>0$ such that for all $n \in \N_{>0}$
	\begin{equation}\label{eq:incrementbound2}
	\p_\beta \left(\frac{a_{k+1}(n) - a_{k}(n) }{n} \geq \alpha\right) \geq 0.5 \text .
	\end{equation}
\end{lemma}
\begin{proof}
	As the differences $(a_{k+1}(n) - a_{k}(n))_{k\geq 0}$ are identically distributed, it suffices to consider the case $k=0$. The proof uses a dyadic decomposition along the $e_1$-axis. Let $n$ be large enough so that $\log_2(n)\leq \frac{n^\dxp}{2}$; this holds for all $n$ sufficiently large. We can make this assumption, as the statement \eqref{eq:incrementbound2} clearly holds for small $n$ by taking $\alpha$ small enough. Consider $\alpha>0$ such that $\alpha n=2^h$ for some $h\in \N$. By our assumption on $n$ we have $h =  \log_2(\alpha n) \leq \log_2(n)\leq \frac{n^\dxp}{2}$. We define the events
	\begin{align*}
	\Omega_l = \bigcap_{j=0}^{2^{h-l}-1} \left\{ D_{V_{je_1}^{2^l}} \left(j2^le_1,(j2^l+2^l-1)e_1\right) \leq \left( 2 \sum_{i=0}^\infty \left(\frac{3}{2}\right)^{-i\dxp} \right)^{-1} n^{\dxp}\left(\frac{3}{2}\right)^{(l-h)\dxp} \right\}
	\end{align*}
	and
	\begin{align*}
	\Omega = \bigcap_{l=0}^{h} \Omega_l \text .
	\end{align*}
	For an $x \in \{0,\ldots,2^h\}$, say $x=\sum_{l=0}^{h} x_l 2^l$, where $x_l\in\{0,1\}$ for all $l$, we consider the path that goes from $\mz$ to $\left(\sum_{l=h}^{h} x_l 2^l\right)e_1$, from there to $\left(\sum_{l=h-1}^{h} x_l 2^l\right)e_1$, and iteratively to $\left(\sum_{l=0}^{h} x_l 2^l\right)e_1=xe_1$. Using this path from $\mz$ to $xe_1$ through the dyadic points of the form $2^le_1$, one gets that on the event $\Omega$ one has for all $x \in \{0, \ldots , \alpha n\}$
	\begin{align*}
	D_{V_\mz^{x+1}}(\mz,xe_1) \leq \left(2 \sum_{i=0}^\infty \left(\frac{3}{2}\right)^{-i\dxp} \right)^{-1}  n^{\dxp} \sum_{l=0}^h \left(\frac{3}{2}\right)^{(l-h)\dxp} + h <
	2^{-1} n^\dxp + h
	\leq 
	n^\dxp\text,
	\end{align*}
	where we used that $h \leq \frac{n^\dxp}{2}$ in the last step.
	Now, we want to estimate the probability of the event $\Omega$. Let us write $C(\dxp)$ for the constant $\left(2 \sum_{i=0}^\infty \left(\frac{3}{2}\right)^{-i\dxp} \right)^{-1}$ and let $C_{\frac{4}{\dxp}}$ be a constant such that
	\begin{align*}
	\E_\beta \left[ D_{V_\mz^n}(\mz,(n-1)e_1)^{4/\dxp} \right] \leq C_\frac{4}{\dxp} \left(n^\dxp\right)^{4/\dxp} = C_\frac{4}{\dxp} n^4
	\end{align*}
	for all $n\in \N$. Such a constant exists by \Cref{coro:allmoments}. By an application of Markov's inequality we get that
	\begin{align}\label{eq:omegaprobbound}
	& \notag \p_{\beta } \left( D_{V_\mz^{2^l}} \left(\mz,(2^l-1)e_1\right) > C(\dxp) n^{\dxp}\left(\frac{3}{2}\right)^{(l-h)\dxp} \right)\\
	&\notag = \p_{\beta } \left( D_{V_\mz^{2^l}} \left(\mz,(2^l-1)e_1\right)^{\frac{4}{\dxp}} > C(\dxp)^{\frac{4}{\dxp}} n^{\dxp \frac{4}{\dxp}} \left(\frac{3}{2}\right)^{(l-h)\dxp \frac{4}{\dxp} } \right)\\
	&\notag \leq \E_\beta \left[ D_{V_\mz^{2^l}} \left(\mz,(2^l-1)e_1\right)^{\frac{4}{\dxp}} \right]
	C(\dxp)^{-\frac{4}{\dxp}}
	n^{-4}
	\left(\frac{3}{2}\right)^{4(h-l)} 
	\leq 
	C(\dxp)^{-\frac{4}{\dxp}}
	C_{\frac{4}{\dxp}} \left(2^{\l \dxp}\right)^{\frac{4}{\dxp}}
	n^{-4 }
	\left(\frac{3}{2}\right)^{4(h-l) }\\
	&
	\leq
	C(\dxp)^{-\frac{4}{\dxp}}
	C_{\frac{4}{\dxp}} 2^{4l}
	\alpha^4 2^{-4h}
	\left(\frac{3}{2}\right)^{4(h-l) }
	\end{align}
	Define $a_k\coloneqq -2$ and define $a_{k+1}$ as in \eqref{eq:ak definition}. Then one has the line of implications
	\begin{align*}
	\Omega  & \Rightarrow \left\{ D_{V_\mz^{x+1}}(\mz,xe_1) \leq n^\dxp \text{ for all } x \in \{0,\ldots,\alpha n\} \right\} \\
	& \Leftrightarrow \left\{ a_{k+1}(n) \geq \alpha n \right\} \Rightarrow \left\{ \frac{a_{k+1}(n)-a_k(n)}{n} > \alpha \right\} \ .
	\end{align*}
	This already gives us that
	\begin{align*}
	& \p_\beta  \left(\frac{a_{k+1}(n) - a_{k}(n) }{n} \leq \alpha\right) 
	\leq \p_\beta\left(\Omega^C\right)
	\leq \sum_{l=0}^h 2^{h-l} \p_{\beta } \left( D_{V_\mz^{2^l}} \left(\mz,(2^l-1)e_1\right) > C(\dxp) n^{\dxp}\left(\frac{3}{2}\right)^{(l-h)\dxp} \right)\\
	& \overset{\eqref{eq:omegaprobbound}}{\leq}  C(\dxp)^{-\frac{4}{\dxp}}  C_{\frac{4}{\dxp}}
	\sum_{l=0}^h 2^{h-l} 2^{4l}
	\alpha^4 2^{-4h}
	\left(\frac{3}{2}\right)^{4(h-l) }
	= \alpha^4 C(\dxp)^{-\frac{4}{\dxp}}  C_{\frac{4}{\dxp}} \sum_{l=0}^h \left( \frac{81}{128} \right)^{h-l}
	< 0.5
	\end{align*}
	for some $\alpha>0$ small enough. So in particular this implies \eqref{eq:incrementbound2}.
\end{proof}

\begin{lemma}\label{lem:bootstrapping}
	Assume that
	\begin{align}\label{eq:bootstrap assumption}
	\sup_{n \in \N} \E_\beta \left[ \exp \left( \left( \frac{\dia \left(\{0,\ldots,n-1\}^d\right)}{n^\dxp}\right)^{\eta}\right) \right] < \infty
	\end{align}
	for some $\eta>0$. Then
	\begin{align}\label{eq:bootstrap conclusion}
	\sup_{n \in \N} \E_\beta \left[ \exp \left( \left( \frac{\dia \left(\{0,\ldots,n-1\}^d\right)}{n^\dxp}\right)^{\bar{\eta}}\right) \right] < \infty
	\end{align}
	for all $\bar{\eta} \in (0, 1+\dxp \eta)$.
\end{lemma}

\begin{proof}
	Assume that \eqref{eq:bootstrap assumption} holds for some $\eta>0$. Then \Cref{lem:distance to dia} implies that 
	\begin{align}\label{eq:diam stretched exp bd}
	\p_\beta \left( \dia \left(V_\mz^{\bar{n}}\right) > S  C_\dxp n^\dxp \text{ for some } \bar{n} \in \{0,\ldots,n\} \right) 
	\leq C \exp\left(- S^{\eta}\right)
	\end{align}
	for some constants $C,C_\dxp <\infty$. As before, we define $a_k(n)$ inductively by $a_0(n)=0$ and
	\begin{align*}
	a_{k+1}(n) = & (a_k(n) + 2) \\
	& + \sup \left\{z \in \N_{>0} : D_{(a_k+2)e_1+\{0,\ldots,z\}^d}\big((a_k(n)+2)e_1, (a_k(n)+2)e_1 + ze_1\big) \leq n^\dxp \right\}.
	\end{align*}
	The differences $\left(a_{k+1}(n)-a_k(n)\right)_{k\in \N_0}$ are independent and identically distributed. 
	For $\alpha \in (0,1)$, we have that
	\begin{align}\label{eq:bootstrap product}
	\notag & \p_{\beta } \left( \frac{a_1(n)-a_0(n)}{n} \leq \alpha \right)
	=
	\p_{\beta } \left( D_{2e_1 +\{0,\ldots,z\}^d }\left(2e_1, (2+z)e_1\right) > n^\dxp \text{ for some } z \in \{2,\ldots, \lfloor \alpha n \rfloor\}\right)\\
	\notag &
	\leq
	\p_{\beta } \left( \dia \left(\{0,\ldots,z\}^d \right) > n^\dxp \text{ for some } z \in \{0,\ldots, \lfloor \alpha n \rfloor\} \right)\\
	\notag &
	=
	\p_{\beta } \left( \dia \left(\{0,\ldots,z\}^d \right) > \frac{1}{\alpha^\dxp C_\dxp} C_\dxp (\alpha n)^\dxp \text{ for some } z \in \{0,\ldots, \lfloor \alpha n \rfloor\} \right)\\
	&
	\overset{\eqref{eq:diam stretched exp bd}}{\leq} C \exp \left(- \left(\frac{1}{\alpha^\dxp C_\dxp}\right)^\eta \right)
	=
	C \exp \left(- C_\dxp^\prime \alpha^{-\dxp \eta} \right)
	\end{align}
	for a constant $C_\dxp^\prime \in \R_{>0}$. Remember that the random variable $K_n$ was defined by
	\begin{equation*}
	K_n= \inf \left\{k \in \N : a_k(n) \geq n \right\} \ 
	\end{equation*}
	and that
	\begin{align}\label{eq:K n bound}
	D_{V_\mz^n}\left(\mz,(n-1)e_1\right) \leq K_n n^\dxp + 2K_n \leq 3K_n n^\dxp.
	\end{align}
	Assume that $K_n > 2L$ for some large integer $L$. Then there needs to exist at least $L$ indices $i \in \{1,\ldots,2L\}$ such that $a_i(n)-a_{i-1}(n) \leq \frac{1}{L}$. Using independence of the random variables $a_i(n)-a_{i-1}(n)$
	\begin{align*}
	&\p_\beta\left(K_n > 2L\right) 
	\leq 
	\p_\beta\left( \bigcup_{\substack{U \subset \{1,\ldots,2L\}: \\ |U|=L}} \left\{ a_i(n)-a_{i-1}(n) \leq \frac{1}{L} \text{ for all } i\in U \right\} \right)\\
	&
	\leq 
	\sum_{\substack{U \subset \{1,\ldots,2L\}: \\ |U|=L}} \ \prod_{i\in U}
	\p_\beta\left(  a_i(n)-a_{i-1}(n) \leq \frac{1}{L}  \right)
	\leq 
	2^{2L} 
	\p_\beta\left(   a_1(n)-a_{0}(n) \leq \frac{1}{L}  \right)^L\\
	& \overset{\eqref{eq:bootstrap product}}{\leq}
	2^{2L} C \exp \left(-C_\dxp^\prime L^{\dxp \eta} \right)^L 
	\leq \bar{C} \exp\left(-\bar{C_\dxp} L^{1+\dxp \eta}\right)
	\end{align*}
	for some constants $\bar{C},\bar{C_\dxp}\in \R_{>0}$ and all $L$ large enough. From \eqref{eq:K n bound} we have the line of implications
	\begin{align*}
	\left\{ D_{V_\mz^n}\left(\mz,(n-1)e_1\right) > 6L n^\dxp \right\} \Rightarrow \left\{K_n > 2L\right\}
	\end{align*}
	and thus we get that for $L$ large enough
	\begin{align*}
		\p_\beta \left(D_{V_\mz^n}\left(\mz,(n-1)e_1\right) > 6L n^\dxp\right)
		\leq
		\p_\beta \left( K_n > 2L \right)
		\leq
		\bar{C} \exp\left(-\bar{C_\dxp} L^{1+\dxp \eta}\right) \text ,
	\end{align*}
	which implies that
	\begin{align*}
	\sup_{n \in \N} \E_\beta \left[ \exp \left(\left(\frac{D_{V_\mz^n}(\mz,(n-1)e_1)}{n^\dxp}\right)^{\bar{\eta}}\right) \right] < \infty
	\end{align*}
	for all $\bar{\eta} \in (0, 1 + \dxp \eta)$. The same technique as in the proof of \Cref{lem:moments for e1} shows that 
	\begin{align*}
	\sup_{n \in \N} \E_\beta \left[ \exp \left(\left(\frac{D_{V_\mz^n}(\mz,(n-1)e)}{n^\dxp}\right)^{\bar{\eta}}\right) \right] < \infty
	\end{align*}
	for all $e \in \{0,1\}^d$ and all $\bar{\eta} < 1+\dxp \eta$. Using \Cref{lem:distance to dia}, we can finally see that
	 this also implies that
	\begin{align*}
	\sup_{n \in \N} \E_\beta \left[ \exp \left( \left( \frac{\dia \left(\{0,\ldots,n-1\}^d\right)}{n^\dxp}\right)^{\bar{\eta}}\right) \right] < \infty
	\end{align*}
	for all $\bar{\eta} \in (0, 1+\dxp \eta)$.
\end{proof}

With this, we are ready to go  to the proof of \Cref{theo:tail behavior}, which works via a bootstrapping argument.

\begin{proof}[Proof of \Cref{theo:tail behavior}]
	\Cref{lem:moments for e1} and \Cref{lem:distance to dia} imply that \eqref{eq:stretched upper bound} holds for $\bar{\eta}=0.4$.
	We define the function $f(x)=1+\dxp\cdot x$. \Cref{lem:bootstrapping} says that if \eqref{eq:stretched upper bound} holds for some $\bar{\eta}>0$, then it actually holds for all $\eta < f(\bar{\eta})$. Iterating this argument, we see that \eqref{eq:stretched upper bound} holds for all $\eta < f^{(k)}(0.4)$, where $k\in \N$ is an arbitrary integer and $f^{(k)}$ is the $k$-fold iteration of $f$. Letting $k$ go to infinity, the value $f^{(k)}(0.4)$ converges to the fixed point $x_0$ of the equation $x=f(x)$, which is given by $x_0 = \frac{1}{1-\dxp}$. So in particular we see that \eqref{eq:stretched upper bound} holds for all $\eta< \frac{1}{1-\dxp}$.
\end{proof}

\section{Comparison with different inclusion probabilities}\label{section:comparison}

In this section, we compare the graph distances that result from percolation with the measure $\p_\beta$ to the graph distances that result from independent bond percolation on $\Z^d$ where two vertices $u,v\in \Z^d$ are connected with probability $p^\star\left(\beta,\{u,v\}\right)$, which is assumed to be close enough to $p \left(\beta,\{u,v\}\right)$. The precise condition required for the function $p^\star\left(\beta,\{u,v\}\right)$ is that for fixed $\beta$ it satisfies that 
\begin{align}\label{eq:bigO}
	p^\star\left(\beta,\{u,v\}\right)=1 \text{ for } \|u-v\| =1 \text{ and }
	p^\star\left(\beta,\{u,v\}\right)= p \left(\beta,\{u,v\}\right) + \mathcal{O}\left(\frac{1}{\|u-v\|^{2d+1}}\right)
\end{align}
as $\|u-v\|\to \infty$. An example of such a set of inclusion probabilities $p^\star\left(\beta,\{u,v\}\right)$ is given by
\begin{align*}
	p^\star\left(\beta,\{u,v\}\right)=\begin{cases}
	1 & \text{for } \|u-v\| =1\\
	\frac{\beta}{\|u-v\|^{2d}}\wedge 1 & \text{for } \|u-v\| > 1
	\end{cases}
\end{align*}
where we prove in \Cref{example} that \eqref{eq:bigO} is satisfied. These inclusion probabilities were for example also used in \cite{ding2013distances} for $d=1$.

We write $\p^{\star}_\beta$ for the probability measure resulting from independent bond percolation with inclusion probabilities $\left(p^\star(\beta,\{u,v\})\right)_{u,v \in \Z^d}$. In the following, we give a proof that both the graph distance $D\left(\mz,x\right)$ and the diameter of a box $\dia \left(\{0,\ldots,n\}^d\right)$ scale like $\|x\|^{\dxp(\beta)}$, respectively $n^{\dxp(\beta)}$, under the measure $\p_\beta^\star$. 

\begin{theorem}
	For fixed $\beta \geq 0$, suppose that $p^\star\left(\beta,\{u,v\}\right)$ satisfies \eqref{eq:bigO}. Then the graph distance between the origin $\mz$ and $x \in \Z^d$ satisfies
	\begin{align}\label{eq:comparison distance}
		\|x\|^{\dxp(\beta)} \approx_P D\left(\mz,x\right) \approx_P \E_\beta^\star \left[ D\left(\mz,x\right) \right]
	\end{align}
	under the measure $\p_\beta^\star$. The diameter of cubes satisfies
	\begin{align}\label{eq:comparison diameter}
	n^{\dxp(\beta)} \approx_P \dia\left(\{0,\ldots,n\}^d\right) \approx_P \E_\beta^\star \left[ \dia\left(\{0,\ldots,n\}^d\right) \right]
	\end{align}
	under the measure $\p_\beta^\star$.
\end{theorem}

For the proof of \eqref{eq:comparison distance}, we follow a technique that was already used in \cite{ding2013distances} in a similar form for a comparison between the discrete and the continuous model of percolation. The proof of \eqref{eq:comparison diameter} needs more involved methods, and is done in Section \ref{subsec:dia compare}.

\begin{proof}[Proof of \eqref{eq:comparison distance}]
	We fix the dimension $d$ and $\beta$ from here on and consider them as constants. 
	We write $E_{u,v}^\star$ for the event when there exists an edge between $u$ and $v$ in the graph sampled with the measure $\p_\beta^\star$, and we write $E_{u,v}$ if there exists an edge between $u$ and $v$ in the graph sampled with the measure $\p_\beta$. With the standard coupling for percolation we can couple the measures $\p_\beta$ and $\p_\beta^\star$ so that uniformly over all $u\in \Z^d$, $v \in \Z^d \setminus \{u\}$
	\begin{align*}
		\p\left( E_{u,v}^\star \setminus E_{u,v}\right) + \p\left( E_{u,v} \setminus E_{u,v}^\star \right) \leq C_1\frac{1}{\|u-v\|^{2d+1}} 
	\end{align*}
	where $C_1 < \infty$ is a constant, and where we write $\p$ for the joint measure. Thus we also get
	\begin{align*}
		\p\left( \left(E_{u,v}^\star \right)^C \big| E_{u,v}\right) + \p\left( \left(E_{u,v}\right)^C \big| E_{u,v}^\star \right) \leq C_2\frac{1}{\|u-v\|} 
	\end{align*}
	for some constant $C_2 < \infty$. We write $\omega^\star$ for the percolation configuration sampled by $\p_\beta^\star$ and $\omega$ for the percolation configuration sampled by $\p_\beta$. 
	For two points $x,y \in \Z^d$, let $P$ be a geodesic between $x$ and $y$ for the environment $\omega$. We construct a path between $x$ and $y$ in the environment $\omega^\star$ with the following two rules. For each edge $\{u,v\} \in P$ we give a path between $u$ and $v$ in the environment $\omega^\star$. Concatenating all these paths then gives a path between $x$ and $y$ in the environment $\omega^\star$. For the individual sub-paths in $\omega^\star$, we use the following two rules:
	\begin{itemize}
		\item For $\{u,v\} \in P$, if $E_{u,v}^\star$ occurs we use the direct edge between $u$ and $v$ as a path on $\omega^\star$.
		\item For $\{u,v\} \in P$, if $E_{u,v}^\star$ does not occur, we consider the path that connects $u$ to $v$ using $ \|u-v\|_1$ many nearest-neighbor edges.
	\end{itemize}
	This gives a path $P^\star$ between $x$ and $y$ in the environment $\omega^\star$. The length of this path equals
	\begin{align*}
		\sum_{\substack{\{u,v\}\in P : \\ E_{u,v}^\star \text{occurs} }} 1 \ + \sum_{\substack{\{u,v\}\in P : \\ (E_{u,v}^\star)^C \text{occurs} }} \|u-v\|_1 
		= \sum_{\{u,v\}\in P} \left( \mathbbm{1}_{E_{u,v}^\star} + \|u-v\|_1 \mathbbm{1}_{(E_{u,v}^\star)^C}  \right)
	\end{align*}
	and thus we get that
	\begin{align}\label{eq:diameter gleich bound}
		&\notag \E\left[ D(x,y;\omega^\star) \ \big| \ \omega \ \right] \leq 
		\sum_{\{u,v\}\in P} \E \left[  1 + \|u-v\|_1 \mathbbm{1}_{(E_{u,v}^\star)^C} \ \big| \ \omega \ \right]\\
		& \leq \sum_{\{u,v\}\in P} \left(1+ \|u-v\|_1 C_2 \frac{1}{\|u-v\|} \right) \leq C_3 D(x,y;\omega)
	\end{align}
	for some constant $C_3 < \infty$. Markov's inequality for the conditional measure $\p\left( \cdot \big| \omega\right)$ gives that for each $\eps>0$ there exists a constant $C_\eps$ such that
	\begin{align*}
		\p\left( D(x,y;\omega^\star) \leq C_\eps D(x,y;\omega) \right) \geq 1-\eps \ .
	\end{align*}
	Interchanging the roles of $\omega$ and $\omega^\star$ one gets that for each $\eps>0$ there exists a constant $C_\eps^\star$ such that
	\begin{align*}
	\p\left( D(x,y;\omega) \leq C_\eps^\star D(x,y;\omega^\star) \right) \geq 1-\eps \text,
	\end{align*}
	which shows that $D(x,y;\omega^\star) \approx_P \|x-y\|^{\dxp(\beta)}$. Inequality \eqref{eq:diameter gleich bound}, and interchanging the roles of $\omega$ and $\omega^\star$, implies that $\E\left[D(x,y;\omega^\star)\right] $ and $\E\left[D(x,y;\omega)\right]$ are at most a constant factor apart. Thus we get that $\|x-y\|^{\dxp(\beta)} \approx_P D(x,y;\omega^\star) \approx_P \E\left[D(x,y;\omega^\star)\right]$, which finishes the proof.
\end{proof}

\begin{example}\label{example}
	The inclusion probabilities given by
	\begin{align*}
	p^\star\left(\beta,\{u,v\}\right)=\begin{cases}
	1 & \text{for } \|u-v\| = 1 \\
	\frac{\beta}{\|u-v\|^{2d}}\wedge 1 & \text{for } \|u-v\| > 1
	\end{cases}
	\end{align*}
	satisfy \eqref{eq:bigO}.
\end{example}

\begin{proof}
	For all $x \in v+\cC$ and $ y \in u+\cC$, we have by the triangle inequality
	\begin{align*}
	\|u-v\| - \sqrt{d} \leq \|x-y\| \leq \|u-v\| + \sqrt{d} \text, 
	\end{align*} 
	and this already implies that for $\|u-v\|>\sqrt{d}$
	\begin{align*}
	\frac{1}{\left(\|u-v\| + \sqrt{d} \right)^{2d}} \leq \int_{v+\cC}\int_{u+\cC} \frac{1}{\|x-y\|^{2d}} \md y \md x
	\leq \frac{1}{\left(\|u-v\| - \sqrt{d} \right)^{2d}} \text .
	\end{align*}
	With a Taylor expansion we see that
	\begin{align*}
	\frac{1}{\|u-v\| \pm \sqrt{d}}&
	= \frac{1}{\|u-v\|} \frac{1}{1 \pm \frac{\sqrt{d}}{\|u-v\|}}
	= \frac{1}{\|u-v\|} \left(1+\mathcal{O}\left(\frac{1}{\|u-v\|}\right)\right)
	\\
	&
	= \frac{1}{\|u-v\|} +  \mathcal{O}\left(\frac{1}{\|u-v\|^2}\right)
	\end{align*}
	and raising this expression to the $2d$-th power already gives that
	\begin{align}\label{eq:taylor}
	\int_{v+\cC}\int_{u+\cC} \frac{1}{\|x-y\|^{2d}} \md y \md x
	= \frac{1}{\|u-v\|^{2d}} +  \mathcal{O}\left(\frac{1}{\|u-v\|^{2d+1}}\right)
	\end{align}
	for $\|u-v\|\to \infty$. With the Taylor expansion of the exponential function we have $1-e^{-s}=s+\mathcal{O}(s^2)$ for small $s$ and thus by inserting \eqref{eq:taylor} into the definition of $p\left(\beta,\{u,v\}\right)$ we get
	\begin{align}\label{eq:taylor2}
	p\left(\beta,\{u,v\}\right) = 1-e^{-\beta 	\int_{v+\cC}\int_{u+\cC} \frac{1}{\|x-y\|^{2d}} \md y \md x} = \frac{\beta}{\|u-v\|^{2d}} +  \mathcal{O}\left(\frac{1}{\|u-v\|^{2d+1}}\right) 
	\end{align}
	which implies that
	\begin{align*}
	p^\star\left(\beta,\{u,v\}\right) = \frac{\beta}{\|u-v\|^{2d}} \wedge 1 = p\left(\beta,\{u,v\}\right) +
	\mathcal{O}\left(\frac{1}{\|u-v\|^{2d+1}}\right) \text.
	\end{align*}
\end{proof}

\begin{example}
	The inclusion probabilities given by
	\begin{align*}
	\tilde{p}\left(\beta,\{u,v\}\right)=\begin{cases}
	1 & \text{for } \|u-v\|=1\\
	1-e^{-\frac{\beta}{\|u-v\|^{2d}}} & \text{for } \|u-v\| > 1
	\end{cases}.
	\end{align*}
	satisfy \eqref{eq:bigO}.
\end{example}

\begin{proof}
	By a Taylor expansion of the exponential function we get
	\begin{align*}
		1-e^{-\frac{\beta}{\|u-v\|^{2d}}} = \frac{\beta}{\|u-v\|^{2d}} + \mathcal{O} \left(\frac{1}{\|u-v\|^{4d}}\right)
		& = p^\star(\beta,\{u,v\}) + \mathcal{O} \left(\frac{1}{\|u-v\|^{2d+1}}\right) \text ,
	\end{align*}
	where $p^\star(\beta,\{u,v\}) = \frac{\beta}{\|u-v\|^{2d}}  \wedge 1$ is the function from \Cref{example}.
	We already know from Example \ref{example} that $p^\star(\beta,\{u,v\})$ satisfies \eqref{eq:bigO}. Thus we directly get that  $\tilde{p}\left(\beta,\{u,v\}\right)$ also satisfies \eqref{eq:bigO}. 
\end{proof}

\subsection{The diameter of boxes}\label{subsec:dia compare}

Before going to the proof of \eqref{eq:comparison diameter}, we prove a technical lemma that we will use later in this section. It follows directly from the Burkholder-Davis-Gundy-inequality \cite{burkholder1972integral}.

\begin{lemma}\label{lem:q moments}
	Let $X_1,\ldots,X_m$ be independent random variables such that $|\E\left[X_i\right]|\leq C$ for all $i\in \{1,\ldots,m\}$. Then for all $p\geq 2$, there exists a constant $C^\prime=C^\prime(p,C)$ such that
	\begin{equation*}
	\E \left[ \left|\sum_{i=1}^{m} X_i\right|^p \right] 
	\leq
	C^\prime m^{p/2}\max_{i}\E\left[|X_i|^p\right] + C^\prime m^p \text .
	\end{equation*}
\end{lemma}
\begin{proof}
	Define $Y_i = X_i - \E\left[X_i\right]$. We clearly have
	\begin{align}\label{eq:bdglemma1}
	\notag \E \left[ \left|\sum_{i=1}^{m} X_i \right|^p \right]
	& = 
	\E \left[ \left|\sum_{i=1}^{m} Y_i + \sum_{i=1}^{m} \E\left[X_i\right]\right|^p \right]
	\leq
	2^p \E \left[ \left|\sum_{i=1}^{m} Y_i \right|^p \right] + 2^p \E \left[ \left| \sum_{i=1}^{m} \E\left[X_i\right]\right|^p \right]\\
	&
	\leq
	2^p \E \left[ \left|\sum_{i=1}^{m} Y_i \right|^p \right] + 2^p |mC|^p.
	\end{align}
	The process $M_t = \sum_{i=1}^{t} Y_i$ is a martingale and thus we get by the BDG-inequality \cite{burkholder1972integral} that there exists a constant $C_p$ such that 
	\begin{align}\label{eq:bdglemma2}
	\notag & \E \left[ \left|\sum_{i=1}^{m} Y_i \right|^p \right] 
	\leq C_p \E \left[ \left(\sum_{i=1}^m Y_i^2\right)^{p/2} \right]
	= C_p m^{p/2} \left\| \frac{1}{m} \sum_{i=1}^m Y_i^2 \right\|_{p/2}^{p/2}\\
	&
	\leq C_p m^{p/2} \max_{i} \left\| Y_i^2 \right\|_{p/2}^{p/2}
	= C_p m^{p/2} \max_{i} \E \left[|Y_i|^p\right] \text .
	\end{align}
	For $i \in \{1,\ldots,m\}$, we have $\E \left[|Y_i|^p\right] \leq 2^p \E \left[|X_i|^p\right] + 2^p |\E \left[X_i\right]|^p \leq 2^p \E \left[|X_i|^p\right] + 2^p|C|^p$. Combining this with \eqref{eq:bdglemma1} and \eqref{eq:bdglemma2}, we finally get that
	\begin{align*}
	&\E \left[ \left|\sum_{i=1}^{m} X_i \right|^p \right]
	\leq 
	2^p \E \left[ \left|\sum_{i=1}^{m} Y_i \right|^p \right] + 2^p |mC|^p
	\leq
	2^p C_p m^{p/2} \max_{i} \E \left[|Y_i|^p\right] + 2^p |mC|^p\\
	&
	\leq
	2^p C_p m^{p/2}\left( \max_{i} 2^p \E \left[|X_i|^p\right] + 2^p|C|^p\right) + 2^p |mC|^p
	\leq C^\prime  m^{p/2} \max_{i} \E \left[|Y_i|^p\right] + C^\prime m^p
	\end{align*}
	for an appropriate constant $C^\prime$ depending on $p$ and $C$ only.
\end{proof}

Assume that $\left(p^\star(\beta,e)\right)_{e\in E}$ satisfies \eqref{eq:bigO}. From \eqref{eq:comparison distance} and the fact that $\dia(\{0,\cdots,n\}^d)\geq D(\mz, n \mo)$ it directly follows that there exists a constant $c>0$, and for all $\eps >0$ there exists a $c_\eps >0$, such that
\begin{equation*}
	\p^\star_\beta \left( \dia \left(\{0,\ldots,n\}^d\right) > c_\eps n^{\dxp(\beta)} \right) > 1-\eps \text{ and } \E^\star_\beta \left[ \dia \left(\{0,\ldots,n\}^d \right)\right] > c n^{\dxp(\beta)} 
\end{equation*}
for all $n\in \N$. So we are left to show that 
\begin{equation}\label{eq:dia compare upper bound}
\p^\star_\beta \left( \dia \left(\{0,\ldots,n\}^d\right) \leq C_\eps n^{\dxp(\beta)} \right) > 1-\eps \text{ and } \E^\star_\beta \left[ \dia \left(\{0,\ldots,n\}^d \right)\right] \leq  C n^{\dxp(\beta)} 
\end{equation}
uniformly over all $n\in \N$, for appropriate constants $C,C_\eps$. In the following, we will show that
\begin{align}\label{eq:dia compare reduction}
	\p_\beta^\star \left( D_{V_\mz^{\bar{n}}} \left( \mz, (\bar{n}-1)e_1 \right) \leq \bar{S} n^{\dxp(\beta)} \text{ for all } \bar{n} \in \{0,\ldots,n\}\right) \geq 0.25
\end{align}
for some constant $\bar{S}$ and all $n\in \N$. From there one can with the same techniques as in \Cref{lem:moments for e1} and \Cref{lem:distance to dia} show that \eqref{eq:dia compare upper bound} holds. Thus, we will focus on \eqref{eq:dia compare reduction} from here on. We will only do the case where $n=2^k$ for $k \in \N$ large enough. The general case follows with \Cref{propo:scaling}.
We couple the measures $\p_\beta^\star$ and $\p_\beta$, using the standard Harris coupling for percolation. For an edge $e\in E$, we say that it is {\sl non-regular} if $\omega(e)=1$, but $\omega^\star(e)=0$. In words, if the edge is open under the measure $\p_\beta$, but closed under the measure $\p_\beta^\star$. Let  $C_1$ be a constant such that 
\begin{align*}
	\p \left(e \text{ is non-regular } | \ \omega(e)=1\right) \leq \tfrac{C_1}{|e|} \text ,
\end{align*}
where we write $|\{x,y\}|=\|x-y\|_\infty$ for an edge $\{x,y\}$.
Such a constant exists by the assumption \eqref{eq:bigO}. We will always use $C_1$ as this constant in the rest of the chapter. 
The rough strategy of the proof of \eqref{eq:comparison distance} above was to fill in the gaps that occurred through non-regularities using edges in the nearest-neighbor lattice. Such an approach does not work for the diameter. Instead, we fill in the gaps using a third percolation configuration $\omega^{-}$, which is contained in $\omega^\star$. For this, we first choose a list of parameters whose origin will be clear later on. We choose $	q = \frac{4}{3\dxp(\beta)}$, and we choose $\beta_{-} \in \left[0,\beta\right)$, $\eps>0$ such that
\begin{align}\label{eq:assumption on beta-}
	\dxp(\beta_{-})q - 1 + \eps q = 0.5 \text{ and } \frac{2^{3 q \left(\dxp(\beta_{-})-\dxp(\beta)\right)}}{2.2} < \frac{1}{2.1}
\end{align}
which is possible, as the function $\beta\mapsto \dxp(\beta)$ is continuous in $\beta$ \cite{baeumler2022behaviour}. These definitions seem quite arbitrary at the moment, but they are chosen in a way so that the proof works.
The third percolation configuration $\omega^-$ is distributed according to the measure $\p_{\beta_{-}}$. So we can couple the three percolation configurations $\omega, \omega^\star$, and $\omega^-$ using the standard Harris coupling for percolation. We write $\p$ for the joint measure. We have that $p(\beta_{-},e) \leq p^\star(\beta,e)$ for all edges $e$ that are sufficiently long, which follows directly from \eqref{eq:bigO}. In the following, we will even assume that $p(\beta_{-},e) \leq p^\star(\beta,e)$ for all edges $e$. Removing this assumption is relatively easy, as all nearest-neighbor edges are open. This already implies that under this coupling $D(x,y;\omega^\star)\leq D(x,y;\omega^{-})$ for all points $x,y\in \Z^d$. With this, we are ready to go to the proof of \eqref{eq:dia compare reduction}, which already implies \eqref{eq:comparison diameter}.

\begin{proof}[Proof of \eqref{eq:dia compare reduction}]
Define the event $\cA$ by
\begin{align*}
	\cA = & \bigcap_{l=0}^k \bigcap_{a \in V_{\mz}^{2^{k-l}}} \left\{\dia \left(V_a^{2^l}; \omega^{-}\right) \leq 2^{l\dxp(\beta_{-})} 2^{\eps k}\right\} \text .
\end{align*}
For $k$ large enough, we have $\p\left(\cA\right)\geq 0.5$, as we will argue now. Using that 
\begin{equation*}
	\sup_{l\in \N}\E \left[\exp\left(\frac{\dia \left(V_\mz^{2^l}; \omega^{-}\right)}{2^{l\dxp(\beta_{-})}}\right)\right] < \infty 
\end{equation*}
by \Cref{theo:tail behavior}, we get that for some constant $C$
\begin{align*}
	&	\p \left(\cA^C\right) = \p \left( \exists l \in \{0,\ldots,k\}, a \in V_\mz^{2^{k-l}} : \dia \left(V_a^{2^l}; \omega^{-}\right) > 2^{l\dxp(\beta_{-})} 2^{\eps k} \right)
	\\
	&
	\leq \sum_{l=0}^{k} \sum_{a \in V_\mz^{2^{k-l}}}  \p \left( \dia \left(V_a^{2^l}; \omega^{-}\right) > 2^{l\dxp(\beta_{-})} 2^{\eps k} \right)
	\leq \sum_{l=0}^{k} 2^{d(k-l)}  \p \left( \frac{\dia \left(V_\mz^{2^l}; \omega^{-}\right)}{2^{l\dxp(\beta_{-})}} >  2^{\eps k} \right)
	\\
	&
	\leq \sum_{l=0}^{k} 2^{d(k-l)} C e^{-2^{\eps k}}
	=
	C e^{-2^{\eps k}} \sum_{l=0}^{k} 2^{dl}  < 0.5
\end{align*}
for $k$ large enough. Assume that $\cA$ holds, and let $a \in V_\mz^{2^{k-l}}$, $u,v \in V_a^{2^l}$. Assume that $2^{m-1}< \|u-v\|\leq 2^m$. Then $u,v$ are either in the same box $V_w^{2^m}$, or in adjacent boxes $V_{w_1}^{2^m}, V_{w_2}^{2^m}$ with $\|w_1-w_2\|_\infty = 1$. This implies that $D_{V_a^{2^l}} (u,v;\omega^{-}) \leq 2 \cdot 2^{m\dxp(\beta_{-})} 2^{\eps k} + 1 \leq 4 \cdot \|u-v\|^{\dxp(\beta_{-})}2^{\eps k} + 1$. So if the event $\cA$ holds, then for all $u,v \in V_a^{2^l}$
\begin{equation*}
	D_{V_a^{2^l}} (u,v;\omega^{-}) \leq 5 \|u-v\|^{\dxp(\beta_{-})}2^{\eps k}.
\end{equation*}
For $a\in V_\mz^{2^{k-l}}$, let $P$ be a geodesic between $x=2^l a$ and $y=2^l a + (2^l-1)e_1$ in the set $V_a^{2^l}$ for the environment $\omega$. We construct a path between $2^l a$ and $2^l a + (2^l-1)e_1$ in the environment $\omega^\star$ as follows. For each edge $\{u,v\} \in P$, we construct an open path between $u$ and $v$ in the environment $\omega^\star$. Concatenating these open paths then gives an open path between $x$ and $y$ in the environment $\omega^\star$. For the individual sub-paths, we use the following two rules:
\begin{itemize}
	\item For $\{u,v\} \in P$, if $E_{u,v}^\star$ occurs we use the direct edge between $u$ and $v$ as a path between $u$ and $v$ in the environment $\omega^\star$.
	\item For $\{u,v\} \in P$, if $E_{u,v}^\star$ does not occur, we take the shortest path between $u$ to $v$ within the set $V_a^{2^l}$ in the environment $\omega^{-}$ as a sub-path.
\end{itemize}
This gives a path $P^\star$ between $2^l a$ and $2^l a + (2^l-1)e_1$ in the environment $\omega^\star$, as we assumed that all edges contained in $\omega^{-}$ are also contained in $\omega^\star$. The path $P^\star$ is also contained in $V_a^{2^l}$. Write $X_{\{u,v\}}$ for the distance $D_{V_a^{2^l}}(u,v;\omega^{\star})$. The random variable $X_{\{u,v\}}$ is either $1$ or $D_{V_a^{2^l}}(u,v;\omega^{-})$. We define the random variable $X^\prime_{\{u,v\}}$ by
\begin{align*}
X^\prime_{\{u,v\}} = \begin{cases}
1 & \text{ if } X_{\{u,v\}} = 1\\
\min\left(\|u-v\| , 5 \|u-v\|^{\dxp(\beta_{-})} 2^{\eps k}\right) & \text{ else }
\end{cases},
\end{align*}
so in particular we have $X_{\{u,v\}} \leq X^\prime_{\{u,v\}}$ on the event $\cA$, and this already implies that
\begin{equation}
D_{V_a^{2^l}}(x,y;\omega^\star) \leq \sum_{\{u,v\}\in P} X_{\{u,v\}}^\prime \text . 
\end{equation} 
The important thing about the random variables $X^\prime_{e}$ is that they are independent for different edges $e \in P$, as it is independent for different edges whether they are non-regular. 
Next, we want to estimate the first and the $q$-th moment of the random variable $X^\prime_{\{u,v\}}$.
For the expectation we get that
\begin{align*}
\E \left[ X^\prime_{\{u,v\}} \ | \ \omega(\{u,v\})=1, \cA\right] \leq 1 + \|u-v\| \frac{C_1}{\|u-v\|} = 1+ C_1 \text ,
\end{align*}
whereas for the $q$-th moment we see that
\begin{align*}
	& \E \left[ \left(X^\prime_{\{u,v\}} \right)^q \big| \ \omega(\{u,v\})=1, \cA\right] \leq 1 + 5\|u-v\|^{\dxp(\beta_{-})q} 2^{\eps k q} \frac{C_1}{\|u-v\|}\\
	& \leq 1+ 5C_1 \|u-v\|^{\dxp(\beta_{-})q - 1} 2^{\eps k q}
	\leq C_2 2^{k \left( \dxp(\beta_{-})q - 1 + \eps q \right)},
\end{align*}
for some constant $C_2$, as $\dxp(\beta_{-})q> 1$. Using \Cref{lem:q moments}, we see that there exists a constant $C<\infty$ such that
\begin{align*}
	\E & \left[ D_{V_a^{2^l}}(x,y;\omega^\star)^q | \cA,\omega\right]  \leq \E \left[\left(\sum_{\{u,v\}\in P} X_{\{u,v\}}^\prime\right)^q \ \Big| \ \cA,\omega \right]\\
	&
	\leq C D_{V_a^{2^l}}  (x,y;\omega)^{q/2} C_2 2^{k \left( \dxp(\beta_{-})q - 1 + \eps q \right)} + C D_{V_a^{2^l}}(x,y;\omega)^{q}
	\\
	&
	=
	C D_{V_a^{2^l}}(x,y;\omega)^{q/2} C_2 2^{0.5 k } + C D_{V_a^{2^l}}(x,y;\omega)^{q},
\end{align*}
and now taking expectation with respect to $\omega$ yields
\begin{align*}
	& \E  \left[ D_{V_a^{2^l}}(x,y;\omega^\star)^q | \cA\right] \leq 
	\E \left[ C D_{V_a^{2^l}}(x,y;\omega)^{q/2} C_2 2^{k \left( \dxp(\beta_{-})q - 1 + \eps \right)} + C D_{V_a^{2^l}} (x,y;\omega)^{q} | \cA \right]\\
	&
	\leq 
	\tilde{C} \|x-y\|^{\dxp(\beta)q/2}  2^{0.5 k} + \tilde{C} \|x-y\|^{\dxp(\beta)q} 
\end{align*}
for some constant $\tilde{C}$. Here we also used that $\p\left(\cA\right) \geq 0.5$, and thus for all $r>0$ the $r$-th moment of $D_{V_a^{2^l}} (x,y;\omega)$ is of order $\|x-y\|^{r\dxp(\beta)}$, under the measure $\p\left(\cdot | \cA\right)$.
 Assume that $\|x-y\|_\infty = 2^{\gamma k}$ with $\gamma>\frac{3}{4}$. Then we get
\begin{align*}
\E & \left[ D_{V_a^{2^l}}(x,y;\omega^\star)^q | \cA\right] 
\leq 
\tilde{C} \|x-y\|^{\dxp(\beta)q/2}  2^{k \left( \dxp(\beta_{-})q - 1 + \eps \right)} + \tilde{C} \|x-y\|^{\dxp(\beta)q} \\
&
 \leq C^\prime \left(2^{k\left(\frac{\gamma  \dxp(\beta) q}{2} + 0.5\right)}  + 2^{\gamma k \dxp(\beta)q}  \right)
 =
 C^\prime \left(2^{k\left(\frac{\gamma 2}{3} + 0.5\right)}  + 2^{k \frac{\gamma 4}{3}}  \right) 
 \leq 
 C^{\prime \prime}  2^{k \frac{\gamma 4}{3}}  
 \leq 
 C^{\prime \prime \prime}  \|x-y\|^{q \dxp(\beta)}
\end{align*}
for some constants $C^\prime, C^{\prime \prime}, C^{\prime \prime \prime} < \infty$. The second last inequality holds as $\frac{\gamma 2}{3} + 0.5 < \frac{\gamma 4}{3}$ for $\gamma>\frac{3}{4}$. Using Markov's inequality we see that there exists a constant $C<\infty$ such that for all $l>\frac{3}{4}k$, $a \in V_\mz^{2^{k-l}}$, and $S \geq 1$
\begin{align*}
	&\p \left( D_{V_a^{2^l}}\left( 2^l a , 2^l a + (2^l-1)e_1 ; \omega^\star \right) > S 2^{k\dxp(\beta)} 1.1^{(l-k)\dxp(\beta)} \ \big| \ \cA\right)\\
	&
	\leq \p \left( \left(\frac{D_{V_a^{2^l}}\left( 2^l a , 2^l a + (2^l-1)e_1 ; \omega^\star \right)}{2^{\dxp(\beta)l}}\right)^q > S^q \left(\frac{2}{1.1}\right)^{(k-l)4/3}  \ \big| \ \cA\right)
	\\
	&
	\leq C S^{-q} \left(\frac{2}{1.1}\right)^{-(k-l)4/3} \leq C S^{-q} \left(\frac{1}{2.2}\right)^{k-l}.
\end{align*}
On the other hand, for $l\leq \frac{3}{4}k$ we have $l \leq 3 (k-l)$, which implies that
\begin{align*}
	&\p \left( D_{V_a^{2^l}}\left( 2^l a , 2^l a + (2^l-1)e_1 ; \omega^\star \right) > S 2^{k\dxp(\beta)} 1.1^{(l-k)\dxp(\beta)} \ \big| \ \cA\right)\\
	&
	\leq \p \left( D_{V_a^{2^l}}\left( 2^l a , 2^l a + (2^l-1)e_1 ; \omega^{-} \right) > S 2^{k\dxp(\beta)-l\dxp(\beta)} 2^{l\dxp(\beta)} 1.1^{(l-k)\dxp(\beta)} \ \big| \ \cA\right)\\
	&
	= \p \left( \left(\frac{D_{V_a^{2^l}}\left( 2^l a , 2^l a + (2^l-1)e_1 ; \omega^{-} \right)}{2^{l\dxp(\beta_{-})}}\right)^q 
	> 
	S^q \left(\frac{2}{1.1}\right)^{(k-l)\frac{4}{3}} 2^{l(\dxp(\beta)-\dxp(\beta_{-}))q}
	\ \big| \ \cA\right) \\
	&
	\leq C S^{-q} \left(\frac{1}{2.2}\right)^{k-l} 2^{ql(\dxp(\beta_{-})-\dxp(\beta))}
	\leq
	 C S^{-q}   \left(\frac{1}{2.2}\right)^{k-l} 2^{q3(k-l)(\dxp(\beta_{-})-\dxp(\beta))} 
	 \leq C  S^{-q} \left(\frac{1}{2.1}\right)^{k-l}
\end{align*}
where the last inequality holds because of our assumption on $\beta_{-}$ \eqref{eq:assumption on beta-}. So in total we see that there exists a constant $C$ such that for all $k\in \N$, $l \in \{0,\ldots,k\}$, and $a \in V_\mz^{2^{k-l}}$ one has
\begin{align*}
	\p \left( D_{V_a^{2^l}}\left( 2^l a , 2^l a + (2^l-1)e_1 ; \omega^\star \right) > S 2^{k\dxp(\beta)} 1.1^{(l-k)\dxp(\beta)}  \ \big| \ \cA   \right) \leq C S^{-q} \left(\frac{1}{2.1}\right)^{k-l} .
\end{align*}
Write $\Omega^S$ for the event 
\begin{align*}
	\Omega^S = \bigcap_{l=0}^k \bigcap_{j=0}^{2^{k-l}-1} \left\{D_{V_{2^l j e_1}^{2^l}}\left( 2^l j e_1 , 2^l j e_1 + (2^l-1)e_1 ; \omega^\star \right) \leq S 2^{k\dxp(\beta)} 1.1^{(l-k)\dxp(\beta)}   \right\}
\end{align*}
we get with a union bound that
\begin{align*}
	\p \left(\left(\Omega^S\right)^C  \big| \ \cA \right) & \leq \sum_{l=0}^k \sum_{j=0}^{2^{k-l}-1}
	\p \left( D_{V_{2^l j e_1}^{2^l}}\left( 2^l j e_1 , 2^l j e_1 + (2^l-1)e_1 ; \omega^\star \right) > S 2^{k\dxp(\beta)} 1.1^{(l-k)\dxp(\beta)}  \ \Big| \ \cA \right)\\
	&
	\leq 
	\sum_{l=0}^k 2^{k-l}  C S^{-q} \left(\frac{1}{2.1}\right)^{k-l} < 0.5 
\end{align*}
for $S$ large enough. Thus we get that $\p \left(\Omega^S  \right) \geq \p \left(\Omega^S \ \big| \ \cA \right) \p \left( \cA \right) > 0.25$. Using a dyadic path between $\mz$ and $(\bar{n}-1)e_1$, one can show that on the event $\Omega^S$ one has $D_{V_\mz^{\bar{n}}}(\mz,(\bar{n}-1)e_1) \leq C(\dxp(\beta))S n^{\dxp(\beta)}$ for some constant $C(\dxp(\beta))$, depending on $\dxp(\beta)$ only. This shows \eqref{eq:dia compare reduction} and thus finishes the proof.
\end{proof}

\section{Appendix: Proofs for $d=1$}

In this appendix, we show a few lemmas for $d=1$, where slightly different techniques compared to $d\geq 2$ are needed. It is well-known that for fixed $\beta < 1$ one has $\E_\beta \left[D(0,n)\right] = \Omega \left(n^{1-\beta}\right)$. The next lemma gives a more uniform bound on the growth of $\E_\beta \left[D(0,n)\right]$ that holds for all $\beta \in \left[0,1\right]$ simultaneously. 

\begin{lemma}\label{lem:comparethedistances}
	There exists a $c>0$ such that for all $M,n \in \N$ and $\beta \in \left[0,1\right]$
	\begin{align}\label{eq:comparethedistances}
	\E_\beta \left[ D_{\left[0,Mn-1\right]} \left(0,M n-1\right) \right] \geq c M^{1-\beta} \E_\beta \left[ D_{\left[0,n-1\right]} \left(0, n-1\right) \right] \text .
	\end{align}
\end{lemma}

\begin{proof}
	First note that the proof of \eqref{eq:set_to_set expectation} does not depend on a uniform bound on the second moment and works as written above. So we can safely apply it in our argumentation here. By \eqref{eq:set_to_set expectation} we can choose $\iota>0$ small enough so that
	\begin{equation*}
	\E_\beta \left[ D_{\left[0,n-1\right]}\left( \left[0,\iota n\right] , \left[n-\iota n-1, n-1\right] \right) \right] \geq \frac{1}{2} \E_\beta \left[ D_{\left[0,n-1\right]} (0,n-1) \right]
	\end{equation*}
	uniformly over $\beta \in \left[0,1\right]$. This implies the existence of a $c^\prime>0$ such that
	\begin{align}\label{eq:near intervals bound}
	\E_\beta \left[ D_{\left[-n,2n-1\right]} \left(V_{-1}^n, V_1^n\right) \right] \geq c^\prime \E_\beta \left[ D_{\left[0,n-1\right]}(0,n-1)\right]
	\end{align}
	uniformly over $\beta \in \left[0,1\right]$ and $n\in \N$ large enough, as we will argue now. For fixed $\iota > 0$ there is a uniform positive probability (in $\beta \in \left[0,1\right]$ and $n\in \N$) that the rightmost vertex incident to $V_{-1}^n$ lies inside $\left[0,\iota n\right]$ and that the leftmost vertex incident to $V_{1}^n$ lies inside $\left[n-\iota n-1, n-1\right]$. Call this event $A$.  Whenever the event $A$ holds, one already has
	\begin{equation*}
		D_{\left[-n,2n-1\right]} \left(V_{-1}^n, V_1^n\right)
		\geq
		D_{\left[0,n-1\right]}\left( \left[0,\iota n\right] , \left[n-\iota n-1, n-1\right] \right) \text ,
	\end{equation*}
	and as both the event $A$ and the distance $D_{\left[0,n-1\right]}\left( \left[0,\iota n\right] , \left[n-\iota n-1, n-1\right] \right)$ are decreasing one has by the FKG inequality
	\begin{align*}
		&\E_\beta \left[ D_{\left[-n,2n-1\right]} \left(V_{-1}^n, V_1^n\right) \right]  \geq
		\E_\beta \left[ D_{\left[-n,2n-1\right]} \left(V_{-1}^n, V_1^n\right) \mathbbm{1}_A \right] \\
		&\geq
		\E_\beta \left[ D_{\left[0,n-1\right]}\left( \left[0,\iota n\right] , \left[n-\iota n-1, n-1\right] \right)  \mathbbm{1}_A \right] \\
		&
		\geq 
		\E_\beta \left[ D_{\left[0,n-1\right]}\left( \left[0,\iota n\right] , \left[n-\iota n-1, n-1\right] \right)  \right] \p_\beta(A)
		\geq
		\frac{\p_{\beta }(A)}{2} \E_\beta \left[ D_{\left[0,n-1\right]} (0,n-1) \right] \text ,
	\end{align*}
	which shows \eqref{eq:near intervals bound}.
	For long-range percolation on the line segment $\{0,\ldots,M-1\}$, we call an odd point $w \in \{1,\ldots, M-2\}$ a {\sl separation point} if $w \nsim \{0,\ldots,w-2\}, w\nsim \{w+2,\ldots,M-1\}$, and $\{0,\ldots,w-1\} \nsim \{w+1,\ldots, M-1\}$; See Figure \ref{fig:separationvertex} for an illustration. Even points can simply never be separation points with our definition. These three events are independent and we can bound the probability of the first event by
	\begin{align*}
	\p_\beta \left(w \nsim \{0,\ldots,w-2\}\right) \geq e^{-\beta \int_{-\infty}^{0} \int_{1}^{2} \frac{1}{|t-s|^2} \md t \md s } \geq e^{-1}.
	\end{align*}
	The same calculation also works for the second event and shows that $\p_\beta \left(w \nsim \{w+2,\ldots,M-1\}\right) \geq e^{-1}$ for all $\beta \in \left[0,1\right]$.
	The probability of the event $\{0,\ldots,w-1\} \nsim \{w+1,\ldots, M-1\}$ can be bounded from below by
	\begin{align*}
	&\prod_{0\leq u < w} \ \prod_{w < v \leq M-1} e^{-\beta \int_{u}^{u+1} \int_{v}^{v+1} \frac{1}{|x-y|^2} \md x \md y }
	= e^{-\beta \int_{0}^{w} \int_{w+1}^{M} \frac{1}{|x-y|^2} \md x \md y } \\
	& \geq  e^{-\beta \int_{0}^{w} \int_{w+1}^{\infty} \frac{1}{|x-y|^2} \md x \md y } 
	= e^{-\beta \int_0^w \frac{1}{w+1-y} \md y}
	= e^{-\beta\log(w+1)} \geq M^{-\beta}.
	\end{align*}
	uniformly over $\beta \in \left[0,1\right]$. Using the independence of the three relevant events, we get that
	\begin{align*}
	\p_\beta&\left( w \text{ is a separation point}\right) =
	\p_\beta \left( w \nsim \{0,\ldots,w-2\}\right) 
	\cdot
	 \p_\beta \left(  w\nsim \{w+2,\ldots,M-1\}\right) \\
	 &
	 \cdot 
	 \p_\beta\left( \{0,\ldots,w-1\} \nsim \{w+1,\ldots, M-1\}\right) 
	 \geq e^{-2}M^{-\beta} \geq 0.1 M^{-\beta} .
	\end{align*}
	\begin{figure}
		\[\begin{tikzpicture}[scale = 1.5]
		\vertex[draw=none] (-3) at (-3,0) {};
		\vertex[ ] (-2) at (-2,0) {$w$$-$$2$};
		\vertex[ ] (-1) at (-1,0) {$w$$-$$1$};
		\vertex[ ] (0) at ( 0,0) {$\ \  w \ \ $};
		\vertex[ ] (1) at ( 1,0) {$w$+$1$};
		\vertex[ ] (2) at ( 2,0) {$w$+$2$};
		\vertex[draw=none] (3) at ( 3,0) {};
		
		\vertex[draw=none] (-A) at (-2.5,0.6) {};
		\vertex[draw=none] (-B) at (-2,0.9) {};
		\vertex[draw=none] (-C) at (-1.8,1.2) {};
		\vertex[draw=none] (-D) at (-1.6,1.5) {};
		
		\vertex[draw=none] (A) at (2.5,0.6) {};
		\vertex[draw=none] (B) at (2,0.9) {};
		\vertex[draw=none] (C) at (1.8,1.2) {};
		\vertex[draw=none] (D) at (1.6,1.5) {};

		\path[very thick]
		(-3) edge (-2)(-2) edge (-1) (-1) edge (0)(0) edge (1)(1) edge (2) (2) edge (3);
		
		\path[lightgray, thick]
		(-2) edge[bend right = 40] (-A) 
		(-1) edge[bend right = 20] (-B)
		(-1) edge[bend right = 20] (-C)
		(-1) edge[bend right = 20] (-D);
		
		\path[lightgray, thick]
		(2) edge[bend left = 40] (A) 
		(1) edge[bend left = 20] (B)
		(1) edge[bend left = 20] (D);
		
		\end{tikzpicture}\]
		\centering
		\parbox{11cm}{\caption{The vertex $w$ is a separation point if all edges $e$ with $|e|\geq 2$ are either strictly to the left or right of $w$, as above.} \label{fig:separationvertex}}
	\end{figure}
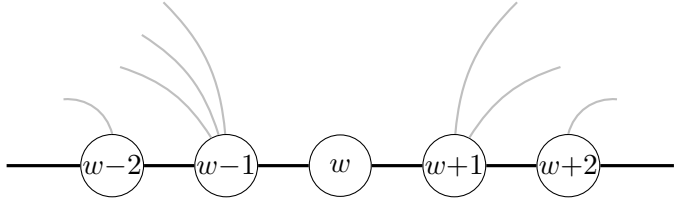
	For odd $w$, we call the set $V_{w}^n$ a {\sl separation  interval} if $V_{w}^n \nsim \left[0,(w-1)n-1\right], V_{w}^n \nsim \left[(w+2)n,Mn-1\right]$, and $\{0,\ldots,wn-1\} \nsim \{(w+1)n,\ldots, Mn - 1\}$. Again, an even $w$ can never define a separation interval. By the scaling invariance of the underlying continuous model, the probability that $V_w^n$ is a separation interval is exactly the probability that  $w$ is a separation point for the line segment $\{0,\ldots, M-1\}$, and this probability is bounded from below by $0.1 M^{-\beta}$. Let $w_1,\ldots,w_l \in \{1,\ldots, M-2\}$ be integers such that $V_{w_i}^n$ is a separation interval for all $i$. Then each path between $0$ and $Mn-1$ in the graph $\left\{0,\ldots,Mn-1\right\}$ needs to cross all separation intervals of this form and in particular
	\begin{align*}
	D_{\left[0,Mn-1\right]}\left(0, Mn-1\right) \geq \sum_{i=1}^l D_{\left[ (w_i-1)n, (w_i+2)n-1 \right]} \left(V_{w_i-1}^n, V_{w_i+1}^n\right) \text .
	\end{align*} 
	The fact that $V_w^n$ is a separation interval reveals no information about the edges with both endpoints in $\left\{(w-1)n,\ldots, (w+2)n-1\right\}$, except that there is no direct edge from $\left\{(w-1)n,\ldots, wn-1\right\}$ to $\left\{(w+1)n,\ldots, (w+2)n-1\right\}$. Thus, by taking expectations in the above inequality and using that both the event $\left\{V_{w}^n \text{ is a sep. int.}\right\}$ and the random distance $D_{\left[ (w-1)n, (w+2)n-1 \right]} \left(V_{w-1}^n, V_{w+1}^n\right)$ are decreasing, we get by the FKG-inequality 
	\begin{align*}
	\E_\beta & \left[ D_{\left[0,Mn-1\right]} \left(0, Mn-1\right)  \right]  \geq \E_\beta \left[ \sum_{w=1}^{M-2} \mathbbm{1}_{\left\{V_{w}^n \text{ is a sep. int.}\right\}}  D_{\left[ (w-1)n, (w+2)n-1 \right]} \left(V_{w-1}^n, V_{w+1}^n\right) \right]\\
	& \geq  \sum_{w=1}^{M-2} \E_\beta \left[ \mathbbm{1}_{\left\{V_{w}^n \text{ is a sep. int.}\right\}} \right]
	\E_\beta \left[ D_{\left[-n,2n-1\right]} \left(V_{-1}^n, V_1^n\right) \right]\\
	& \overset{\eqref{eq:near intervals bound}}{\geq}  \sum_{\substack{w \in \{1,\ldots, M-2\}:\\ w \text{ odd}}} 0.1 M^{-\beta} c^\prime \E_\beta \left[ D_{\left[0,n-1\right]}(0,n-1)\right] \geq c M^{1-\beta } \E_\beta \left[ D_{\left[0,n-1\right]}(0,n-1)\right]
	\end{align*}
	for some small constant $c>0$ and $M$ large enough. For $M$ small, one can take $c$ small enough so that \eqref{eq:comparethedistances} holds for such $M$, by \Cref{propo:scaling}.
\end{proof}

With this we are now ready to go to the proof of Lemma \ref{lem:secondmomentbound} for $d=1$.

\begin{proof}[Proof of Lemma \ref{lem:secondmomentbound} for $d=1$]
	\hypertarget{target:d=1}{} We say that the vertex $w\in\{1,\ldots,m-2\}$ is a cut point (for the interval $\{0,\ldots, m-1\}$) if there exists no edge of the form $\{u,v\}$ with $0\leq u < w <v \leq m-1$. For $w<\frac{m}{2}$ and $\beta \leq 2$ we have
	\begin{align*}
	\p_\beta\left(w \text{ is a cut point}\right) & = \prod_{0\leq u < w} \ \prod_{w < v \leq m-1} e^{-\beta \int_{u}^{u+1} \int_{v}^{v+1} \frac{1}{|x-y|^2} \md x \md y }
	= e^{-\beta \int_{0}^{w} \int_{w+1}^{m} \frac{1}{|x-y|^2} \md x \md y } \\
	& \leq e^{-\beta \int_{0}^{w}  \int_{w+1}^{2w+1} \frac{1}{|x-y|^2} \md x \md y }  
	= e^{-\beta \int_{0}^{w}  \frac{1}{w+1-y} - \frac{1}{2w+1-y} \md y } \\
	& = e^{-\beta \left(-\log(1) +2\log(w+1) - \log(2w+1)\right)}
	= e^{-\beta \log\left(\frac{(w+1)^2}{2w+1}\right)} 
	\leq  e^{-\beta \log\left(\frac{w+1}{2}\right)} \\
	& \leq e^{-\beta \log(w+1)} e^{\beta \log(2)} \leq 4 w^{-\beta}
	\end{align*}
	and with this we get, by linearity of expectation and symmetry of the process, that
	\begin{align*}
	\E_\beta & \left[ \left|\left\{ w \in \{1,\ldots,m-2\} : w \text{ is a cut point} \right\}\right| \right] \leq 1 + 2 \sum_{w=1}^{\lfloor m/2 \rfloor} \p_\beta \left( w \text{ is a cut point} \right) \\
	& \leq 1 + 8 \sum_{w=1}^{\lfloor m/2 \rfloor} w^{-\beta} \leq 10 + 8 \int_{1}^{m} w^{-\beta} \md w = 
	\begin{cases}
	10 + 8\left[\frac{w^{-\beta +1}}{-\beta + 1}\right]_1^m & \beta \in \left[0,2\right]\setminus \{1\} \\
	10 + 8 \log(m) & \beta = 1
	\end{cases}.
	\end{align*}
	As the expected number of cut points is monotone decreasing in $\beta$, we get that for the function $f(\beta, m) \coloneqq \E_\beta  \left[ \left|\left\{ w \in \{1,\ldots,m-2\} : w \text{ is a cut point} \right\}\right| \right]$ we have the upper bound
	\begin{align}\label{eq:fbound}
	f(\beta, m) \leq
	\begin{cases}
	\frac{20}{1-\beta} m^{1-\beta} & \beta < 1 \\
	10 + 8 \log(m) & 1\leq \beta \leq 2 \\
	20 & \beta > 2 
	\end{cases}.
	\end{align}

	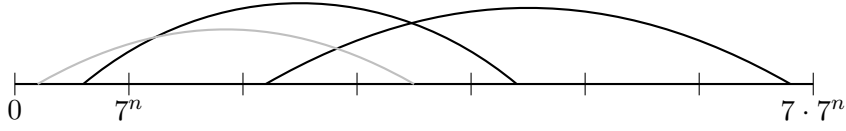
\begin{figure}
		\[\begin{tikzpicture}[scale = 1.5]
		
		\vertex[draw=none] (a) at (0,-0.23) {$0$};
		\vertex[draw=none] (a) at (1,-0.23) {$7^n$};
		\vertex[draw=none] (a) at (7,-0.23) {$7 \cdot 7^n$};
		
		\path[thick]
		(0,0) edge (7,0);
		
		\path
		(0,-0.1) edge (0,0.1)
		(1,-0.1) edge (1,0.1)
		(2,-0.1) edge (2,0.1)
		(3,-0.1) edge (3,0.1)
		(4,-0.1) edge (4,0.1)
		(5,-0.1) edge (5,0.1)
		(6,-0.1) edge (6,0.1)
		(7,-0.1) edge (7,0.1);
		
		\path[thick]
		(0.6,0) edge[bend left = 40] (4.4,0)
		(2.2,0) edge[bend left = 30] (6.8,0);
		
		\path[thick, lightgray]
		(0.2,0) edge[bend left = 30] (3.5,0);
		
		\end{tikzpicture}\]
		\centering
		\parbox{11cm}{\caption{The long edges inside the box $\{0,\ldots,7^{n+1}-1\}$. The set $\mathcal{B}$ are the two bold black edges.} \label{fig:theedgesetB}}
		
	\end{figure}

	We now use a method (that was already used in \cite{ding2013distances} in a similar form for the continuous model) in order to bound the second moment of $D_{V_0^{m^{n+1}}}\left(0,m^{n+1}-1\right)$. We say that an interval $V_k^{m^n}$ is {\sl unbridged} if there exists no edge $\{u,v\}$ with both endpoints in $\left\{ 0,\ldots, m^{n+1}-1 \right\}$ and $u < k m^n, v \geq (k+1) m^n$; Contrary, if there exists such an edge we say that the interval is {\sl bridged}. In this case, we also say that the interval is bridged by the edge $\{u,v\}$. So clearly the intervals $V_0^{m^n}, V_{m-1}^{m^n}$ are unbridged, and the probability that $V_w^{m^n}$ is unbridged for $w \in \left\{1,\ldots,m-2\right\}$ is exactly the probability that $w$ is a cut point for the interval $\{0,\ldots,m-1\}$. We now define a set of edges $\mathcal{B}$ as follows: Let $i<j\in \{0,\ldots,m-1\}$ with $|i-j|>1$ satisfy $V_i^{m^n}\sim V_j^{m^n}$ and $V_{i-l_1}^{m^n} \nsim V_{j+l_2}^{m^n}$ for all $(l_1,l_2) \in \{0,\ldots,i\} \times \{0,\ldots,m-1-j\} \setminus \{(0,0)\}$. In this situation, we add one edge between $V_i^{m^n}$ and $V_j^{m^n}$ to $\mathcal{B}$. If there are several edges between $V_i^{m^n}$ and $V_j^{m^n}$ we choose the left-most shortest such edge (this rule is arbitrary, any deterministic rule would work here). An example of this construction is given in Figure \ref{fig:theedgesetB}. So the set $\mathcal{B}$ is the set of possible bridges where we already delete edges that are furthermore bridged by even longer edges.	 With this construction, we get $|\mathcal{B}|\leq m$, as each interval $V_j^{m^n}$ can be adjacent to at most two edges in $\mathcal{B}$, and each edge in $\mathcal{B}$ touches two intervals. Furthermore, if an interval $V_j^{m^n}$ is bridged, then there exists an edge $e \in \mathcal{B}$ so that $V_j^{m^n}$ is bridged by $e$. Let $\mathcal{U}^\prime$ be the set of endpoints of edges in $\mathcal{B}$ and let 
	\begin{align*}
		\mathcal{U}\coloneqq \mathcal{U}^\prime \cup \left\{0,m^n,\ldots, (m-1)m^n\right\} \cup \left\{m^n-1,2m^n-1,\ldots, m^{n+1}-1\right\}\text.
	\end{align*}
	Let $\mathcal{U}= \left\{x_0,x_1,\ldots, x_u\right\}$, where $x_0 < \ldots < x_u$. By the construction we have $|\mathcal{U}|\leq 4 m$ and $|x_{i-1}-x_i| \leq m^n-1$. For $x_{i-1},x_i$ with $(x_{i-1},x_i) \neq (km^n-1, km^n)$ for all $k$, we say that $\left[x_{i-1},x_{i}\right]$ is bridged, if there exists an edge $\{u,v\}\in \mathcal{B}$ with $u\leq x_{i-1} < x_{i} \leq v$. Assume we have $(x_{i-1},x_i) $ which is not of the form $ (km^n-1, km^n)$, say with $\left[x_{i-1},x_i\right] \subset V_j^{m^n}$ for some $j \in \{0,\ldots,m-1\}$, and  $\left[x_{i-1},x_{i}\right]$ is not bridged. Then also $V_j^{m^n}$ is not bridged. On the other hand, if $\left[x_{i-1},x_{i}\right]$ is bridged, then also $V_j^{m^n}$ is bridged by some edge in $\mathcal{B}$.	
	In each interval $V_j^{m^n}$ there are at most two points in $V_j^{m^n}\cap \mathcal{U}$ that come from endpoints of edges in $\mathcal{B}$; Furthermore, the two endpoints of the interval are also in $V_j^{m^n}\cap \mathcal{U}$. So in total there are at most $4$ points in $V_j^{m^n}\cap \mathcal{U}$ for all $j\in \{0,\ldots,m-1\}$, and thus there are at most three intervals of the form $\left[x_{i-1},x_{i}\right]$ inside each $V_j^{m^n}$. This already implies that
	\begin{align}\label{eq:unbridgednumbercompare}
	\left| \left\{ i \in \{1,\ldots,u\} : \left[x_{i-1},x_i\right] \text{ is not bridged}  \right\} \right| \leq 3
	\left| \left\{ j \in \{0,\ldots,m-1\} : V_j^{m^n} \text{ is not bridged}  \right\} \right| \text .
	\end{align}
	Now we want to construct a path between $0$ and $m^{n+1}-1$. Let
	\begin{equation*}
	\tau = \arg\max_{i \in \{1,\ldots, u\} } D_{\left[x_{i-1},x_i\right]} \left(x_{i-1},x_i\right) \text .
	\end{equation*}
	If there are multiple maximizers, we pick one with $x_i \neq k m^n$ for all $k$, and with minimal $x_i$ among those maximizers. So in particular $\left[x_{i-1},x_i\right]$ always lies inside some interval $V_j^{m^n}$. If $\left[x_{\tau-1},x_\tau\right]$ is bridged by some edge $e= \{x_{\tau_1}, x_{\tau_2}\} \in \mathcal{B}$, say with $x_{\tau_1} < x_{\tau_2}$, then we consider the path that goes from $0=x_0$ to $x_{\tau_1}$, then directly jumps to $x_{\tau_2}$ and from there goes to $x_u=m^{n+1}-1$. This implies that
	\begin{align*}
	D_{[0,m^{n+1}-1]}\left(0,m^{n+1}-1\right) & \leq \sum_{i=1}^{\tau_1} D_{\left[x_{i-1},x_i\right]} \left(x_{i-1},x_i \right) + 1 + \sum_{i=\tau_2 + 1}^{u} D_{\left[x_{i-1},x_i\right]} \left(x_{i-1},x_i \right) \\
	& \leq u \max_{i \neq \tau} D_{\left[x_{i-1},x_i\right]} \left(x_{i-1},x_i \right)
	\leq 4m \max_{i \neq \tau} D_{\left[x_{i-1},x_i\right]} \left(x_{i-1},x_i \right)
	\end{align*}
	in this case. For the case where $\left[x_{\tau-1},x_\tau\right]$ is not bridged, we consider the path that goes iteratively from $x_0$ to $x_u$. Here we have
	\begin{align*}
	 &D_{[0,m^{n+1}-1]}\left(0,m^{n+1}-1\right)
	\\
	\notag &
	\leq \sum_{i=1}^{\tau-1} D_{\left[x_{i-1},x_i\right]} \left(x_{i-1},x_i \right) + D_{\left[x_{\tau-1},x_\tau\right]} \left(x_{\tau-1},x_\tau\right) + \sum_{i=\tau + 1}^{u} D_{\left[x_{i-1},x_i\right]} \left(x_{i-1},x_i \right) \\
	& \leq 4m \max_{i \neq \tau} D_{\left[x_{i-1},x_i\right]} \left(x_{i-1},x_i \right)
	+ \  \max_{\left[x_{i-1},x_i\right] \text{ not bridged} } D_{\left[x_{i-1},x_i\right]} \left(x_{i-1},x_i \right) \text ,
	\end{align*}
	and thus we have in both cases that
	\begin{align}\label{eq:second part}
	\notag &\left(D_{[0,m^{n+1}-1]}\left(0,m^{n+1}-1\right)\right)^2
	\\
	\notag &
	\leq  2  \left(4m \max_{i \neq \tau} D_{\left[x_{i-1},x_i\right]} \left(x_{i-1},x_i \right)\right)^2
	+ 2 \left( \max_{\left[x_{i-1},x_i\right] \text{ not bridged} } D_{\left[x_{i-1},x_i\right]} \left(x_{i-1},x_i \right)\right)^2\\
	& \leq  32m^2 \left( \max_{i \neq \tau} D_{\left[x_{i-1},x_i\right]} \left(x_{i-1},x_i \right)\right)^2
	+ 2 \sum_{\left[x_{i-1},x_i\right] \text{ not bridged} } \left(  D_{\left[x_{i-1},x_i\right]} \left(x_{i-1},x_i \right)\right)^2 .
	\end{align}
	Next, we want to bound both terms in the above sum in expectation. To bound the first term, we use the following observation: If $X_1,\ldots,X_{\tilde{m}}$ are independent non-negative random variables and $\tau = \arg\max_{i \in \{1,\ldots,\tilde{m}\}}$, then
	\begin{align}\label{eq:2ndmax 2nd moment}
	\E\left[\left(\max_{i \neq \tau} X_i\right)^2\right] & \leq \E \left[ \sum_{i=1}^{\tilde{m}}  X_i \left(\sum_{j\neq i} X_j\right) \right] = \sum_{i=1}^{\tilde{m}} \sum_{j\neq i} \E \left[ X_i \right] \E\left[ X_j \right] \leq \tilde{m}^2 \max_{i} \E\left[X_i\right]^2  \text .
	\end{align}
	Conditioned on $\mathcal{U}$, the random variables $\left(D_{\left[x_{i-1},x_i\right]} \left(x_{i-1},x_i\right)\right)_{i=1}^u$ are independent and by Proposition \ref{propo:scaling} their expectation is bounded by the expectation of $\left(D_{V_0^{m^n}}(0,m^n-1)\right)$, up to a factor of $3$. As $u\leq 4m$, we get with \eqref{eq:2ndmax 2nd moment} and Proposition \ref{propo:scaling} that
	\begin{align}\label{eq:expectation bound1}
	\notag \E_\beta & \left[  \max_{i \neq \tau} D_{\left[x_{i-1},x_i\right]} \left(x_{i-1},x_i \right)^2 \right]
	= \E_\beta \left[ \E\left[ \max_{i \neq \tau} D_{\left[x_{i-1},x_i\right]} \left(x_{i-1},x_i \right)^2 \big| \ \mathcal{U} \right] \right] \\
	&
	\leq \E_\beta \left[ 16 m^2  \max_{i} \E\left[  D_{\left[x_{i-1},x_i\right]} \left(x_{i-1},x_i \right) \big| \ \mathcal{U} \right]^2 \right] 
	\leq 144 m^2 \E_\beta \left[D_{V_0^{m^n}} \left(0,m^n-1\right)\right]^2.
	\end{align}
	In order to bound the second summand in \eqref{eq:second part} in expectation, we use the bound on the number of unbridged segments \eqref{eq:unbridgednumbercompare}. Also note that the second moment of $D_{\left[x_{i-1},x_i\right]} \left(x_{i-1},x_i\right)$ is, by Proposition \ref{propo:scaling}, bounded by the second moment of $\left(D_{V_0^{m^n}}(0,m^n-1)\right)$, up to a factor of $9=3^2$.
	With this we get that
	\begin{align}\label{eq:expectation bound2}
	\notag \E_\beta & \left[\sum_{\left[x_{i-1},x_i\right] \text{ not bridged} } \left(  D_{\left[x_{i-1},x_i\right]} \left(x_{i-1},x_i \right)\right)^2 \right] \\
	\notag &
	= \E_\beta \left[ \E_\beta  \left[\sum_{\left[x_{i-1},x_i\right] \text{ not bridged} } \left(  D_{\left[x_{i-1},x_i\right]} \left(x_{i-1},x_i \right)\right)^2 \Big| \ \mathcal{U} \right] \right]\\
	\notag & \leq 9 \E_\beta \left[D_{V_0^{m^n}} \left(0, m^n -1\right)^2 \right] \E_\beta \left[\sum_{\left[x_{i-1},x_i\right] \text{ not bridged} } 1 \right]
	\\ 
	\notag & \leq 27 \E_\beta \left[D_{V_0^{m^n}} \left(0, m^n -1\right)^2 \right] \E_\beta \left[ \left|\left\{ j \in \{0,\ldots,m-1\} : V_j^{m^n} \text{ unbridged}\right\} \right| \right]\\
	& \leq 27 \E_\beta \left[D_{V_0^{m^n}} \left(0, m^n -1\right)^2 \right] \left(2 + f(\beta,m)\right) \eqqcolon   \E_\beta \left[D\left(0, m^n -1\right)^2 \right]  \tilde{f}(\beta,m) \text ,
	\end{align}
	where $\tilde{f}(\beta,m) = 27(2+f(\beta,m))$. Combining \eqref{eq:expectation bound1} and \eqref{eq:expectation bound2}, and taking expectations in \eqref{eq:second part}, we obtain that
	\begin{align*}
	\E_\beta \left[ D_{V_0^{m^{n+1}}}\left(0, m^{n+1}-1 \right)^2 \right] & \leq  5000 m^4 \E_\beta \left[D_{V_0^{m^n}} \left(0,m^n-1\right)\right]^2 \\
	& \ \ \ + 2 \tilde{f}(\beta,m) \E_\beta \left[D_{V_0^{m^n}} \left(0, m^n -1\right)^2 \right] \text .
	\end{align*}
	Iterating this inequality over all $k=1,\ldots,n$, we get
	\begin{align}\label{eq:iterated moment bound}
	\E_\beta \left[ D_{V_0^{m^{n+1}}} \left(0, m^{n+1}-1 \right)^2 \right] 
	& 
	\leq  5000 m^4 \sum_{k=1}^{n} \left(2 \tilde{f} (\beta,m)\right)^{n+1-k} \E_\beta \left[D_{V_0^{m^k}} \left(0, m^k -1 \right)\right]^2  .
	\end{align}
	Using the bounds on $f(\beta,m)$ from \eqref{eq:fbound}, we see that function $\tilde{f}(\beta,m)$ satisfies
	\begin{align}\label{eq:f tilde bound}
	\tilde{f}(\beta, m) = 27 \left(2+ f(\beta,m)\right) \leq
	\begin{cases}
	\frac{600}{1-\beta} m^{1-\beta} & \beta < 1 \\
	600 \left(1+ \log(m)\right) & 1\leq \beta \leq 2 \\
	600 & \beta > 2 
	\end{cases} \ .
	\end{align}
	By compactness of each interval $\left[\beta, \beta + 1\right]$, it suffices to show that the uniform bound on the second moment \eqref{eq:second moment bound} holds for all $\beta>0$ and $\eps \in \left(-c_\beta, c_\beta\right)$ for some $c_\beta>0$ small enough, respectively for $\beta = 0$ and $\eps \in \left[ 0 , c_\beta \right)$.
	To extend the inequality from open sets to compact intervals, one can cover each compact interval with finitely many open sets and then take the largest among these finitely many constants that arose from this procedure. 
	 So we are left to show that
	for all $\beta \geq 0$, there exist a constants $c_\beta >0$ and $C_\beta<\infty$ such that for all $n \in \N$, all $\eps \in (-c_\beta,c_\beta)$, respectively all $\eps \in \left[0,c_\beta\right)$ for $\beta=0$, and all $u,v \in V_{0}^n$
	\begin{align}\label{eq:second moment bound2}
	\E_{\beta+\eps} \left[ D_{V_{0}^n} (u,v)^2 \right] \leq C_\beta \Lambda(n,\beta+\eps)^2\text.
	\end{align}
	We start with the case $\beta \geq 1$. By \Cref{remark:Lambda grows}, there exists a $\dxp^\prime=\dxp^\prime(\beta) > 0 $ such that 
	\begin{align}\label{eq:theta prime}
		\E_{\beta + \eps} \left[D_{V_0^{m^{k+1}}}  \left(0,m\cdot m^k-1\right)\right] \geq m^{\dxp^\prime(\beta) } \E_{\beta + \eps} \left[D_{V_0^{m^{k}}}  \left(0, m^k-1\right)\right]
	\end{align}
	for all $k \in \N$, $m$ large enough, and $|\eps| \leq \frac{1}{2}$. Inserting this into \eqref{eq:iterated moment bound}, we get 
	\begin{align*}
	& \E_\beta  \left[ D_{V_0^{m^{n+1}}} \left(0, m^{n+1}-1 \right)^2 \right] 
	 \\
	&
	\leq  5000 m^4 \sum_{k=1}^{n} \left(2 \tilde{f} (\beta,m)\right)^{n+1-k} \left( m^{-2\dxp^\prime}\right)^{n-k} \E_\beta \left[D_{V_0^{m^n}} \left(0, m^n -1 \right)\right]^2  .
	\end{align*}
	Now choose $m \in \N$ large enough and $c_\beta \in (0,0.1)$ small enough so that $2 \tilde{f}(\beta + \eps, m)m^{-2\dxp^\prime(\beta)} \leq 0.5$ for all $\eps \in \left(-c_\beta, c_\beta \right)$. This is clearly possible for $\beta > 1$. For $\beta = 1$, we can choose $c_\beta$ small enough so that $c_\beta < \dxp^\prime (1)$, where $\theta^\prime(1)$ is the one defined in \eqref{eq:theta prime}. By monotonicity in the first argument of the function $\tilde{f}(\cdot,\cdot)$ one then has $\tilde{f}(1+\eps,m) \leq \frac{600}{c_\beta}m^{c_\beta}$ for all $\eps \in \left(-c_\beta, c_\beta \right)$, which shows that one can find $m,c_\beta$ so that  $2 \tilde{f}(1 + \eps, m)m^{-2\dxp^\prime(1)} \leq 0.5$ for all $\eps \in \left(-c_\beta, c_\beta \right)$. This then gives that
	\begin{align*}
	&\E_{\beta+\eps} \left[ D_{V_0^{m^{n+1}}} \left(0, m^{n+1}-1 \right)^2 \right]  
	\leq  10000 \tilde{f}(\beta-c_\beta, m) m^4 \sum_{k=1}^{n} 0.5^{n-k} \ \E_{\beta + \eps } \left[D_{V_0^{m^{n}}} \left(0, m^n -1 \right)\right]^2 \\
	& \leq 20000 \tilde{f}(\beta-c_\beta, m) m^4  \E_{\beta + \eps } \left[D_{V_0^{m^{n}}} \left(0, m^n -1 \right)\right]^2
	\leq 20000 \tilde{f}(\beta-c_\beta, m) m^4  \Lambda\left(m^n, \beta+\eps\right)^2
	\end{align*}
	for all $\eps \in \left(-c_\beta, c_\beta\right)$. This shows \eqref{eq:second moment bound2} along the subsequence $m,m^2,m^3,\ldots$ To extend inequality \eqref{eq:second moment bound2} from this subsequence to all integers, use \Cref{propo:scaling}.
	
	 Next, we consider the case where $\beta \in (0,1)$. Using Lemma \ref{lem:comparethedistances}, we know that there is a constant $c\in (0,1)$ such that
	\begin{align*}
		\E_{\beta } \left[ D_{V_0^{m^{n}}} \left(0,m^{n}-1\right) \right] 
		& \geq
		c m^{(n-k)(1-\beta)} \E_{\beta } \left[ D_{V_0^{m^{n-k}}} \left(0,m^{n-k}-1\right) \right]\\
		&
		\geq
		\left(c m^{1-\beta}\right)^{n-k} \E_{\beta } \left[ D_{V_0^{m^{n-k}}} \left(0,m^{n-k}-1\right) \right]
	\end{align*}
	for all $n\geq k$ and $m\in \N$. Now take $m$ large enough and $c_\beta$ small enough so that $\frac{1200 \ m^{\beta+\eps-1} }{c (1-\beta-\eps)} < 0.5$ for all $\eps \in \left(-c_\beta,c_\beta\right)$, and that $0<\beta-c_\beta < \beta +c_\beta < 1$. Using \eqref{eq:iterated moment bound} we get that for such $m$ and $\eps \in (-c_\beta, c_\beta)$
	\begin{align*}
	& \E_{\beta + \eps}  \left[ D_{V_0^{m^{n+1}}}\left(0, m^{n+1}-1 \right)^2 \right] 
	\leq  5000 m^4 \sum_{k=1}^{n} \left(2 \tilde{f} (\beta,m)\right)^{n+1-k} \E_\beta \left[D_{V_0^{m^{k}}} \left(0, m^k -1 \right)\right]^2\\
	&
	\leq  5000 m^4 \sum_{k=1}^{n} \left(2 \tilde{f} (\beta+\eps,m)\right)^{n+1-k} \left(cm^{1-\beta-\eps}\right)^{2(k-n)} \E_{\beta + \eps} \left[D_{V_0^{m^{n}}} \left(0, m^n -1 \right)\right]^2 \\
	& \overset{\eqref{eq:f tilde bound}}{\leq} 10000 \tilde{f}(\beta-c_\beta,m)m^4 \sum_{k=1}^{n}  \left( \frac{1200 \ m^{\beta+\eps-1} }{c (1-\beta-\eps)} \right)^{n-k} \E_{\beta + \eps} \left[D_{V_0^{m^{n}}} \left(0, m^n -1 \right)\right]^2 \\
	& \leq 10000 \tilde{f}(\beta-c_\beta,m)m^4 \sum_{k=1}^{n}  0.5^{n-k} \E_{\beta + \eps} \left[D_{V_0^{m^{n}}} \left(0, m^n -1 \right)\right]^2 \\
	&
	\leq 20000 \tilde{f}(\beta-c_\beta,m)m^4  \E_{\beta + \eps} \left[D_{V_0^{m^{n}}} \left(0, m^n -1 \right)\right]^2
	\leq 10^6 m^5 \E_{\beta + \eps} \left[D_{V_0^{m^{n}}} \left(0, m^n -1 \right)\right]^2,
	\end{align*}
	which shows \eqref{eq:second moment bound2} for numbers of the form $m,m^2,m^3,\ldots$ Here, we used that $f(\beta,m) \leq m-2$ for all $\beta \in \R_{\geq 0}$, and thus $\tilde{f}(\beta,m) = 27(2+f(\beta,m)) \leq 27 m$ for the last inequality. To extend inequality \eqref{eq:second moment bound2} from this subsequence to all integers, use \Cref{propo:scaling}. The proof for $\beta=0$ works analogous to the case $\beta \in (0,1)$, and we omit it.
\end{proof}

\begin{proof}[Proof of \Cref{coro:allmoments} for $d=1$]
	We use the same notation as in the proof of \Cref{lem:secondmomentbound} for $d=1$ above. We have that
	\begin{align*}
	D_{[0,m^{n+1}-1]}\left(0,m^{n+1}-1\right)
	 \leq 4m \max_{i \neq \tau} D_{\left[x_{i-1},x_i\right]} \left(x_{i-1},x_i \right)
	+ \  \max_{\left[x_{i-1},x_i\right] \text{ not bridged} } D_{\left[x_{i-1},x_i\right]} \left(x_{i-1},x_i \right)
	\end{align*}
	and this implies that for any $r > 0$
	\begin{align*}
	&\left(D_{[0,m^{n+1}-1]}\left(0,m^{n+1}-1\right)\right)^{2r}\\
	&
	 \leq  2^{2r} 4^{2r} m^{2r} \left( \max_{i \neq \tau} D_{\left[x_{i-1},x_i\right]} \left(x_{i-1},x_i \right)\right)^{2r}
	+ 2^r \sum_{\left[x_{i-1},x_i\right] \text{ not bridged} } \left(  D_{\left[x_{i-1},x_i\right]} \left(x_{i-1},x_i \right)\right)^{2r} \text .
	\end{align*}
	Taking expectations and the same arguments as in the proof of \Cref{lem:secondmomentbound} yield
	\begin{align*}
	& \E_\beta \left[\left( \max_{i \neq \tau} D_{\left[x_{i-1},x_i\right]} \left(x_{i-1},x_i \right)^r\right)^2\right]
	\leq
	\E_\beta \left[16m^2 \max_{i} \E_\beta  \left[  D_{\left[x_{i-1},x_i\right]} \left(x_{i-1},x_i \right)^r \big| \mathcal{U}\right]^2 \right] \\
	&
	\leq
	16 m^2 3^{2r}  \E_\beta \left[ D_{\left[0,m^n-1\right]} \left(0,m^n-1 \right)^r  \right] \text .
	\end{align*}
	From here, the same proof as before gives that $\E_\beta\left[D_{\left[0,n\right]}(0,n)^r\right] \leq C(r) \E_\beta\left[D_{\left[0,n\right]}(0,n)\right]^r$
	for a constant $C(r)$, and $r$ of the form $r=2^k$ with natural $k$. Extending this to all $r>0$ works with H\"older's inequality.
\end{proof}

\end{document}